\documentclass[12pt]{amsart}
\usepackage{amssymb}
\begin{document}
\theoremstyle{plain}
\newtheorem{thm}{Theorem}[section]
\newtheorem*{thm1}{Theorem 1}
\newtheorem*{thm2}{Theorem 2}
\newtheorem{lemma}[thm]{Lemma}
\newtheorem{lem}[thm]{Lemma}
\newtheorem{cor}[thm]{Corollary}
\newtheorem{prop}[thm]{Proposition}
\newtheorem{propose}[thm]{Proposition}
\newtheorem{variant}[thm]{Variant}
\theoremstyle{definition}
\newtheorem{notations}[thm]{Notations}
\newtheorem{rem}[thm]{Remark}
\newtheorem{rmk}[thm]{Remark}
\newtheorem{rmks}[thm]{Remarks}
\newtheorem{defn}[thm]{Definition}
\newtheorem{ex}[thm]{Example}
\newtheorem{claim}[thm]{Claim}
\newtheorem{ass}[thm]{Assumption}
\numberwithin{equation}{section}
\newcounter{elno}                
\def\points{\list
{\hss\llap{\upshape{(\roman{elno})}}}{\usecounter{elno}}} 
\let\endpoints=\endlist


\catcode`\@=11
%
%
\def\opn#1#2{\def#1{\mathop{\kern0pt\fam0#2}\nolimits}} 
\def\bold#1{{\bf #1}}%
\def\underrightarrow{\mathpalette\underrightarrow@}
\def\underrightarrow@#1#2{\vtop{\ialign{$##$\cr
 \hfil#1#2\hfil\cr\noalign{\nointerlineskip}%
 #1{-}\mkern-6mu\cleaders\hbox{$#1\mkern-2mu{-}\mkern-2mu$}\hfill
 \mkern-6mu{\to}\cr}}}
\let\underarrow\underrightarrow
\def\underleftarrow{\mathpalette\underleftarrow@}
\def\underleftarrow@#1#2{\vtop{\ialign{$##$\cr
 \hfil#1#2\hfil\cr\noalign{\nointerlineskip}#1{\leftarrow}\mkern-6mu
 \cleaders\hbox{$#1\mkern-2mu{-}\mkern-2mu$}\hfill
 \mkern-6mu{-}\cr}}}
%
%

%
\def\:{\colon}
\let\oldtilde=\tilde
\def\tilde#1{\mathchoice{\widetilde{#1}}{\widetilde{#1}}%
{\indextil{#1}}{\oldtilde{#1}}}
\def\indextil#1{\lower2pt\hbox{$\textstyle{\oldtilde{\raise2pt%
\hbox{$\scriptstyle{#1}$}}}$}}
\def\pnt{{\raise1.1pt\hbox{$\textstyle.$}}}
%

%
\let\amp@rs@nd@\relax
\newdimen\ex@
\ex@.2326ex
\newdimen\bigaw@
\newdimen\minaw@
\minaw@16.08739\ex@
\newdimen\minCDaw@
\minCDaw@2.5pc
\newif\ifCD@
\def\minCDarrowwidth#1{\minCDaw@#1}
\newenvironment{CD}{\@CD}{\@endCD}
\def\@CD{\def\A##1A##2A{\llap{$\vcenter{\hbox
 {$\scriptstyle##1$}}$}\Big\uparrow\rlap{$\vcenter{\hbox{%
$\scriptstyle##2$}}$}&&}%
\def\V##1V##2V{\llap{$\vcenter{\hbox
 {$\scriptstyle##1$}}$}\Big\downarrow\rlap{$\vcenter{\hbox{%
$\scriptstyle##2$}}$}&&}%
\def\={&\hskip.5em\mathrel
 {\vbox{\hrule width\minCDaw@\vskip3\ex@\hrule width
 \minCDaw@}}\hskip.5em&}%
\def\verteq{\Big\Vert&&}%
\def\noarr{&&}%
\def\vspace##1{\noalign{\vskip##1\relax}}\relax\let\amp@rs@nd@&\iffalse}\fi
 \CD@true\vcenter\bgroup\relax\let\\=\cr\iffalse}\fi\tabskip\z@skip\baselineskip20\ex@
 \lineskip3\ex@\lineskiplimit3\ex@\halign\bgroup
 &\hfill$\m@th##$\hfill\cr}
\def\@endCD{\cr\egroup\egroup}
%
\def\>#1>#2>{\amp@rs@nd@\setbox\z@\hbox{$\scriptstyle
 \;{#1}\;\;$}\setbox\@ne\hbox{$\scriptstyle\;{#2}\;\;$}\setbox\tw@
 \hbox{$#2$}\ifCD@
 \global\bigaw@\minCDaw@\else\global\bigaw@\minaw@\fi
 \ifdim\wd\z@>\bigaw@\global\bigaw@\wd\z@\fi
 \ifdim\wd\@ne>\bigaw@\global\bigaw@\wd\@ne\fi
 \ifCD@\hskip.5em\fi
 \ifdim\wd\tw@>\z@
 \mathrel{\mathop{\hbox to\bigaw@{\rightarrowfill}}\limits^{#1}_{#2}}\else
 \mathrel{\mathop{\hbox to\bigaw@{\rightarrowfill}}\limits^{#1}}\fi
 \ifCD@\hskip.5em\fi\amp@rs@nd@}
\def\<#1<#2<{\amp@rs@nd@\setbox\z@\hbox{$\scriptstyle
 \;\;{#1}\;$}\setbox\@ne\hbox{$\scriptstyle\;\;{#2}\;$}\setbox\tw@
 \hbox{$#2$}\ifCD@
 \global\bigaw@\minCDaw@\else\global\bigaw@\minaw@\fi
 \ifdim\wd\z@>\bigaw@\global\bigaw@\wd\z@\fi
 \ifdim\wd\@ne>\bigaw@\global\bigaw@\wd\@ne\fi
 \ifCD@\hskip.5em\fi
 \ifdim\wd\tw@>\z@
 \mathrel{\mathop{\hbox to\bigaw@{\leftarrowfill}}\limits^{#1}_{#2}}\else
 \mathrel{\mathop{\hbox to\bigaw@{\leftarrowfill}}\limits^{#1}}\fi
 \ifCD@\hskip.5em\fi\amp@rs@nd@}
%
%
\newenvironment{CDS}{\@CDS}{\@endCDS}
\def\@CDS{\def\A##1A##2A{\llap{$\vcenter{\hbox
 {$\scriptstyle##1$}}$}\Big\uparrow\rlap{$\vcenter{\hbox{%
$\scriptstyle##2$}}$}&}%
\def\V##1V##2V{\llap{$\vcenter{\hbox
 {$\scriptstyle##1$}}$}\Big\downarrow\rlap{$\vcenter{\hbox{%
$\scriptstyle##2$}}$}&}%
\def\={&\hskip.5em\mathrel
 {\vbox{\hrule width\minCDaw@\vskip3\ex@\hrule width
 \minCDaw@}}\hskip.5em&}
\def\verteq{\Big\Vert&}
\def\novarr{&}
\def\noharr{&&}
\def\SE##1E##2E{\slantedarrow(0,18)(4,-3){##1}{##2}&}
\def\SW##1W##2W{\slantedarrow(24,18)(-4,-3){##1}{##2}&}
\def\NE##1E##2E{\slantedarrow(0,0)(4,3){##1}{##2}&}
\def\NW##1W##2W{\slantedarrow(24,0)(-4,3){##1}{##2}&}
\def\slantedarrow(##1)(##2)##3##4{%
\thinlines\unitlength1pt\lower 6.5pt\hbox{\begin{picture}(24,18)%
\put(##1){\vector(##2){24}}%
\put(0,8){$\scriptstyle##3$}%
\put(20,8){$\scriptstyle##4$}%
\end{picture}}}
\def\vspace##1{\noalign{\vskip##1\relax}}\relax\let\amp@rs@nd@&\iffalse}\fi
 \CD@true\vcenter\bgroup\relax\let\\=\cr\iffalse}\fi\tabskip\z@skip\baselineskip20\ex@
 \lineskip3\ex@\lineskiplimit3\ex@\halign\bgroup
 &\hfill$\m@th##$\hfill\cr}
\def\@endCDS{\cr\egroup\egroup}
%
\newdimen\TriCDarrw@
\newif\ifTriV@
\newenvironment{TriCDV}{\@TriCDV}{\@endTriCD}
\newenvironment{TriCDA}{\@TriCDA}{\@endTriCD}
\def\@TriCDV{\TriV@true\def\TriCDpos@{6}\@TriCD}
\def\@TriCDA{\TriV@false\def\TriCDpos@{10}\@TriCD}
\def\@TriCD#1#2#3#4#5#6{%
\setbox0\hbox{$\ifTriV@#6\else#1\fi$}
\TriCDarrw@=\wd0 \advance\TriCDarrw@ 24pt
\advance\TriCDarrw@ -1em
\def\SE##1E##2E{\slantedarrow(0,18)(2,-3){##1}{##2}&}
\def\SW##1W##2W{\slantedarrow(12,18)(-2,-3){##1}{##2}&}
\def\NE##1E##2E{\slantedarrow(0,0)(2,3){##1}{##2}&}
\def\NW##1W##2W{\slantedarrow(12,0)(-2,3){##1}{##2}&}
\def\slantedarrow(##1)(##2)##3##4{\thinlines\unitlength1pt
\lower 6.5pt\hbox{\begin{picture}(12,18)%
\put(##1){\vector(##2){12}}%
\put(-4,\TriCDpos@){$\scriptstyle##3$}%
\put(12,\TriCDpos@){$\scriptstyle##4$}%
\end{picture}}}
\def\={\mathrel {\vbox{\hrule
   width\TriCDarrw@\vskip3\ex@\hrule width
   \TriCDarrw@}}}
\def\>##1>>{\setbox\z@\hbox{$\scriptstyle
 \;{##1}\;\;$}\global\bigaw@\TriCDarrw@
 \ifdim\wd\z@>\bigaw@\global\bigaw@\wd\z@\fi
 \hskip.5em
 \mathrel{\mathop{\hbox to \TriCDarrw@
{\rightarrowfill}}\limits^{##1}}
 \hskip.5em}
\def\<##1<<{\setbox\z@\hbox{$\scriptstyle
 \;{##1}\;\;$}\global\bigaw@\TriCDarrw@
 \ifdim\wd\z@>\bigaw@\global\bigaw@\wd\z@\fi
 \mathrel{\mathop{\hbox to\bigaw@{\leftarrowfill}}\limits^{##1}}
 }
 \CD@true\vcenter\bgroup\relax\let\\=\cr\iffalse}\fi
 \tabskip\z@skip\baselineskip20\ex@
 \lineskip3\ex@\lineskiplimit3\ex@
 \ifTriV@
 \halign\bgroup
 &\hfill$\m@th##$\hfill\cr
#1&\multispan3\hfill$#2$\hfill&#3\\
&#4&#5\\
&&#6\cr\egroup%
\else
 \halign\bgroup
 &\hfill$\m@th##$\hfill\cr
&&#1\\%
&#2&#3\\
#4&\multispan3\hfill$#5$\hfill&#6\cr\egroup
\fi}
\def\@endTriCD{\egroup}
\newcommand{\mc}{\mathcal}
\newcommand{\mb}{\mathbb}
\newcommand{\surj}{\twoheadrightarrow}
\newcommand{\inj}{\hookrightarrow}
\newcommand{\zar}{{\rm zar}}
\newcommand{\an}{{\rm an}} 
\newcommand{\red}{{\rm red}}
\newcommand{\codim}{{\rm codim}}
\newcommand{\rank}{{\rm rank}}
\newcommand{\Ker}{{\rm Ker \ }}
\newcommand{\Pic}{{\rm Pic}}
\newcommand{\Div}{{\rm Div}}
\newcommand{\Hom}{{\rm Hom}}
\newcommand{\im}{{\rm im}}
\newcommand{\Spec}{{\rm Spec \,}}
\newcommand{\Sing}{{\rm Sing}}
\newcommand{\sing}{{\rm sing}}
\newcommand{\reg}{{\rm reg}}
\newcommand{\Char}{{\rm char}}
\newcommand{\Tr}{{\rm Tr}}
\newcommand{\Gal}{{\rm Gal}}
\newcommand{\Min}{{\rm Min \ }}
\newcommand{\Max}{{\rm Max \ }}
\newcommand{\Alb}{{\rm Alb}\,}
\newcommand{\GL}{{\rm GL}\,}        
\newcommand{\ie}{{\it i.e.\/},\ }
\newcommand{\niso}{\not\cong}
\newcommand{\nin}{\not\in}
\newcommand{\soplus}[1]{\stackrel{#1}{\oplus}}
\newcommand{\by}[1]{\stackrel{#1}{\rightarrow}}
\newcommand{\longby}[1]{\stackrel{#1}{\longrightarrow}}
\newcommand{\vlongby}[1]{\stackrel{#1}{\mbox{\large{$\longrightarrow$}}}}
\newcommand{\ldownarrow}{\mbox{\Large{\Large{$\downarrow$}}}}
\newcommand{\lsearrow}{\mbox{\Large{$\searrow$}}}
\renewcommand{\d}{\stackrel{\mbox{\scriptsize{$\bullet$}}}{}}
\newcommand{\dlog}{{\rm dlog}\,}    
\newcommand{\longto}{\longrightarrow}
\newcommand{\vlongto}{\mbox{{\Large{$\longto$}}}}
\newcommand{\limdir}[1]{{\displaystyle{\mathop{\rm lim}_{\buildrel\longrightarrow\over{#1}}}}\,}
\newcommand{\liminv}[1]{{\displaystyle{\mathop{\rm lim}_{\buildrel\longleftarrow\over{#1}}}}\,}
\newcommand{\norm}[1]{\mbox{$\parallel{#1}\parallel$}}
\newcommand{\boxtensor}{{\Box\kern-9.03pt\raise1.42pt\hbox{$\times$}}}
\newcommand{\into}{\hookrightarrow}
\newcommand{\image}{{\rm image}\,}
\newcommand{\Lie}{{\rm Lie}\,}      
\newcommand{\CM}{\rm CM}
\newcommand{\sext}{\mbox{${\mathcal E}xt\,$}}  
\newcommand{\shom}{\mbox{${\mathcal H}om\,$}}  
\newcommand{\coker}{{\rm coker}\,}  
\newcommand{\sm}{{\rm sm}}
\newcommand{\tensor}{\otimes}
\renewcommand{\iff}{\mbox{ $\Longleftrightarrow$ }}
\newcommand{\supp}{{\rm supp}\,}
\newcommand{\ext}[1]{\stackrel{#1}{\wedge}}
\newcommand{\onto}{\mbox{$\,\>>>\hspace{-.5cm}\to\hspace{.15cm}$}}
\newcommand{\propsubset}
{\mbox{$\textstyle{
\subseteq_{\kern-5pt\raise-1pt\hbox{\mbox{\tiny{$/$}}}}}$}}
\newcommand{\sA}{{\mathcal A}}
\newcommand{\sB}{{\mathcal B}}
\newcommand{\sC}{{\mathcal C}}
\newcommand{\sD}{{\mathcal D}}
\newcommand{\sE}{{\mathcal E}}
\newcommand{\sF}{{\mathcal F}}
\newcommand{\sG}{{\mathcal G}}
\newcommand{\sH}{{\mathcal H}}
\newcommand{\sI}{{\mathcal I}}
\newcommand{\sJ}{{\mathcal J}}
\newcommand{\sK}{{\mathcal K}}
\newcommand{\sL}{{\mathcal L}}
\newcommand{\sM}{{\mathcal M}}
\newcommand{\sN}{{\mathcal N}}
\newcommand{\sO}{{\mathcal O}}
\newcommand{\sP}{{\mathcal P}}
\newcommand{\sQ}{{\mathcal Q}}
\newcommand{\sR}{{\mathcal R}}
\newcommand{\sS}{{\mathcal S}}
\newcommand{\sT}{{\mathcal T}}
\newcommand{\sU}{{\mathcal U}}
\newcommand{\sV}{{\mathcal V}}
\newcommand{\sW}{{\mathcal W}}
\newcommand{\sX}{{\mathcal X}}
\newcommand{\sY}{{\mathcal Y}}
\newcommand{\sZ}{{\mathcal Z}}
\newcommand{\A}{{\mathbb A}}
\newcommand{\B}{{\mathbb B}}
\newcommand{\C}{{\mathbb C}}
\newcommand{\D}{{\mathbb D}}
\newcommand{\E}{{\mathbb E}}
\newcommand{\F}{{\mathbb F}}
\newcommand{\G}{{\mathbb G}}
\newcommand{\HH}{{\mathbb H}}
\newcommand{\I}{{\mathbb I}}
\newcommand{\J}{{\mathbb J}}
\newcommand{\M}{{\mathbb M}}
\newcommand{\N}{{\mathbb N}}
\renewcommand{\P}{{\mathbb P}}
\newcommand{\Q}{{\mathbb Q}}
\newcommand{\R}{{\mathbb R}}
\newcommand{\T}{{\mathbb T}}
\newcommand{\U}{{\mathbb U}}
\newcommand{\V}{{\mathbb V}}
\newcommand{\W}{{\mathbb W}}
\newcommand{\X}{{\mathbb X}}
\newcommand{\Y}{{\mathbb Y}}
\newcommand{\Z}{{\mathbb Z}}
\title[The universal regular quotient of the Chow group of 0-cycles]  
{The universal regular quotient of the Chow group of 0-cycles on a
singular projective variety}  
\author{H\'el\`ene Esnault} 
\address{Universit\"at Essen, FB6 Mathematik, 45 117 Essen, Germany}
\email{esnault@uni-essen.de}
\thanks{This work has been partly supported by the DFG Forschergruppe
''Arithmetik und Geometrie''}
\author{V. Srinivas} 
\address{School of Mathematics, TIFR,
Homi Bhabha Road, Bombay-400005, India}
\email{srinivas@math.tifr.res.in}
\author{Eckart Viehweg} 
\address{Universit\"at Essen, FB6 Mathematik, 45 117 Essen, Germany}
\email{viehweg@uni-essen.de}
\maketitle
Let $X$ be a projective variety of dimension $n$ defined over an
algebraically closed field $k$. For $X$ irreducible and non-singular,
Matsusaka \cite{Matsusaka} constructed an abelian variety $\Alb(X)$ and a
morphism $\alpha:X\to \Alb(X)$ (called the Albanese variety and mapping
respectively), depending on the choice of a base-point on $X$, which is
universal among the morphisms to abelian varieties (see Lang \cite{La},
Serre \cite{Se} for other constructions). Over the field of complex
numbers the existence of $\Alb(X)$ and $\alpha$ was known before, and has
a purely Hodge-theoretic description (see Igusa \cite{Igusa} for the Hodge
theoretic construction). Incidentally, the terminology ``Albanese
variety'' was introduced by A.~Weil, for reasons explained in his
commentary on the article [1950a] of Volume~I of his collected works (see
\cite{Weil}), one of which is that the paper \cite{Albanese} of Albanese
defines it (for a surface) as a quotient of the group of 0-cycles of
degree 0 modulo an equivalence relation.

Let $CH^n(X)_{\deg 0}$ denote the Chow group of 0-cycles of degree 0 on
$X$ modulo rational equivalence. When $X$ is irreducible and non-singular,
a remarkable feature of the Albanese morphism $\alpha$ is that it factors
through a regular homomorphism $\varphi: CH^n(X)_{\deg 0} \to \Alb(X)$,
that is a homomorphism, which when composed with 
the cycle map $\gamma: X \to CH^n(X)_{\deg 0}$, gives an algebraic
morphism. This follows immediately from the fact that
an abelian variety does not contain any rational curve.  Thus one can
reformulate Matsusaka's theorem as the statement that there is a universal
regular quotient of $CH^n(X)_{\deg 0}$ as an abelian variety.

For $n=1$, $\varphi$ is an isomorphism, and moreover
$CH^1(X)_{\deg 0} \cong \Pic^0(X)$, the Picard variety of $X$.
Since $\Pic$ is a well understood functor, one defines the Chow
group $CH^1(X)_{\deg 0}$ and the generalized Albanese variety $A^1(X)$
by 
$$
CH^1(X)_{\deg 0} \cong A^1(X) \cong \Pic^0(X)
$$
even in the singular case.
This has several consequences for the expected structure of the
generalized Albanese variety $A^n(X)$ of a projective reduced variety $X$
of dimension $n$. First, it forces the correct definition of the Chow
group $CH^n(X)_{\deg 0}$, as proposed by Levine and Weibel in \cite{LW}.
Second, it shows that $A^n(X)$ should be a smooth commutative algebraic
group, that is an extension of an abelian variety by a linear group,
where the latter is a product of additive and multiplicative factors.
Third, the cycle map to $CH^n(X)_{\deg 0}$ is only defined on the
regular locus $X_{\reg}$ of $X$. Consequently the expected generalized
Albanese mapping should be only defined on $X_{\reg}$. But already for
curves, a morphism from $X_{\reg}$ to a smooth commutative algebraic group
$G$ need not factor through $CH^1(X)_{\deg 0}$, as $G$ contains rational
subvarieties. Therefore, the expected $A^n(X)$ should be constructed as a
regular quotient of $CH^n(X)_{\deg 0}$ in the category of smooth
commutative algebraic groups.  Note here that the difficulty comes from
the non-normality of $X$.  In fact for normal surfaces and for irreducible
normal varieties in characteristic zero (see
\cite{S}), it is known that $A^n(X) = \Alb(\tilde{X})$, where
$\tilde{X}$ is a resolution of singularities, and that the cycle
map does factor through $CH^n(X)_{\deg 0}$.

Roitman \cite{R} proved that, when $k$ is a universal domain, $\varphi$ is
an isomorphism precisely when the Chow group of 0-cycles is {\em finite
dimensional} in the sense of Mumford \cite{M}.  This result was
generalized for irreducible normal varieties in characteristic zero for
$A^n(X) = \Alb(\tilde{X})$ (see \cite{S}).

In this article we prove the existence of a universal regular
quotient of the Chow group of $0$-cycles 
for singular projective varieties. We note that the term
``variety'' is used to mean a reduced quasi-projective 
scheme of finite type over a field; in particular  it need not
be irreducible, or equidimensional. A regular
homomorphism is defined in \ref{def-reg} and the finite
dimensionality of the Chow group in \ref{fin-dim}.
\begin{thm1}\label{thm1}
Let $X$ be a projective variety of dimension $n$, defined over
an algebraically closed field $k$. 
\begin{points}
\item There exists a smooth
connected commutative algebraic group $A^n(X)$, together with a regular
homomorphism $\varphi:CH^n(X)_{\deg 0}\to A^{n}(X)$, such that $\varphi$
is universal among regular homomorphisms from $CH^n(X)_{\deg 0}$ to smooth
commutative algebraic groups. 
\item Over a universal domain $k$ the Chow group is finite
dimensional precisely when $\varphi$ is an isomorphism.
\item $A^n(X\times_kK)=A^n(X)\times_kK$, for all algebraically
closed fields $K$ containing $k$.
\end{points}
\end{thm1}

We also give a second construction of $A^n(X)$ and $\varphi$ using
transcendental arguments when $k=\C$.  Over $k=\C$, there is a natural
semi-abelian variety   
\[J^n(X)=\frac{H^{2n-1}(X,\C(n))}{F^0H^{2n-1}(X,\C(n))+
{\rm image}\,H^{2n-1}(X,\Z(n))}, \]
that is a commutative algebraic group without additive factors,
whose construction is implicit in Deligne's article \cite{D}. 
(For $k\propsubset \C$, one can in fact define $J^n(X)$ over $k$, see
\cite{BS} and also \cite{FW}). 
{From} the discussion above, one sees that $A^n(X)$ cannot be
isomorphic to $J^n(X)$. However, there an Abel-Jacobi mapping 
$$ AJ^n: CH^n(X)_{\deg 0} \>>> J^n(X)$$
with very good properties.
For example,  if $X$ is irreducible and non-singular, then Roitman
\cite{R2} proved that the Albanese mapping $\varphi$ is an isomorphism on
torsion subgroups. For any reduced, projective $X$ of dimension $n$, the
Abel-Jacobi map $AJ^n$ induces an isomorphism on torsion subgroups as well 
(see \cite{BiS} for the general result, and \cite{L},
\cite{C}, \cite{BPW} for earlier partial results). This indicates that
$J^n(X)$ should differ from $A^n(X)$ only by additive factors. 
This, together with the classical theory for curves, was the main motivation
for our construction. For related constructions using
K\"ahler differentials in the singular case, see \cite{L}.
 
\begin{thm2}\label{thm2}
Let $X$ be a projective variety over $\C$. For any $m\geq 0$, define the
{\em Deligne complex} 
\[\sD(m)_X=\left(0\>>> \Z_X(m)\>>>\sO_X\>>>\Omega^{1}_{X/\C}\>>>\cdots\>>>
\Omega^{m-1}_{X/\C}\>>> 0\right),\]
and associated cohomology group $D^{m}(X)=\HH^{2m}(X,\sD(m)_X).$
For $n=\dim X$, let
\[A^{n}(X)=\ker\left(D^{n}(X)\>>> H^{2n}(X,\Z(n))\right)\]
be the kernel of the map induced by the natural surjection 
$\sD(n)_X\to\Z(n)_X$ of complexes. Then
\begin{points}
\item
the analytic group $A^n(X)$ has an underlying algebraic
structure and for some $s\geq 0$ a presentation via an exact sequence of
commutative algebraic groups
$$0 \>>> (\G_a)^s \>>>A^n(X) \>>> J^n(X) \>>> 0$$
\item
there is a cycle class homomorphism $CH^{n}(X)\to D^{n}(X),$
such that the composite $CH^{n}(X)\to H^{2n}(X,\Z(n))$
is the degree homomorphism, induced by the topological cycle
class map, and giving rise to a commutative diagram
\[\begin{CD}
CH^{n}(X)_{\deg 0} \>{\varphi}>> A^{n}(X) \\
\V V V \V V V \\
CH^{n}(X) \>>> D^{n}(X)
\end{CD}
\]
\item
the homomorphism $\varphi$ is the universal regular homomorphism
from the Chow group $CH^{n}(X)_{\deg 0}$ to commutative
algebraic groups over $\C$
\item $\varphi$ is an isomorphism precisely when the Chow group is finite
dimensional.
\end{points}
\end{thm2}

We give examples in \ref{example} illustrating two pathological properties
of $A^n(X)$. First we give examples of irreducible projective varieties
$X$ and $Y$ of dimensions $n$ and $m$, respectively, for which 
$\dim(A^{n+m}(X\times Y)) > \dim(A^n(X) \times A^m(Y))$.
Next, we exhibit a flat family $\sX\to S$, with geometrically integral fibres, 
for which the dimension of $A^n(X_s)$ is not locally constant on $S$.

We do not study  higher dimensional cycles in this article.
For this reason we do not analyze $D^m(X)$, $A^m(X)=\ker(D^m(X)\to
H^{2m}(X,\Z(m))$ and their relation to the Chow ring
$CH^\bullet(X,X_{\rm sing})$ considered by Levine in \cite{L2}.

After recalling the definition (\cite{LW}) of the Chow group
$CH^n(X)_{\deg 0}$ and its relation to the Picard group of curves, and a
moving lemma from \cite{BiS}, we construct the analytic cycle class map in
section 2 and prove theorem 2, (ii). In section 3 we prove theorem 2, (i)
and give a cohomological description for the Lie algebra of $A^n(X)$. The
regularity of $\varphi:CH^n(X)_{\deg 0}\to A^n(X)$ and the universal
property (iii) for $\varphi$ are shown by analytic methods in section 4.

The next two sections, independent of the transcendental arguments
used before, contain the algebraic part of the article.
We recall in section 5 some of the properties of Picard groups
of curves and apply them to general curves in $X$. The main
technical tool in section 6 is the boundedness 
of the dimension of a regular quotient. 
This being guaranteed $A^n(X)$ is constructed essentially
using Lang's arguments \cite{La}. 

In section 7 we give two slightly different proofs for part
(iv) of theorem 2 and and (ii) of theorem 1, the first one building up on
the transcendental methods, the second one using the algebraic
arguments developed in sections 1, 5 and 6.
\section{Chow groups and regular homomorphisms}
We begin by recalling the definition of the Chow group of
0-cycles $CH^n(X)$, as given in \cite{LW} (see also \cite{BiS}). As in
\cite{BiS}, we adopt the convention that a point lying on a lower
dimensional component of $X$ is deemed to be singular. Let $X_{\sing}$
denote the (closed) subset of singular points, and $X_{\reg}=X-X_{\sing}$ the
complementary open set. The closure of $X_{\reg}$ is the union of the
$n$-dimensional components of $X$.

The group $Z^n(X)$ of 0-cycles is defined to be the free abelian
group on the closed points of $X_{\reg}$. The subgroup $R^n(X)$ of cycles
rationally equivalent to 0 is defined using the notion of a
Cartier curve.
\begin{defn}\label{ccurve}
A {\em Cartier curve} is a subscheme $C\subset X$, defined over
$k$, such that 
\begin{points}
\item $C$ is pure of dimension 1
\item no component of $C$ is contained in $X_{\sing}$ 
\item if $x\in C\cap X_{\sing}$, then the ideal of $C$ in $\sO_{x,X}$ is
generated by a regular sequence (consisting of $n-1$ elements).
\end{points}
\end{defn}
If $C$ is a Cartier curve on $X$, with generic points
$\eta_1,\ldots,\eta_s$, and $\sO_{S,C}$ is the
semilocal ring on $C$ of the points of $S=(C\cap
X_{\sing})\cup\{\eta_1,\ldots,\eta_s\}$, there is a natural map on
unit groups
\[\theta_{C,X}:\sO_{S,C}^*\>>> \bigoplus_{i=0}^s\sO_{\eta_i,C}^*.\]
Define $R(C,X)=\image\theta_{C,X}$. For $f\in R(C,X)$,
define the divisor of $(f)_C$ as 
follows: let $C_i$ denote the maximal Cohen-Macaulay subscheme
of $C$ supported on the component with generic point $\eta_i$.
Then for any $x\in C_i$ the map
\[\sO_{x,C_i}\>>>\sO_{\eta_i,C_i}=\sO_{\eta_i,C}\]
is the injection of a Cohen-Macaulay local ring of dimension 1
into its total quotient ring. If $f_i$ is the component of $f$ in
$\sO_{\eta_i,C}$, then $f_i=a_x/b_x$ for some non zero-divisors 
$a_x,b_x\in\sO_{x,C_i}$.  Define 
\[(f)_C=\sum_{i=1}^s(f_i)_{C_i}=\sum_{i=1}^s\sum_{x\in
C_i}(\ell(\sO_{x,C_i}/a_x\sO_{x,C_i})-\ell(\sO_{x,C_i}/b_x\sO_{x,C_i}))
\cdot [x].\]
Standard arguments imply that this is well-defined (\ie the
coefficient of $[x]$ is independent of the choice of the
representation $f_i=a_x/b_x$, and vanishes for all but a finite
number of $x$). 

Suppose $C$ is reduced. Then in the above considerations, $\sO_{x,C_i}$ is
an integral domain with quotient field $\sO_{\eta_i,C}$. If $v_1,
\ldots,v_m$ are the discrete valuations of $\sO_{\eta_i,C}$ centered at
$x$, then the multiplicity of $x$ in $(f)_{C_i}$ is
\begin{equation}\label{val}
\ell(\sO_{x,C_i}/a_x\sO_{x,C_i})-\ell(\sO_{x,C_i}/b_x\sO_{x,C_i})=
\sum_{j=1}^m v_j(f_i)
\end{equation}
(compare \cite{Fulton}, Example~A.3.1.).
In fact, let $R$ be the integral closure of $O=\sO_{x,C_i}$ in
$\sO_{\eta_i,C}$. The Chinese remainder theorem implies that
\[\ell(R/a_xR) = \sum_{j=1}^m v_j(a_x),\]
 and similarly for $b_x$.  Multiplying $a_x$ and $b_x$ by the same 
element of $O$ we may assume that both $a_xR$ and $b_xR$ are contained 
in $O$, and $$\ell(O/a_xO)+\ell(R/O) = \ell(R/a_xR) + \ell(a_xR/a_xO).$$
Since $a_x \neq 0$ the second terms on both sides are equal. 
\begin{defn}\label{rat-equ}
Let $U \subset X_\reg$ be an open dense subscheme.
$R^n(X,U)$ is defined to be the subgroup of $Z^n(U)$
generated by elements $(f)_C$ as $C$ ranges over all Cartier
curves with $C\cap U$ dense in $C$, and $f\in R(C,X)$
with $(f)_C \in Z^n(U)$. 
For $U=X_{\reg}$ we write $R^n(X)$ instead of $R^n(X,X_{\reg})$
and define 
$$CH^n(X) = Z^n(X)/R^n(X).$$
Mapping a point $x \in X_{\reg}$ to its rational equivalence
class defines a map 
$$\gamma :X_{\reg} \>>> CH^n(X).$$
If $U_1, \ldots , U_r$ denote the irreducible components of $X_\reg$
then $Z^n(X)_{\deg 0}$ and $CH^n(X)_{\deg 0}$
denote the subgroup of $Z^n(X)$ and $CH^n(X)$, respectively, of
cycles $\delta$ with $\deg(\delta|_{U_i})=0$
for $i = 1, \ldots r$.
\end{defn}
As noted in \cite{BiS}, lemma~1.3 of \cite{LW} allows one to
restrict to considering only curves $C$ such that $C\cap
X_{\reg}$ has no embedded points, and any irreducible component
$C'$ of $C$ which lies entirely in $X_{\reg}$ occurs in $C$ with
multiplicity 1. The moving lemmas 2.2.2 and 2.2.3 of \cite{BiS}
allow stronger restrictions on $C$:
\begin{lemma}\label{moving}
Let $A \subset X_{\sing}$ be a closed subset of dimension $\leq n-2$,
and let $D \subset X$ be a closed subset of dimension $\leq n-1$. Then any 
element $\delta\in R^n(X)$ can be written in the form 
$\delta=(f)_C$
for a single (possibly reducible) Cartier curve $C$, such that
\begin{enumerate}
\item[(a)] $C$ is reduced
\item[(b)] $C\cap A=\emptyset$
\item[(c)] $C\cap D$ is empty or consists of finitely many points.
\end{enumerate}
\end{lemma}
\begin{cor}\label{equidim} 
If $U \subset X_\reg$ is an open and dense subscheme, then
$$CH^n(X)=Z^n(U)/R^n(X,U) \mbox{ \ \ and \ \ }CH^n(X)_{\deg
0}=Z^n(U)_{\deg 0}/R^n(X,U).$$ 
\end{cor}
\begin{proof}
First note that the zero cycles supported on $U$ generate
$CH^n(X)$ since the corresponding assertion holds true for
curves. The moving lemma \ref{moving} for $D=X-U$
implies that $R^n(X)\cap Z^n(U)=R^n(X,U)$. 
\end{proof}
\begin{rmk}
Let $X^{(n)}$ denote the union of the $n$-dimensional irreducible
components of $X$, and let $X^{<n}$ be the union of the lower dimensional
components. Applying the corollary to $X^{(n)}$ and the open subset 
$U=X^{(n)}-X_{\sing}=X-X_{sing}$, we see that the natural map from
$CH^n(X)$ to $CH^n(X^{(n)})$ is surjective. It seems plausible that a
stronger form of lemma~\ref{moving} holds, where $A$ is allowed to be any
closed subset of $X$ of codimension $\geq 2$ which is disjoint from
$\supp(\delta)$. If this is true, then applying it to $X^{(n)}$ with
$A=X^{(n)}\cap X^{<n}$,  one sees that for any $\delta\in R^n(X^{(n)})\cap
Z^n(X)$ there exists a reduced  Cartier curve $C$ in $X^{(n)}$, disjoint from
$A$, and $f \in R(C,X^{(n)})$ with $\delta=(f)_C$. Then $C$ is also a
Cartier curve on $X$, and $\delta\in R^n(X)$. We deduce that $CH^n(X)\to
CH^n(X^{(n)})$ is an isomorphism. We have as yet been unable to prove this. 
\end{rmk}

\begin{rmk}\label{rmk-equi} Keeping the notation from the previous remark,
we note further that for $k=\C$, the natural maps 
\[H^{2n}(X,\Z(n))\>>> H^{2n}(X^{(n)},\Z(n)),\;\; D^n(X)\>>> 
D^n(X^{(n)}),\;\; A^n(X)\>>> A^n(X^{(n)})\]
are isomorphisms, since $X^{<n}$ has constructible cohomological
dimension $\leq 2(n-1)$ and coherent cohomological dimension $\leq n-1$.
\end{rmk}

As reflected by the notation, $R(C,X)$ depends on the pair
$(C,X)$, and is not necessarily intrinsic to $C$. In fact,
since we have not imposed any unit condition at singular points
of $C$ which lie in $X_{\reg}$, the functions $f \in R(C,X)$ are
defined on some curve $C'$, birational to $C$.
\begin{defn}\label{admissible}
Let $C'$ be a reduced projective curve and $\iota:C' \to X$ be a
morphism. Then $(C',\iota)$ will be called {\it admissible}
if $\iota: C' \to C=\iota(C')$ is birational, if $C$ is a
reduced Cartier curve and if for some open neighbourhood $W$ of
$X_{\sing}$ the restriction of $\iota$ to $\iota^{-1}(W)$ is a
closed embedding.
\end{defn}
If $(C',\iota)$ is admissible one has an inclusion $R(C',C')
\subset R(\iota(C),X)$ which is an equality if 
$\iota^{-1}(X_{\reg})$ is non-singular. 
\begin{lemma}\label{gysin}
Let $(C',\iota)$ be admissible. Then there exists a homomorphism
(of abstract groups)
$$ \eta: \Pic^0(C') \cong CH^1(C')_{\deg 0} \>>> CH^n(X)_{\deg 0}$$
which maps the isomorphism class of $\sO_{C'}(p-p')$ to $\gamma
(\iota(p))-\gamma(\iota(p'))$. 
\end{lemma}
\begin{proof}
By definition $\Pic(C')= Z^1(C'_\reg)/R(C',C')$ and one has
a map
$$\gamma \circ \iota : Z^1(C'_\reg) \>>> CH^n(X).
$$
The equality (\ref{val}) shows that for $f\in R(C',C')$
the image of $(f)_{C'}$ in $CH^n(X)$ is zero.
\end{proof}
\begin{notations}\label{difference} 
Let $Y$ be a non-singular scheme with irreducible components $Y_1,
\ldots ,Y_s$, let $G$ be an abstract or an algebraic group, and let
$\pi: Y \to G$ a map or morphism. 
\begin{points} 
\item After choosing base points $p_i \in Y_i$ a map
$$
\pi_m : S^m(Y):= S^m( \bigcup_{i=1}^s Y_i) \>>> G,
$$
is defined by $\pi_m(y_1, \ldots, y_m)=\sum_{j=1}^m
(\pi(y_j) - \pi(p_{\rho(j)}))$, where $\rho(j) = i$ if $y_j \in Y_i$. 
\item To avoid the reference to base points, we will frequently use
different maps:
$$
\pi^{(-)} : \Pi_Y = \bigcup_{i=1}^s Y_i\times Y_i \>>> G
$$
is defined by $\pi^{(-)}(y,y') = \pi(y)-\pi(y')$, and
$\pi^{(-)}_m: S^m(\Pi_Y) \to G$ is the composite
$S^m(\Pi_Y) \to S^m(G) \> {\rm sum} >> G$. 
\end{points}
\end{notations}
If $G$ is an algebraic group, then the images of
$\pi^{(-)}_m$ lie in the connected component of $0$.
In particular for $U$ open and dense
in $X_\reg$ we will frequently consider
\begin{align*}
&\gamma^{(-)} = \gamma_U^{(-)}: \Pi_U \>>> CH^n(X)_{\deg 0}\\
\mbox{and \ \ \ \ }
&\gamma_m = \gamma_{U,m} : S^m(U) \>>> CH^n(X)_{\deg 0}.
\end{align*}
\begin{lemma} \label{generators}
Let $G$ be a $d$-dimensional
smooth connected commutative algebraic group
and let $\Gamma \subset G$ be a constructible subset which
generates $G$ as an abstract group. 
Then 
\begin{points}
\item the image of the
composite map \ $S^d(\Gamma) \>>> S^d(G) \>{\rm sum} >> G $ \
is dense
\item 
$S^{2d}(\Gamma) \>>> S^{2d}(G) \>{\rm sum} >> G$
\  
is surjective
\item if $B$ is a non-singular scheme with connected components
$B_1, \ldots ,B_s$ and if $\vartheta: B \to G$ is a morphism
with image $\Gamma$ then the morphism 
$$
\vartheta_{d}^{(-)}: S^d(\Pi_B)=S^d \bigl( \bigcup_{i=1}^s
B_i \times B_i \bigr) \>>> G
$$
with $\vartheta_{d}^{(-)}((b_1,b_1'),\ldots (b_d,b_d')) = \sum_{i=1}^d
(\vartheta(b_i)-\vartheta(b_i'))$ is surjective.
\end{points}
\end{lemma}
\begin{proof}
Let $\bar{\Gamma}_1, \ldots, \bar{\Gamma}_s$ be the irreducible
components of the closure $\bar{\Gamma}$ of $\Gamma$, and let
$\Gamma_i= \bar{\Gamma_i}\cap \Gamma$. 
It is sufficient to find non-negative integers
$d_1, \ldots ,d_s$ with $\sum_{i=1}^s d_i \leq d$ such that
the image of
$S^{d_1}(\Gamma_1)\times \cdots \times S^{d_s}(\Gamma_s)$
is dense in $G$. To this end, we may assume that the identity
of $G$ lies on each $\Gamma_i$. 

Let $\bar{\Gamma}_1^\nu$ be the closure of the image of
$S^\nu(\Gamma_1)$ in $G$. Since
$\bar{\Gamma}_1^\nu \subset \bar{\Gamma}_1^{\nu+1}$
there exists some $d_1 \leq d$ with 
$\bar{\Gamma}_1^{d_1} = \bar{\Gamma}_1^{d_1+1}$,
and $d_1$ is minimal with this property.
Hence $\bar{\Gamma}_1^{d_1} = \bar{\Gamma}_1^{2 d_1}$ 
and $\bar{\Gamma}_1^{d_1}$ is a subgroup of $G$
of dimension larger than or equal to $d_1$. 
If $s=1$, i.e. if $\bar{\Gamma}$ is irreducible, then
$\bar{\Gamma}_1^{d_1}=G$.

In general, replacing $G$ by $G/\bar{\Gamma}_1^{d_1}$ one
obtains \ref{generators} (i) by induction on $s$.

The second part is an easy consequence of (i). Let $U$
be an open dense subset of $G$, contained in
$S^{d_1}(\Gamma_1)\times \cdots \times S^{d_s}(\Gamma_s)$.
Given $p \in G$ the intersection of the two open sets $U$ and
$p-U$ is non-empty and hence there are points $a,b \in U$ with
$p-b=a$.

For (iii) we may assume that for some point $b_i' \in B_i$, the image of
$B_i\times\{b'_i\}$ in $\Gamma_i$ is dense, for each $i$.
By (i) one finds
$d_1, \ldots ,d_s$ with $\sum_{i=1}^s d_i = d$ such that
the image 
$$
\vartheta^{(-)}_d(S^{d_1}(B_1\times \{b_1'\})\times \cdots \times
S^{d_s}(B_s\times \{b_s'\}))
$$
contains a subset $U$ which is open in $G$. Given $p \in G$
the intersection of $U$ and of $p+U$ is non empty and hence
$p=a-b$ for two points $a$ and $b$ in $U$. Obviously $a-b$ lies
in the image of $\vartheta^{(-)}_d$.
\end{proof}
\begin{cor}\label{generators2} Let $C'$ be a reduced curve,
let $B_1, \ldots B_s$ be the connected components of $B=C'_\reg$,
let $b'_j \in B_j$ be base points and let $\vartheta: B \to
\Pic^0(C')$ be the morphism with $\vartheta|_{B_j}(b)=\sO_{C'}(b-b_j')$.
Then there exists some open connected subscheme $W$ of $S^g(B)$,
for $g=\dim_k(\Pic^0(C'))$, such that 
$\vartheta_W := \vartheta_g|_W$ 
is an open embedding.
\end{cor}
\begin{proof}
By \ref{generators} we find some $W$ with $\vartheta_W(W)$ open
and $\vartheta_W$ finite over its image. On the other hand, any
fibre of $\vartheta_g$ is an open subset of $\P(H^0(C',\sO_{C'}(D)))$
for some divisor $D$ on $C'$; hence the projective spaces corresponding to
points of $\vartheta_W(W)$ must be 0-dimensional.
\end{proof}
\begin{lemma}\label{pic}
Let $G$ be a smooth commutative algebraic group, $U \subset X_\reg$ an
open and dense subset, and $\pi: U \to G$ a morphism. 
Then the following two conditions are equivalent.
\begin{enumerate}
\item[(a)] There exists a homomorphism (of abstract groups)
$\phi:CH^n(X)_{\deg 0} \to G$ such that $\pi^{(-)} = \phi \circ
\gamma^{(-)}$ (as maps on the closed points).
\item[(b)] For all admissible pairs $(C',\iota)$ with
$B=(\iota^{-1}(U))_\reg$ dense in $C'$ there exists a
homomorphism of algebraic groups $\psi:\Pic^0(C') \to G$ such that
the diagram 
$$
\begin{CD}
\Pi_B \>\vartheta^{(-)} >> \Pic^0(C')\\
\V \iota V V  \V V \psi V \\
\Pi_U \> \pi^{(-)} >> G
\end{CD}
$$
commutes. Here $\vartheta: B \to \Pic(C')$ denotes the natural
morphism, mapping a point $p$ to the isomorphism class of the
invertible sheaf $\sO_{C'}(p)$.
\end{enumerate}
Moreover, if the equivalent conditions (a) and (b) are true, the
morphism $\psi$ in (b) factors as
$$ \begin{TriCDV}
{\Pic^0(C')}{\> \eta >>}{CH^n(X)_{\deg 0}}
{\SE \psi E E }{\SW W \phi W}{G}
\end{TriCDV}
$$
and the image of $\phi:CH^n(X)_{\deg 0}\to G$ is contained in the
connected component of the identity of $G$.
\end{lemma}
\begin{proof}
Assume (a) and let $(C',\iota)$ be admissible and $g=\dim(\Pic^0(C'))$.
Choosing base points $b'_j \in B_j$, one finds by \ref{generators2} an
open subscheme $W$ of $S^{g}(B)$ such that the morphism $\vartheta_W:W \to
\Pic^0(C')$ is an open embedding. By \ref{gysin} one obtains a
homomorphism $$
\psi:\Pic^0(C')\>\eta>>CH^n(X)_{\deg 0} \> \phi >> G
$$
of abstract groups. By assumption $\pi^{(-)}$ is a morphism of
schemes and the same holds true for $\pi^{(-)}\circ\iota: \Pi_B \to G.$
Thereby the restriction of $\psi$ to the open subscheme
$W \subset \Pic^0(C')$ is a morphism of schemes, and
being a homomorphism of abstract groups $\psi$ is a morphism of
algebraic groups. 

Since each point of $X_\reg$ lies on some Cartier curve, 
the images of the connected algebraic groups $\Pic^0(C')$
generate $CH^n(X)_{\deg 0}$ and the image $\phi(CH^n(X)_{\deg
0})$ lies in the connected component of $G$, which contains the
identity. 

The morphism $\pi^{(-)}$ induces a map
$\tilde{\phi}:Z^n(U)_{\deg 0} \to G$ and it remains to verify
that (b) implies that $\tilde{\phi}(R^n(X,U))=0$. By
\ref{moving} each $\delta\in R^n(X,U)$ is of the form $(f)_C$
for a reduced Cartier curve $C$. There exists an admissible pair
$(C',\iota)$ with $\iota(C')=C$ and with $\iota^{-1}(X_\reg)$
non-singular. 
$(f)_C$ is the image of $(f)_{C'}$ in $Z^n(U)$ and 
by assumption $\iota\circ\pi^{(-)}$ factors through $\Pic^0(C')$.
\end{proof}
\begin{cor}\label{equ-reg}
Let $\phi: CH^n(X)_{\deg 0} \to G$ be a homomorphism to
a smooth commutative algebraic group $G$. Then the following
conditions are equivalent: 
\begin{points}
\item $\phi\circ\gamma^{(-)}: \Pi_{X_\reg} \>>> G$
is a morphism of schemes.
\item There exists an open dense subscheme $U$ of $X_\reg$
such that $\phi\circ\gamma^{(-)}|_{\Pi_U}$ is a morphism of
schemes. 
\item Given a base point $p_i$ on each irreducible component
$U'_i$ of some open dense subscheme $U$ of $X_\reg$, the map
$\pi: U \to G$ with $\pi|_{U'_i}(x) = \phi(x-p_i)$ is a morphism
of schemes. 
\item Given any $m>0$ and base points $p_i$ on each irreducible
component $U'_i$ of some open dense subscheme $U$ of $X_\reg$,
$\phi\circ\gamma_m : S^m(U) \to G$ is a morphism of schemes. 
\end{points}
\end{cor}
Of course, ``$\pi$ is a morphism of schemes'' stands for ``there
exists a morphism of schemes whose restriction to closed points
coincides with $\pi$'', an abuse of terminology which we will repeat
throughout this article.
\begin{proof}
Obviously (i) implies (ii). For $U \subset X_\reg$ given, the
equivalence of (ii), (iii), and (iv) is an easy exercise. In
fact, the morphism $\pi$ in \ref{equ-reg} (iii) is just
$\phi\circ\gamma_1$. 

Assume that (iii) holds true for some $U$.
We will show that the corresponding property holds true for
$X_\reg$ itself. To this aim consider the map
$\bar{\pi}: X_\reg \to G$ with $\bar{\pi}|_{U'_i}(x) =
\phi(x-p_i)$ and the graph $\Gamma_{\bar{\pi}}$ of $\bar{\pi}$
in $X_\reg\times G$. By definition, $\Gamma_{\bar{\pi}} \cap
U\times G$ is the graph $\Gamma_{\pi}$. Let $Z$ be the closure
of $\Gamma_\pi$ in $X_\reg\times G$. 

$\Gamma_{\bar{\pi}}$ is contained in $Z$.
In fact, given a point $x \in X_\reg$ one can find a Cartier
curve $C$ through $x$ with $U\cap C$ dense in $C$ and with
$B=C \cap X_\reg$ non-singular. By lemma \ref{pic} the
morphism $(\pi|_{C\cap U})^{(-)} : \Pi_{C\cap U} \to G$
factors through a morphism $\Pic^0(C)\to G$ of algebraic groups
and, in particular, it extends to a morphism 
$\Pi_{B} \to G$.
Again this implies that the restriction of $\bar{\pi}$ to
$B$ is a morphism, hence $\Gamma_{\bar{\pi}}\cap B\times G$
is closed and therefore contained in $Z$.

By construction the morphism $p_1: Z \to X_\reg$ induced by the
projection is birational and surjective. Let $V \subset X_\reg$
be the largest open subscheme with $p_1|_{p_1^{-1}(V)}$ an
isomorphism. Then $\bar{\pi}|_V$ is a morphism of schemes and 
$\codim_{X_\reg}(X_\reg-V) \geq 2$. By theorem 1 in
\cite{BLR}, 4.4, $\bar{\pi}|_V$ extends to a morphism $X_\reg
\to G$. The graph of this morphism is contained in $Z$, hence it
is equal to $Z$ and $\bar{\pi}$ is a morphism.
\end{proof}

We end this section by giving the definition of a
regular homomorphism, used already in the formulation of the
main theorems in the introduction. 

\begin{defn}\label{def-reg}
Let $G$ be a smooth commutative algebraic group.
A homomorphism $\phi:CH^n(X)_{\deg 0} \to G$ (of abstract groups) is
called {\it a regular homomorphism}, if one of the 
equivalent conditions in \ref{equ-reg} holds true.
\end{defn}

\begin{lemma}\label{reg_sur}
The image of a regular homomorphism $\phi:CH^n(X)_{\deg 0} \to
G$ is a connected algebraic subgroup of $G$.
\end{lemma}
\begin{proof}
Let $G'$ denote the Zariski closure of $\phi (CH^n(X)_{\deg
0})$. By \ref{pic}, $G'$ is connected and it is generated by the image of 
$\Pi_{X_\reg} \to G'$. Hence \ref{reg_sur} follows from the third
part of Lemma \ref{generators}.
\end{proof}

\section{The cycle class map}
Throughout the next three sections we will assume that the
ground field $k$ is the field of complex numbers.
$\sO_X$ and $\Omega^m_{X/\C}$ will respectively denote the sheaves of 
holomorphic functions and (analytic K\"ahler) differential $m$-forms.
As in the introduction consider the Deligne complex 
\[\sD(n)_X=\left(0\to
\Z_X(n)\to\sO_X\to\Omega^{1}_{X/\C}\to\cdots\to 
\Omega^{n-1}_{X/\C}\to 0\right),\]
and associated cohomology group $D^{n}(X)=\HH^{2n}(X,\sD(n)_X).$
In this section we construct the cycle class homomorphism
$CH^n(X) \to D^n(X)$, using Cartier curves $C$ in $X$.

By the moving lemma \ref{moving} it will be sufficient to
consider reduced Cartier curves $C$ in $X$.
Note, however, that we do not have that $C$ is a local complete
intersection in $X$, in general; this is only given to hold at points of
$C\cap X_{\sing}$. This leads to a slight technical difficulty. We will
need to define `Gysin' maps for Cartier curves $C$ in
$X$. These are directly defined in case $C$ is a local complete
intersection, and in general one has first to make a sequence of point
blow ups centered in $X_{\reg}$ to reduce to this special case. Indeed,
even to show that the cycle homomorphism $Z^n(X)\to D^n(X)$ respects
rational equivalence, a similar procedure needs to be followed.

Note that the exterior derivative yields a map of complexes
$\sD(n)_X\to \Omega^n_{X/\C}[-n],$ 
and there is an obvious map $\sD(n)_X\to\Z(n)_X$. 
\begin{lemma}\label{local}
For $x\in X_{\reg}$, there is a unique element
$[x]\in H^{2n}_{\{x\}}(X,\sD(n)_X)$
which maps to the topological cycle class of $x$ in
$H^{2n}_{\{x\}}(X,\Z(n))$ as well as to the ``Hodge cycle class''
of $x$ in $H^n_{\{x\}}(X,\Omega^n_{X/\C})$. 

This gives rise to a well-defined cycle class homomorphism
$Z^n(X)\to D^n(X)$, whose composition with $D^n(X)\to
H^{2n}(X,\Z(n))$ is the topological cycle class homomorphism.
\end{lemma}
\begin{proof}  The element $[x]$ exists because the topological and Hodge
cycle classes both map to the de Rham cycle class of $x$ in
$H^{2n}_{\{x\}}(X,\C)=\HH^{2n}_{\{x\}}(X,\Omega^{\d}_{X/\C})$, by a standard
local computation. See \cite{EV}, \S7, for example (though $X$ is singular, the
terms in the above computation depend only on a neighbourhood of $x$ in
$X$, and we have $x\in X_{\reg}$; hence \cite{EV}, \S7 is applicable).
\end{proof}
\begin{lemma}\label{curves1} 
In the above situation, if $\dim X=1$, then there is a natural
quasi-isomorphism $\sD(1)_X\cong\sO_X^*[-1]$, yielding an 
identification $\Pic(X)\cong D^1(X)$ (and hence also $\Pic^0(X)\cong
A^1(X)$). Under the identification, the class of a smooth point $[x]\in
D^1(X)$ corresponds to the class of the invertible sheaf $\sO_X(x)$. 
\end{lemma}
\begin{proof} The natural quasi-isomorphism is equivalent to the exactness
of the exponential sequence. The description of the class of a
point $x$ as the class of the invertible sheaf $\sO_X(x)$ is
also a standard local computation.
\end{proof} 

Now we argue as in \cite{BiS}, in order to show that the map
$Z^n(X)\to D^n(X)$ factors through $CH^n(X)$. We follow the convention
that the truncated de Rham complex of K\"ahler differentials
$$
\Omega_{X/\C}^{< n}=(0\to\sO_X\to\cdots\to\Omega^{n-1}_{X/\C}\to 0)
$$
has $\sO_X$ placed in degree
0; thus we have an exact sequence of complexes
\[0\>>> \Omega_{X/\C}^{< n}[-1]\>>> \sD(n)_X\>>> \Z(n)_X \>>> 0.\]
\begin{lemma}\label{curves2} Let $X$ be a projective variety of
dimension $n$ over $\C$, and $C\subset X$ be a reduced Cartier
curve which is a local complete intersection in $X$.  Then there
is a commutative diagram
\[
\begin{CD}
Z^1(C) \>>> Z^n(X)\\
\V V V \V V V\\
D^1(C) \> {\rm Gysin}>> D^n(X)\\
\V V V \V V V \\
H^2(C,\Z(1)) \>{\rm Gysin}>> H^{2n}(X,\Z(n))
\end{CD}
\]
\end{lemma}
\begin{proof}
Consider the local (hyper) cohomology sheaves $\sH^j_C(\sD(n)_X)$ of the
complex $\sD(n)_X$ with support in $C$. We claim that for any point $x\in
C$, the stalks $\sH^j_C(\sD(n)_X)_x$ vanish for $j\neq 2n-1$, unless $x$
is a singular point of $C$. Indeed, if $x\in C$ is a non-singular point
(so that $x\in X_{\reg}$ as well), then there is a long exact sequence of
stalks
\begin{multline*}
\cdots\>>>\sH^{j-1}_C(\Z(n)_X)_x\oplus
\sH^{j-1-n}_C(\Omega^n_{X/\C})_x\>>>
\sH^{j-1}_C(\C_X)_x\>>> \sH^j_C(\sD(n)_X)_x\\
\>>> \sH^j_C(\Z(n)_X)_x\oplus\sH^{j-n}_C(\Omega^n_{X/\C})_x\>>>
\sH^j_C(\C_X)_x\>>>\cdots
\end{multline*}
However $\sH^i_C(\Z(n)_X)_x=\sH^i_C(\C_X)_x=0$ for $i \neq2n-2$,
$\sH^{2n-2}_C(\Z(n)_X)$ injects into $\sH^{2n-2}_C(\C_X)$, and
$\sH^i_C(\Omega^n_{X/\C})_x=0$ unless $i=n-1$, for a non-singular point
$x\in C$ as above. This implies that $\sH^j_C(\sD(n)_X)_x=0$ for $j\neq
2n-1$, for such $x$. Also $\sH^{2n-1}_C(\sD(n)_X)_x$ fits into an exact
sequence 
\[0\>>>\sH^{2n-2}_C(\C_X/\Z(n)_X)_x\>>>\sH^{2n-1}_C(\sD(n)_X)_x\>>>
\sH^{n-1}_C(\Omega^n_{X/\C})_x\>>>
0,\]
with $\sH^{2n-2}_C(\C_X/\Z(n)_X)_x\cong\C/\Z(1)=\C^*$. 

Thus, $\sH^j_C(\sD(n)_X)$ is supported at a finite set of points, if
$j\neq 2n-1$. Hence in the local-to-global spectral sequence 
\[E_2^{p,q}=H^p(C,\sH^q_C(\sD(n)_X))\Longrightarrow
\HH^{p+q}_C(X,\sD(n)_X)\] 
we have $E_2^{p,q}=0$ for $p>0$, $q\neq
2n-1$. In particular, there is a well-defined injective map
\[\alpha:H^1(C,\sH^{2n-1}_C(\sD(n)_X))\>>> \HH^{2n}_C(X,\sD(n)_X).\]

We will next construct a natural map of sheaves on $C$
\[\sO_C^*\>>>\sH^{2n-1}_C(\sD(n)_X).\]
The desired Gysin map $D^1(C)\to D^n(X)$ is then defined to be the
composition
\begin{multline*}
H^1(C,\sO_C^*)\>>> H^1(C,\sH^{2n-1}_C(\sD(n)_X))\>{\alpha}>>
\HH^{2n}_C(X,\sD(n)_X)\>>>\\
\HH^{2n}(X,\sD(n)_X)=D^n(X)
\end{multline*}
To construct the map on sheaves $\sO^*_C\to
\sH^{2n-1}_C(\sD(n)_X)$, we argue locally, as follows. Let $U$ be
an affine neighbourhood in $X$ of a point $x\in C$, on which
the ideal of $C$ is generated by a regular sequence of functions
$f_1,\ldots,f_{n-1}$, determining a morphism $f:U\to \A^{n-1}_{\C}$ such
that $f^{-1}(0)=C\cap U$. Note that there are well-defined sections (of
the skyscraper sheaves) 
\[\alpha\in
\Gamma(\sH^{2n-2}_{\{0\}}(\Z(n)_{\A^{n-1}_{\C}}))=\Z(1),\;\;\;\;\beta\in 
\Gamma(\sH^{n-1}_{\{0\}}(\Omega^{n-1}_{\A^{n-1}_{\C}/\C})),\] 
which have the same image 
$\gamma\in \Gamma(\sH^{2n-2}_{\{0\}}(\Omega_{\A^{n-1}_\C/\C}^{\bullet})),$
under the obvious maps, and such that $\beta$ is annihilated by the ideal of
$0$ in $\Gamma(\sO_{\A^{n-1}_{\C}})$, for the natural module structure
on $\Gamma(\sH^{n-1}_{\{0\}}(\Omega^{n-1}_{\A^{n-1}_{\C}/\C}))$. In fact,
these conditions uniquely determine such a pair of sections
$(\alpha,\beta)$ up to sign, and there is a standard choice, with
$\beta$ determined by 
$\dlog(z_1)\wedge\cdots\wedge\dlog(z_{n-1}),$
where $z_j$ are the coordinate functions, so that $\beta$ is the cup
product of the local divisor classes
\begin{equation}\label{ext}
\dlog(z_j)\in
\Gamma(\A^{n-1}_{\C},\sext^1_{\A^{n-1}_\C}(\sO_{\{z_j=0\}},
\Omega^1_{\A^{n-1}_{\C}/\C})) 
\subset \Gamma(\A^{n-1}_{\C},
\sH^1_{\{z_j=0\}}(\Omega^1_{\A^{n-1}_{\C}/\C})).
\end{equation}
Hence $\gamma$ is also determined. Now consider
\begin{gather*}
f^*\alpha\in\Gamma(U,\sH^{2n-2}_C(\Z(n)_X)), \;\;\;\; 
f^*\beta\in \Gamma(U,\sH^{n-1}_C(\Omega^{n-1}_{X/\C})),\\
\mbox{and} \;\;\;\; f^*\gamma\in\Gamma(U,\sH^{2n-1}_C(\Omega^{<n}_{X/\C})),
\end{gather*}
where $f^*\alpha$ and $f^*\beta$ both map to $f^*\gamma$, and $f^*\beta$
is annihilated by any section of the ideal sheaf of $C\cap U$ in $U$. Thus
$f^*\alpha$ and $f^*\beta$ yield maps of sheaves
\[\Z(1)_C\mid_U\>>>
\sH^{2n-2}_C(\Z(n)_X)\mid_U,\;\;\;\;
\sO_C\mid_U\>>> \sH^{n-1}_C(\Omega^{n-1}_{X/\C})\mid_U,\]
giving rise to a commutative diagram of sheaves
$$
\begin{CD}
\Z(1)_C\mid_U \>>> \sH^{2n-2}_C(\Z(n)_X)\mid_U\\
\V V V \V V V\\
\sO_C\mid_U \>>> \sH^{2n-2}_C(\Omega^{<n}_{X/\C})\mid_U.
\end{CD}
$$
There is a long exact sequence of sheaves on $C$
$$
\cdots\>>>\sH^j_C(\Z(n)_X)\>>>\sH^j_C(\Omega^{<n}_{X/\C})\>>>
\sH^{j+1}_C(\sD(n)_X)\>>>\sH^{j+1}_C(\Z(n)_X)\>>>\cdots
$$
Hence from the exponential sequence 
\[0\>>> \Z(1)_C\>>> \sO_C\>{\rm exp}>>\sO_C^*\>>> 0,\]
and the above commutative diagram, we deduce that there is a well-defined
map of sheaves
\[\sO_C^*\mid_U\>>> \sH^{2n-1}_C(\sD(n)_X)\mid_U.\]

We will now show that these locally defined maps patch together to give 
well-defined sheaf maps
\begin{equation}\label{sheafmap}
\sO_C\>>> \sH^{2n-2}_C(\Omega^{n-1}_{X/\C}) \mbox{ \ \ \ and \ \
\ }\sO_C^*\>>> \sH^{2n-1}_C(\sD(n)_X).
\end{equation}
To do this, it suffices to show that the classes $f^*\alpha$,
$f^*\beta$ and $f^*\gamma$ defined above are in fact independent
of the map $f$, \ie of the choice of generators for the ideal of
$C$ in $U$. This too can be seen ``universally''.  Since the
ideal sheaf of $C$ in $X$ is locally generated by a regular
sequence, any two such sets of local generators for $\sI_C$ 
on the affine open set $U$ differ by the operation of an element of
$\GL_{n-1}(\sO_X(V))$, for some neighbourhood $V$ of $C\cap U$ in $U$.
Hence it suffices to show that if
$p:\GL_{n-1}(\C)\times\C^{n-1}\to\C^{n-1}$ is the projection, and
$m:\GL_{n-1}(\C)\times\C^{n-1}\to\C^{n-1}$ the map given by the operation
of $\GL_{n-1}(\C)$ on $\C^{n-1}$ by invertible linear transformations,
then $p^*\alpha=m^*\alpha$, $p^*\beta=m^*\beta$, 
and hence also $p^*\gamma=m^*\gamma$. We leave the verification of this to
the reader, as a simple application of the K\"unneth formula. 

Finally, note that for $U=X_{\reg}$, we have a commutative diagram with
exact rows
\[
\minCDarrowwidth=.8cm
\begin{CD}
0\>>> \Z(1)_C\mid_U \>>> \sO_C\mid_U\>>>\sO_C^*\mid_U\>>> 0\\ 
\noarr \V\cong V V \V V V \V V V\\
0\>>>
\sH^{2n-2}_C(\Z(n)_X)\mid_U\>>>\sH^{2n-2}_C(\Omega^{<n}_{X/\C})\mid_U
\>>> \sH^{2n-1}_C(\sD(n)_X)\mid_U\>>> 0 
\end{CD}\]
where the left vertical arrow is an isomorphism. 
For a smooth point $x\in C$, apply the functors
$\sH^j_{\{x\}}$ to the rows of the above diagram, and note that
$\sH^j_{\{x\}}(\Z(1)_C)=0$ for $j\neq 2$, and
$\sH^j_{\{x\}}(\sO_C)=0$ for $j\neq 1$.
We then obtain another diagram with exact rows 
\[
\minCDarrowwidth=.6cm
\hspace*{-.4cm}
\begin{CD}
0\>>> \sH^1_{\{x\}}(\sO_C) \>>> \sH^1_{\{x\}}(\sO_C^*) \>>>
\sH^2_{\{x\}}(\Z(1)_C)\>>> 0\\
\noarr \V V V \V V V \V V \cong V\\
0\>>> \sH^1_{\{x\}}(\sH^{2n-2}_C(\Omega^{<n}_{X/\C}))
\>>> \sH^1_{\{x\}}(\sH^{2n-1}_C(\sD(n)_X)) \>>> 
\sH^2_{\{x\}}(\sH^{2n-2}_C(\Z(n)_X))\>>> 0
\end{CD}\]
The bottom row may be identified (see \cite{Ha2}, III, Ex.~8.7, pg.~161) with
the exact sequence 
\[0\>>> \sH^{2n-1}_{\{x\}}(\Omega^{<n}_{X/\C})\>>> 
\sH^{2n}_{\{x\}}(\sD(n)_X)\>>> \sH^{2n}_{\{x\}}(\Z(n)_X)\>>> 0.\]
We claim that, under the above identification, the local cycle class of
$x$ in $\sH^2_{\{x\}}(\sD(1)_C)=\sH^1_{\{x\}}(\sO_C^*)$ maps to the
corresponding local cycle class of $x$ in $\sH^{2n}_{\{x\}}(\sD(n)_X)$.
Choosing a suitable regular system of parameters on $X$ at $x$, we reduce
to checking this in the special case when $x\in X$ is the origin $0\in
\C^n$, and the curve $C$ is the $z_n$-axis, given by the vanishing of the
first $n-1$ coordinates. We again leave this verification to the reader.

This means that, in the commutative diagram
\[
\begin{CD}
H^1_{\{x\}}(C,\sO_C^*) \>>> H^1(C,\sO_C^*)=D^1(C)\\
\V V V \V V V\\
\HH^{2n}_{\{x\}}(X,\sD(n)_X) \>>> D^n(X)
\end{CD}\]
the cycle class of $x$ in $D^1(C)$ maps to that of $x$ in $D^n(X)$. Hence
we have shown that there is a commutative diagram
\[\begin{CD}
Z^1(C) \>>> Z^n(X)\\
\V V V \V V V\\
D^1(C) \>{\rm Gysin}>> D^n(X)
\end{CD}
\]

It remains to show that the Gysin map $\Pic(C)=D^1(C)\to D^n(X)$ is
compatible with the topological Gysin map $H^2(C,\Z(1))\to
H^{2n}(X,\Z(n))$. Since $Z^1(C)\to D^1(C)$ is surjective, the
compatibility of the two Gysin maps is clear from the fact that each one
maps the class of $x$ on $C$ to the corresponding class on $X$.
\end{proof}
\begin{rmk}\label{gys_rem}
Assume that the local complete intersection curve $C$ lies in the
Cohen-Macaulay locus $X_{\CM}$ of $X$. Then the 
first sheaf map in (\ref{sheafmap}) factors as
\begin{equation}\label{gys_fac}
\sO_C\>>> \sext^{n-1}_X( \sO_C,\Omega^{n-1}_{X/\C}) \>>>
\sH^{2n-2}_C(\Omega^{n-1}_{X/\C}).
\end{equation}
To see this, note that
$$
\beta \in
\Gamma(\A^{n-1}_{\C},\sext^{n-1}_{\A^{n-1}_\C}(\sO_{\{0\}},
\Omega^{n-1}_{\A^{n-1}_{\C}/\C})) 
\subset \Gamma(\A^{n-1}_{\C},
\sH^{n-1}_{\{0\}}(\Omega^{n-1}_{\A^{n-1}_{\C}/\C}))
$$
as it is the product of the classes $\dlog(z_j)$ in (\ref{ext}).
Further the map $U \to \A^{n-1}_\C$ is flat in a neigbourhood
of $C\cap U$, as it is equidimensional. Thus $f^*\beta$ defines
a class in
$\Gamma(U,\sext^{n-1}_X(\sO_{C},\Omega^{n-1}_{X/\C}))$ 
mapping to $\Gamma(U,\sH^{n-1}_{C}(\Omega^{n-1}_{X/\C}))$.
As
\begin{multline*}
{\rm Ext}^{n-1}_{\GL_{n-1}(\C)\times
\C^{n-1}}(\sO_{\GL_{n-1}(\C)\times \{0\}},
\Omega^{n-1}_{\GL_{n-1}(\C)\times \C^{n-1}/\C})\\
\subset H^{n-1}_{\GL_{n-1}(\C)\times
\{0\}}(\GL_{n-1}(\C)\times \C^{n-1}, \Omega^{n-1}_{\GL_{n-1}(\C)\times
\C^{n-1}/\C})),
\end{multline*}
the class $f^*\beta$ defines the factorization (\ref{gys_fac}).
\end{rmk}
\begin{lemma}\label{blowup} Let $X$ be projective of dimension
$n$ over $\C$, $f:Y\to X$ the blow up of a smooth point $x\in
X$. Then the natural maps $f_*:CH^n(Y)\to CH^n(X)$ and
$f^*:D^n(X)\to D^n(Y)$ are isomorphisms, and there is a
commutative diagram
\[\begin{CD}
Z^n(Y) \>>> D^n(Y)\\
\V f_* V V \A \cong A f^* A\\
Z^n(X) \>>> D^n(X)
\end{CD}
\]
\end{lemma}
\begin{proof} The isomorphism on Chow groups is easy to prove, using the
fact that the exceptional divisor $E$ is a projective space (the
details are in \cite{BiS}). That $f^*:D^n(X)\to D^n(Y)$ is an
isomorphism is also easy to see, for the same reason, using also
the exact sequence
\[0\>>> f^*\Omega^1_{X/\C}\>>> \Omega^1_{Y/\C}\>>>
\Omega^1_{E/\C}\>>> 0.\] 
So we need to prove that if $y\in Y$ is any smooth point, then
its class in $D^n(Y)$ is the inverse image of that of $f(y)$ in
$D^n(X)$. This is clear if $f(y)\neq x$. If $f(y)=x$, we may
argue as follows. There is a commutative diagram with exact
rows
\[
\minCDarrowwidth=.6cm
\hspace*{-.5cm}
\begin{CD}
0 \>>> \HH^{2n}_{\{y\}}(Y,\sD(n)_Y)\>>>
H^{2n}_{\{y\}}(Y,\Z(n))\oplus H^n_{\{y\}}(Y,\Omega^n_{Y/\C})
\>>> H^{2n}_{\{y\}}(Y,\C(n))\>>> 0\\ 
\noarr \V V V \V V V \V V V\\
0 \>>> \HH^{2n}_{E}(Y,\sD(n)_Y)\>>>
H^{2n}_{E}(Y,\Z(n))\oplus H^n_{E}(Y,\Omega^n_{Y/\C})
\>>> H^{2n}_{E}(Y,\C(n))\>>> 0\\ 
\noarr \A f^* A A \A f^* A A \A f^* A A\\ 
0 \>>> \HH^{2n}_{\{x\}}(X,\sD(n)_X)\>>> 
H^{2n}_{\{x\}}(X,\Z(n))\oplus H^n_{\{x\}}(X,\Omega^n_{X/\C})
\>>> H^{2n}_{\{y\}}(X,\C(n))\>>> 0
\end{CD}\]
Here the downward vertical arrows are the natural maps
(``increase support'').  It is standard that the topological
local cycle classes of $x$ and $y$ have the same images in
$H^{2n}_E(Y,\Z(n))$. Similarly, the images in
$H^n_E(Y,\Omega^n_{Y/\C})$ of the local cycle classes of $x$
and $y$ in Hodge cohomology are also known to be equal; for
example, this follows from the existence of a Gysin map
$f_*:H^n_{\{y\}}(Y,\Omega^n_{Y/\C})\to
H^n_{\{x\}}(X,\Omega^n_{X/\C})$, which maps the local class of
$y$ to that of $x$, and which factors through
$$H^n_E(Y,\Omega^n_{Y/\C})\>{(f^*)^{-1}}>>H^n_{\{x\}}
(X,\Omega^n_{X/\C}).$$
Thus $f^*[x]=[y]\in 
\HH^{2n}_E(Y,\sD(n)_Y)$, and hence a similar equality is valid
in $D^n(Y)$ as claimed.
\end{proof}

\begin{lemma}\label{equivalence} The map $Z^n(X)\to D^n(X)$
factors through $CH^n(X)$, and hence determines a homomorphism
$\varphi:CH^n(X)_{\deg 0}\to A^n(X)$.
\end{lemma}
\begin{proof} This is similar to the corresponding proof in
\cite{BiS}. Let $C\subset X$ be a reduced Cartier curve, and $f\in
R(C,X)$. Let $\pi:Y\to X$ be a composition of blow ups at
smooth points so that the strict transform $\tilde{C}$ of $C$ in
$Y$ satisfies $\tilde{C}_{\sing}=\tilde{C}\cap Y_{\sing}\cong
C\cap X_{\sing}$. Then
\[R(C,X)=R(\tilde{C},Y)=R(\tilde{C},\tilde{C}),\]
and $\pi_*(f)_{\tilde{C}}=(f)_C\in Z^n(X).$ 
Now from lemma~\ref{curves1}, $(f)_{\tilde{C}}\mapsto 0\in D^n(Y),$
and so from lemma~\ref{blowup},
$(f)_C=\pi_*(f)_{\tilde{C}}=0\in D^n(X).$
\end{proof}
\begin{cor}\label{corblowup}
If $f:Y\to X$ is a composition of blow ups at smooth points,
then we have a diagram
\[\begin{CD}
CH^n(Y) \>>> D^n(Y)\\
\V f_* V \cong V \A \cong A f^* A\\
CH^n(X) \>>> D^n(X) 
\end{CD}
\]
\end{cor}
\begin{cor}\label{curves3}
For any reduced Cartier curve $C\subset X$, there are commutative
diagrams 
\[\begin{CD}
Z^1(C) \>>> Z^n(X)\\
\V V V \V V V\\
\Pic(C) \>{\rm Gysin}>>CH^n(X)\\
\V \cong V V \V V V\\
D^1(C) \>{\rm Gysin}>> D^n(X)\\
\V V V \V V V \\
H^2(C,\Z(1)) \>{\rm Gysin}>> H^{2n}(X,\Z(n))
\end{CD}
\mbox{ \ \ \ \ and \ \ \ \ \ }
\begin{CD}
\Pic^0(C) \>>> CH^n(X)_{\deg 0}\\
\V \cong V V \V V V \\
A^1(C) \>>> A^n(X)
\end{CD}\]
\end{cor}
\begin{proof} As in the proof of lemma~\ref{equivalence}, by a compositon
of blow-ups at smooth points, we reduce to the case when $C$ is a local
complete intersection in $X$. Then lemma~\ref{curves2} implies the
corollary. 
\end{proof}
Considering embedded resolution of singularities one obtains from
\ref{curves3} and \ref{corblowup} a second construction of the
Gysin map in \ref{gysin} over $\C$. At the same time, it gives the 
compatibility of this map with the Gysin map for the Deligne cohomology,
constructed in \ref{curves2}.
\section{Some general properties of $A^n(X)$ over $\C$}
It is shown in \cite{BiS} that if $X$ is projective over $\C$ of
dimension $n$, then there is a natural surjection (which is
referred to in \cite{BiS} as the {\em Abel-Jacobi map})
\[AJ^n_X:CH^n(X)_{\deg 0}\>>>
J^n(X):=\frac{H^{2n-1}(X,\C(n))}{F^0H^{2n-1}(X,\C(n))+{\rm
image}\,H^{2n-1}(X,\Z(n))},\] where by results of Deligne,
$J^n(X)$ is a semi-abelian variety (since the non-zero Hodge numbers of
$H^{2n-1}(X,\Z(n))$ lie in the set
$\{(-1,0),(0,-1),(-1,-1)\}$).
\begin{lemma}\label{algebraicity}
Let $X$ be projective of dimension $n$ over $\C$. Then there is
a natural surjection $\psi:A^n(X)\to J^n(X),$
whose kernel is a $\C$-vector space. $A^n(X)$ has a unique structure as an
algebraic group such that $\psi$ is a morphism of algebraic groups, with
additive kernel (\ie with kernel isomorphic to a direct sum of copies of
$\G_a$). 
\end{lemma}
\begin{proof}
By a result of Bloom and Herrera \cite{BH}, the natural map
\[H^{2n-1}(X,\C(n))\>>>
H^{2n-1}_{DR}(X/\C)=\HH^{2n-1}(X,\Omega^{\d}_{X/\C})\] 
is split injective. As explained in \cite{D}~(9.3.2), if $X_{\d}\to X$ is a
suitable hypercovering by a smooth proper simplicial scheme, the splitting
may be given by the composition 
$$H^{2n-1}_{DR}(X/\C)\to
H^{2n-1}_{DR}(X_{\d}/\C)\cong H^{2n-1}(X_{\d},\C(n))\cong
H^{2n-1}(X,\C(n)).
$$
{From} this description, the splitting is a map of
filtered vector spaces, where $H^{2n-1}(X,\C(n))$ has the
Hodge filtration for the mixed Hodge structure while
$H^{2n-1}_{DR}(X/\C)$ has the truncation filtration (\ie the
filtration b\^ete). 

Hence we obtain a commutative diagram
$$
\begin{TriCDV}
{H^{2n-1}(X,\C(n))}{\>>>}{H^{2n-1}(X,\C)/F^0H^{2n-1}(X,\C(n))}
{\SE E E}{\NE E {\vartheta}E}
{\HH^{2n-1}(X,\Omega^{< n}_{X/\C})}
\end{TriCDV}
$$
The map $\vartheta$ induces the map $\psi$ taking quotients modulo
$H^{2n-1}(X,\Z(n))$. Note that by weight considerations, the natural map
\[
H^{2n-1}(X,\Z(n))\>>> H^{2n-1}(X,\C(n))/F^0H^{2n-1}(X,\C(n))
\]
has a torsion kernel. Hence the kernels of $\psi$ and $\vartheta$ are the
same (and the latter is a $\C$-vector space).  This represents $A^n(X)$ as
an analytic group extension of the semi-abelian variety $J^n(X)$ by an
additive group $\G_a^r$, for some $r$, and hence as an analytic group
extension of an abelian variety by a group $\G_a^r\times \G_m^s$. As noted
in \cite{D}, (10.1.3.3), for any abelian variety $A$ over $\C$, the isomorphism
classes of analytic and algebraic groups extensions of $A$ by either
$\G_a$ or by $\G_m$ coincide (as a consequence of GAGA); hence a similar 
property is valid for extensions by $\G_a^r\times \G_m^s$.  This implies
that $A^n(X)$ has a unique algebraic structure such that $\psi$ is a
homomorphism of algebraic groups over $\C$, as claimed.
\end{proof}

For $X$ a curve, as in (i) of the next corollary, note that
$\Pic^0(X)$ has the natural algebraic structure obtained by
representing a suitable Picard functor. In 
particular, given an algebraic family of divisors (of degree 0)
on $X$ parametrized by a variety (or scheme) $T$, the induced
map $T\to\Pic^0(X)$ is automatically a morphism. On the other
hand, $A^1(X)$ has the algebraic structure given by
lemma~\ref{algebraicity}. Hence, a priori, the induced 
map $T\to A^1(X)$ obtained from such a family is only analytic,
since it is essentially given by integration. From (i) of the
corollary, it will follow that it is in fact algebraic. The content
of (ii) of the corollary is similar. 
\begin{cor}\label{algcurve}\ \ 
\begin{points} 
\item If $X$ is a curve, then the natural isomorphism $\Pic^0(X)\cong
A^1(X)$ is an isomorphism of algebraic groups.
\item In general, if $C\subset X$ is a reduced Cartier curve,
then the induced homomorphism $A^1(C)\to A^n(X)$ of
corollary~\ref{curves3} is algebraic. 
\end{points}
\end{cor}
\begin{proof} 
(i) The identification is certainly analytic, and in both cases,
when one represents the algebraic group as an extension of an
abelian variety by a commutative affine group, the abelian
variety in question is just
$\Pic^0(\tilde{X})=J(\tilde{X})=D^1(\tilde{X})$, the Jacobian of
the normalized curve $\tilde{X}$ (by which we mean the product of the
Jacobians of the connected components of $\tilde{X}$). Now one argues that
the identification must be algebraic as well, since one has the
one-one correspondence between analytic and algebraic extensions
of an abelian variety by $\G_a^r\times\G_m^s$.

(ii) Let $\tilde{X} \to X^{(n)}$ be a desingularization of $X^{(n)}$ such 
that the proper transform $\tilde{C}$ of $C$ is the normalization of
$C$. First note that one has a factorization 
$$
\begin{CD}
A^1(C) \>>> A^n(X)\\
\V V V \V V V\\
J^1(C) \>>> J^n(X)\\
\V V V \V V V\\
A^1(\tilde{C}) \>>> A^n(\tilde{X})
\end{CD}
$$
where all maps are analytic group homomorphisms, and
the vertical ones are algebraic (lemma \ref{algebraicity}). 
 Indeed the map $C \to X$ induces a morphism of mixed Hodge
structures $H^1(C) \to H^{2n-1}(X)$, and therefore an analytic
group homomorphim $J^1(C) \to J^n(X)$, which has to be algebraic
as it is compatible with its abelian part
$J^1(\tilde{C}) \to A^n(\tilde{X})$ and all analytic group
homomorphisms $\G_m^s \to \G_m^{s'}$ are algebraic. Similarly,
all group homomorphisms $\G_a^r \to \G_a^{r'}$ are algebraic,
and therefore $A^1(C) \to A^n(X)$ is algebraic as well.
\end{proof}
\begin{defn}\label{deflie}
For any commutative algebraic group $A$ over $\C$, let
$\Omega(A)$ denote the dual vector space to the Lie algebra
$\Lie(A)$. We may then identify $\Omega(A)$ with the vector
space of (closed) translation invariant regular 1-forms on $A$.
\end{defn}
Our next goal is to give a description of $\Omega(A^n(X))$, generalizing
the fact that for a non-singular projective variety $X$, $\Omega(A^n(X))$
is the space of holomorphic 1-forms on $X$ (since in that case, $A^n(X)$
is the Albanese variety of $X$). 
\begin{lemma}\label{purity}
Let $X$ be projective of dimension $n$ over $\C$, and let
$\omega_X$ denote the dualizing module of $X$ (in the sense of \cite{Ha}, 
Ch.~III, \S7). Let $X^{(n)}$ be the union of the $n$-dimensional 
components of $X$, and let $\omega_{X^{(n)}}$ denote its dualizing module. 
\begin{points} 
\item $\omega_X$ is annihilated by the ideal sheaf of $X^{(n)}$ in $X$. 
With its natural induced structure as an $\sO_{X^{(n)}}$-module, 
$\omega_X\cong\omega_{X^{(n)}}$, and is a torsion-free   
$\sO_{X^{(n)}}$-module. Hence for any coherent $\sO_X$-module $\sF$, the
sheaf $\shom_{\sO_X}(\sF, \omega_X)$ is also naturally an 
$\sO_{X^{(n)}}$-module, which is $\sO_{X^{(n)}}$-torsion free, and 
for any dense open set $U\subset X^{(n)}$, the restriction map 
\[\Hom_X(\sF,\omega_X)\>>> \Hom_U(\sF\mid_U,\omega_X\mid_U)\]
is injective. In particular, taking $U=X_{\reg}$, so that 
$\omega_X\mid_U=\Omega^n_{U/\C}$, and taking $\sF=\Omega^{n-1}_{X/\C}$, 
we have that  
$\Hom_X(\Omega^{n-1}_{X/\C},\omega_X)$
may be identified with a $\C$-subspace of the vector space of holomorphic
1-forms on $X_{\reg}$ which are meromorphic on $X^{(n)}$. 
\item $\Omega(A^{n}(X))$ is naturally identified with the subspace of 
$\Hom_X(\Omega^{n-1}_{X/\C},\omega_X)$ consisting of {\em closed} 1-forms.
\item When $n=1$,
\[\Omega(A^{1}(X))=\Omega(\Pic^{0}(X)))=H^{0}(X,\omega_X).\]
\item Let $j:X_{\CM}\to X$ be the inclusion of the open subset of
Cohen-Macaulay points. The natural map
\[\Omega(A^{n}(X))\>>>
\left(\mbox{closed 1-forms in 
$\Hom_X(\Omega^{n-1}_{X/\C},j^{m}_*j^{*}\omega_X)$}\right)\]
is an isomorphism, where $j^{m}_*$ denotes the meromorphic
direct image.
\end{points}
\end{lemma}
\begin{proof}  (i) We note first that $\omega_X\cong\omega_{X^{(n)}}$, and
the latter is a torsion-free $\sO_{X^{(n)}}$-module. Indeed, if we fix a
projective embedding $X\into\P^N_{\C}$, then
\[\omega_X=\sext_{\P^N}^{N-n}(\sO_X,\omega_{\P^N_{\C}}),\] 
and there is an analogous formula for $\omega_{X^{(n)}}$. As in \cite{Ha}, 
we see by Serre duality on $\P^N_{\C}$ that  
$\sext^i(\sF,\omega_{\P^N_{\C}})=0$ for all
$i\leq N-n$ for any coherent sheaf $\sF$ supported in dimension $<n$. This
gives the desired isomorphism, and implies that any local section of
$\sO_{X^{(n)}}$, which is a non zero-divisor, is also a non zero-divisor
on $\sext^{N-n}_{\P^N_{\C}}(\sO_X,\omega_{\P^N_{\C}})$. This means exactly
that $\omega_{X^{(n)}}$ is torsion-free.   

We conclude that for any coherent $\sO_X$-module $\sF$, the
sheaf $\shom_X(\sF,\omega_X)$ is a torsion-free $\sO_{X^{(n)}}$-module as
well. Applying this to $\sF=\Omega^{n-1}_{X/\C}$ gives (i).

(iii) is a special case of (ii). To prove (ii), first
note that from the definition of $A^{n}(X)$, we have 
\begin{equation}\label{lie_descr}
\Lie(A^{n}(X))=\coker(d:H^{n}(X,\Omega^{n-2}_{X/\C})\to
H^{n}(X,\Omega^{n-1}_{X/\C})).
\end{equation} 
{From} Serre duality for $H^{n}$ and $\Hom$, as in the definition of the
dualizing sheaf in \cite{Ha}, we have an identification of the dual vector
space
\[H^{n}(X,\Omega^i_{X/\C})^{*}=
\Hom_{\sO_X}(\Omega^i_{X/\C},\omega_X),\]
for any $i$. Thus $\Omega(A^n(X))$ is identified with the subspace of
$\Hom_{\sO_X}(\Omega^{n-1}_{X/\C},\omega_X)$ of elements $\varphi$ such that
the composition
\[\ell: H^n(X,\Omega^{n-2}_{X/\C})\>{d}>> H^n(X,\Omega^{n-1}_{X/\C})
\>{\varphi}>> H^n(X,\omega_X)\cong \C\]
is 0. It remains to show that, identifying elements $\varphi\in
\Hom_{\sO_X}(\Omega^{n-1}_{X/\C},\omega_X)$ with certain holomorphic
1-forms on $X_{\reg}$, $\Omega(A^n(X))$ is just the subspace of closed
1-forms.  

To see this, since we may consider $\varphi$ as a meromorphic 1-form on
$X$ which is holomorphic  on $X_{\reg}$, we can find a coherent sheaf of
ideals $\sJ$, defining the Zariski closed subset $X_{\sing}\subset X$ (\ie
the subscheme determined by $\sJ$ has $X_{\sing}$ as its underlying reduced
scheme), such that 
\begin{points}
\item $\eta\mapsto \eta\wedge\varphi$ defines an element of
$\Hom_{\sO_X}(\sJ\Omega^{n-2}_{X/\C},\Omega^{n-1}_{X/\C})$
\item $\eta\mapsto \eta\wedge d\varphi$ defines an element of 
$\Hom_{\sO_X}(\sJ\Omega^{n-2}_{X/\C},\omega_X)$, where we view
$\omega_X$ as a certain coherent extension of $\Omega^n_{X_{\reg}/\C}$ to
$X$.
\end{points} 
(Here $\sJ\sF$ denotes ${\rm image}\,(\sJ\tensor\sF\to\sF)$, for any ideal
sheaf $\sJ$ and coherent sheaf $\sF$). Since $\sJ$ defines $X_{\sing}$
within $X$, the natural map
\[H^n(X,\sJ\Omega^{n-2}_{X/\C})\>>> H^n(X,\Omega^{n-2}_{X/\C})\]
is surjective, and for any
$\varphi\in\Hom_{\sO_X}(\Omega^{n-1}_{X/\C},\omega_X)$, the composition 
\[\ell_1:H^n(X,\sJ\Omega^{n-2}_{X/\C})\>{d}>>
H^n(X,\Omega^{n-1}_{X/\C})\>{\varphi}>>  
H^n(X,\omega_X)\cong \C\]
factors through $\ell$. Thus
\[\varphi\in \Omega(A^n(X))\;\;\;\iff\;\;\; \ell_1=0.\]
We have 2 other related linear functionals 
\[\ell_2:H^n(X,\sJ\Omega^{n-2}_{X/\C})\>>>\C,\;\;\;\;
\ell_3:H^n(X,\sJ\Omega^{n-2}_{X/\C})\>>>\C,\] 
defined by
\begin{gather*}
\ell_2:H^n(X,\sJ\Omega^{n-2}_{X\C})\>{\wedge d\varphi}>>
H^n(X,\omega_X)\cong \C,\\
\ell_3: H^n(X,\sJ\Omega^{n-2}_{X/\C})\>{\wedge \varphi}>>
H^n(X,\Omega^{n-1}_{X/\C})\>{d}>> H^n(X,\omega_X)\cong\C,
\end{gather*}
where in the definition of $\ell_3$, we have let $d$ also denote the
composite of the exterior derivative
$\Omega^{n-1}_{X/\C}\to\Omega^n_{X/\C}$ with the natural map
$\Omega^n_{X/\C}\to\omega_X$. The formula 
\[d(\eta\wedge\varphi)=d\eta\wedge\varphi+(-1)^{n-2}\eta\wedge d\varphi,\]
for any $n-2$ form $\eta$, implies that
$\ell_3=\ell_1+(-1)^{n-2}\ell_2.$

Now by Serre duality and the $\sO_{X^{(n)}}$-torsion freeness of
$\shom_{\sO_X}(\sJ\Omega^{n-2}_{X/\C}, \omega_X)$ (see (i)), 
$\ell_2$ vanishes precisely when $d\varphi=0$ as a 2-form on $X_{reg}$. On
the other hand, we claim that for any $\varphi\in
\Hom_{\sO_X}(\Omega^{n-1}_{X/\C},\omega_X)$, the map $\ell_3$ constructed
as above is always 0. This will imply that $\ell_1=0$ \iff $\varphi$ is a
closed meromorphic 1-form. 

To prove that $\ell_3$ vanishes, it suffices to prove that the
map
\[H^n(X,\Omega^{n-1}_{X/\C})\>{d}>> H^n(X,\omega_X)\]
vanishes. One way to understand this is to note that if $\pi:Y\to X$ is a
resolution of singularities, then there is a commutative diagram
\[\begin{CD}
H^n(X,\Omega^{n-1}_{X/\C}) \>{d}>> H^n(X,\omega_X)\\
\V \pi^* VV \A A \pi_* A\\ 
H^n(Y,\Omega^{n-1}_{Y/\C}) \>{d}>> H^n(Y,\omega_Y)
\end{CD}
\]
which reduces us to proving that $H^n(Y,\Omega^{n-1}_{Y/\C})\>{d}>>
H^n(Y,\omega_Y)$ vanishes. This follows from Hodge theory, or alternately
may be proved as in \cite{Ha2}, III, lemma~8.4.

Proof of (iv):\quad We begin by recalling that since $X$ is
reduced, it is Cohen-Macaulay in codimension 1, so that
$Z=X-j(X_{\CM})$ has codimension $\geq 2$ in $X$.  Let $\sI$
denote the ideal sheaf of $Z$ in $X$. Let $\sD_m$ be the complex
of sheaves
\[\sD_m=(0\to j_!\Z_{X_{\CM}}(n)\to
\sI^{m+n-1}\by{d}\sI^{m+n-2}\Omega^{1}_{X/\C}\by{d} \cdots
\by{d} \sI^{m}\Omega^{n-1}_{X/\C}).\]
Then $\sD_m$ is a subcomplex of $\sD(n)_X$, whose cokernel
complex consists of sheaves supported on $Z$; the 0-th term of
the cokernel is $\Z(n)_Z$, while the other terms are coherent
sheaves supported on $Z$.  Since $\dim Z\leq n-2$, we see that
$\HH^{i}$ of this cokernel complex vanishes for $i\geq 2n-1$.
Hence $\HH^{2n}(X,\sD_m)\to
\HH^{2n}(X,\sD(n)_X)$ is an isomorphism, for all $m$.  Now as
in the proof of (i), one uses duality, to conclude that for all
$m$, there are isomorphisms
\begin{align*}
\Hom(\Omega^{n-1}_{X/\C},\omega_X)&\>>>
\Hom(\sI^m\Omega^{n-1}_{X/\C},\omega_X),\\ 
\Hom(\Omega^{n-2}_{X/\C},\omega_X)&\>>>
\Hom(\sI^m\Omega^{n-2}_{X/\C},\omega_X),
\end{align*} 
and taking the direct limit over all $m$, we obtain (iv).
\end{proof}

Our next goal is the proof of proposition~\ref{basic}, which gives us
another useful way to recognize elements of the vector space
$\Omega(A^n(X))$. We make use of two lemmas.

\begin{lemma}\label{basic1}
Let $X\subset\P^N_{\C}$ be a reduced projective variety of dimension $n$.
Then we can find a finite number of linear projections $\pi_i:X\to
\P^n_{\C}$, each of which is a finite morphism, such that the induced
sheaf map 
\[\bigoplus_i\pi_i^*\Omega^{n-1}_{\P^n/\C}\>>> \Omega^{n-1}_{X/\C}\]
is surjective.
\end{lemma}
\begin{proof} For any linear projection $\pi:X\to\P^n_{\C}$, there is a
factorization 
\[\pi^*\Omega^{n-1}_{\P^n_{\C}/\C}\>>>
\Omega^{n-1}_{\P^N_{\C}/\C}\tensor\sO_X\>{\psi}>>\Omega^{n-1}_{X/\C},\]
where the natural map $\psi$ is surjective.

So it suffices to prove the stronger assertion that there are projections
$\pi_i$ as above such that the induced sheaf map
\[\bigoplus_i\pi_i^*\Omega^{n-1}_{\P^n/\C}\>>>
\Omega^{n-1}_{\P^N_{\C}/\C}\tensor\sO_X\] 
is surjective.

We claim that for any $x\in X$, we can find a finite set of such
projections $\pi_i:X\to\P^n_{\C}$ such that the map of $\C$-vector spaces
\[\bigoplus_i\pi_i^*\Omega^{n-1}_{\P^n/\C}\tensor\C(\pi_i(x))\>>>
\Omega^{n-1}_{\P^N_{\C}/\C}\tensor\C(x)\]
is surjective. Indeed, the Grassmannian $\G_{\C}(n+1,N+1)$ (of $n+1$
dimensional subspaces of $\C^{N+1}$) parametrizes linear projections from
$\P^N_{\C}$ to $\P^n_{\C}$, and it contains a dense Zariski open subset
corresponding to projections which are finite morphisms on $X$. Hence the
$n$-dimensional vector subspaces 
\[\pi^*\Omega^1_{\P^n_{\C}/\C}\tensor\C(\pi(x))\subset
\Omega^1_{\P^N_{\C}/\C}\tensor\C(x)\] 
also range over a Zariski open subset of the Grassmannian of
$n$-dimensional subspaces of the cotangent space of $\P^N_{\C}$ at $x$.
In particular, we can find a finite number of them whose $(n-1)$-th exterior
powers span the $(n-1)$-th exterior power of this cotangent space, namely
$\Omega^{n-1}_{\P^N_{\C}/\C}\tensor\C(x)$.  

Now suppose $\pi_1,\ldots,\pi_r$ are chosen finite linear projections
$X\to\P^n_{\C}$, and that 
\[\bigoplus_{i=1}^r\pi_i^*\Omega^{n-1}_{\P^n/\C}\>>>
\Omega^{n-1}_{\P^N_{\C}/\C}\tensor\sO_X\] 
is not surjective. We can then find a point $x\in X$ at which the cokernel
is non-zero. By the above claim, we can augment the set of projections to
$\pi_1,\ldots,\pi_r,\pi_{r+1},\ldots,\pi_s$ so that the cokernel of the
new map
\[\bigoplus_{i=1}^{r+s}\pi_i^*\Omega^{n-1}_{\P^n/\C}\>>>
\Omega^{n-1}_{\P^N_{\C}/\C}\tensor\sO_X\]
does not have $x$ in its support. Thus the support of the cokernel has
strictly decreased. Now the lemma follows by Noetherian induction. 
\end{proof}

\begin{lemma}\label{basic2}
Let $\sF$ be a reflexive coherent sheaf on $\P^n_{\C}$, and $\omega$ a
meromorphic section of $\sF\tensor\Omega^1_{\P^n_{\C}/\C}$, which is
regular on some given (nonempty) Zariski open subset $W\subset\P^n_{\C}$. 
Suppose there is a non-empty open set $V$ in $\G_{\C}(2,n+1)$, the
Grassmannian of lines in $\P^n_{\C}$, such that 
\begin{points}
\item each line $L\in V$ meets $W$, and is disjoint from the non
locally-free locus of $\sF$
\item for each $L\in V$, the image of $\omega$ in
$(\sF\tensor\Omega^1_{L/\C})\mid_{L\cap W}$ extends to a regular section
of $\sF\tensor\Omega^1_{L/\C}$ on $L$. 
\end{points}
Then $\omega$ extends (uniquely) to a regular section on $\P^n_{\C}$ of
$\sF\tensor\Omega^1_{\P^n_{\C}/\C}$. 
\end{lemma}
\begin{proof}
Since $\sF$ is reflexive, it is locally free outside a Zariski closed 
set $A$ (of codimension $\geq 3$), and any section of
$\sF\tensor\Omega^1_{\P^n_{\C}/\C}$ defined in the complement of $A$
extends uniquely to a section on all of $\P^n_{\C}$. Since $\omega$ is a
meromorphic section, it determines a (unique) regular section of some twist 
$\sF\tensor\Omega^1_{\P^n_{\C}/\C}(D)$, for an effective divisor $D$;
there is a unique such twist $D$ which is minimal with respect to the
partial order on effective divisors (determined by inclusion of
subschemes). Our goal is to show that $D=0$. 

If $F$ is an irreducible component of $\supp D$ which appears
in $D$ with multiplicity $r>0$, then we can find a point $x\in F$ such that
\begin{points} 
\item $x$ is a non-singular point of $F$, and does not lie on any other
component of $D$; further, $\sF$ is locally free near $x$
\item $V$ contains a line through $x$
\item there is a regular parameter $t\in\sO_{x,\P^n_{\C}}$ (\ie $t$ is
part of a regular system of parameters) such that $t$ defines the
ideal of $F$ at $x$, and such that $t^r\omega$ determines a regular,
non-vanishing section of $\sF\tensor\Omega^1_{\P^n_{\C}/\C}$ in a
neighbourhood of $x$. 
\end{points}
It then follows that for a non-empty Zariski open set of lines $L$ through
$x$, we have $L\in V$, and $t^r \omega$ maps to a regular, non-vanishing
section of $\sF\tensor\Omega^1_{L/\C}$ near $x$, while $\omega$ itself
maps to a regular section of $\sF\tensor\Omega^1_{L/\C}$. However, $t$
vanishes at $x$. This is a contradiction.
\end{proof}

  If $C\subset X$ is a reduced, local complete intersection Cartier curve,
then in fact $C\subset X_{\CM}\cap X^{(n)}$ (recall that $X_{\CM}$ denotes
the (dense) Zariski open subset of Cohen-Macaulay points of $X$). The
sheaf map $\sO_C\to \sH^{2n-1}_C(\Omega^{n-1}_{X/\C})$ in (\ref{sheafmap})
induces a composite map
$$\alpha_C:H^1(C,\sO_C)\>>> H^1(C,\sH^{2n-1}_C(\Omega^{n-1}_{X/\C}))\>>>
H^n_C(X,\Omega^{n-1}_{X/\C})\onto \HH^{2n-1}(X,\Omega^{<n}_{X/\C}).$$
This is just the map $\Lie(\Pic^0(C))\to \Lie(A^n(X))$ on Lie algebras
induced by the composition of the group homomorphisms $\Pic^0(C)\to
A^1(C)$ and the Gysin map $A^1(C)\to A^n(X)$. 

\begin{propose}\label{basic} \ \
\begin{enumerate}
\item[(a)] Let $C\subset X$ be a reduced, local complete intersection
Cartier curve, and let $U\subset X_{reg}$ be a dense open subset such that
$U\cap C$ is dense in $C_{reg}$. Then the dual $\alpha_C^{\vee}$ of
$\alpha_C:H^1(C,\sO_C)\to \HH^{2n-1}(X,\Omega^{<n}_{X/\C})$  
(\ie of $\Lie(\Pic^0(C))\to \Lie(A^n(X))$) 
 fits into a commutative diagram
$$
\begin{CD}
\Omega(A^n(X)) \> \subset >> H^0(X_\reg, \Omega^1_{X_\reg/\C})\\
\V V {\alpha_C^{\vee}}V \V V {\rm restriction}V\\
H^0(C,\omega_C) \> \subset >>
H^0(C\cap U,\Omega^1_{C\cap U/\C}).
\end{CD}
$$
(Here the right hand vertical arrow is given by restriction of 1-forms.)
\item[(b)]
Let $U\subset X_{\reg}$ be a dense Zariski open set, and let $\omega\in
\Gamma(U,\Omega^1_{U/\C})$ be closed. Then $\omega\in
\Omega(A^n(X))$ if and only if
\begin{enumerate} 
\item[(i)] $\omega$ yields a meromorphic section on $X$ of $\Omega^1_{X/\C}$
\item[(ii)] for any reduced, local complete intersection Cartier curve
$C\subset X$ such that $C\cap U$ is dense in $C$, the
restriction of $\omega$ to $B=C_{reg}\cap U$ is in the image of the
natural injective map   
\[H^0(C,\omega_C) \>>> H^0(B,\Omega^1_{B/\C}).\]
\end{enumerate}
\end{enumerate}
\end{propose}
\begin{proof} 
First we prove (a). From lemma~\ref{purity}, it suffices to prove that if
$\beta_C:H^1(C,\sO_C)\to H^n(X,\Omega^{n-1}_{X/\C})$ is the obvious map
through which $\alpha_C$ factors, then the dual map $\beta_C^{\vee}$ fits into
a commutative diagram  
$$
\begin{CD}
H^0(X,\shom_X(\Omega^{n-1}_{X/\C},\omega_X)) \> \subset >> H^0(X_\reg,
\Omega^1_{X_\reg/\C})\\ 
\V V {\beta_C^{\vee}}V \V V {\rm restriction}V\\
H^0(C,\omega_C) \> \subset >>
H^0(C\cap X_\reg,\omega_{C\cap X_{\reg}}).
\end{CD}
$$ 
Here we have used Serre duality on $X$ and $C$ to make the identifications 
$$H^n(X,\Omega^{n-1}_{X/\C})^{\vee}=
H^0(X,\shom_X(\Omega^{n-1}_{X/\C},\omega_X)),$$
$$H^1(C,\sO_C)^{\vee}=H^0(C,\omega_C).$$

Since $C$ is a reduced, local complete intersection Cartier curve in $X$
(so that $C\subset X_{\CM}\cap X^{(n)}$), we have the adjunction formula 
\[\omega_C=\shom_C(\bigwedge^{n-1}\sI_C/\sI_C^2,\omega_X\tensor\sO_C).\] 
Hence there is a natural sheaf map
\[\psi_C:\shom_X(\Omega^{n-1}_{X/\C},\omega_X)\>>>
\shom_C(\bigwedge^{n-1}\sI_C/\sI_C^2,\omega_X\tensor\sO_C)=\omega_C\]
induced by restriction to $C$, and composition with the natural map
\begin{align*}
\bigwedge^{n-1}\sI_C/\sI_C^2&\>>>\Omega^{n-1}_{X/\C}\tensor\sO_C,\\  
f_1\wedge\cdots\wedge f_{n-1}&\longmapsto df_1\wedge\cdots\cdots
df_{n-1}. 
\end{align*} 
On any open set $U\subset X_{\reg}$ with $U\cap C\subset C_{\reg}$, one
verifies at once, from the explicit description, that the map
$\psi_C\mid_U$ is just the restriction map on 1-forms
$\Omega^1_{U/\C}\to\Omega^1_{C\cap U/\C}$. 

Hence the desired commutativity (which implies (a)) follows from:
\begin{claim}\label{duality}
$\beta^{\vee}$ is the map induced by $\psi_C$ on global sections. 
\end{claim}

To prove the claim, first note that for the local complete intersection
curve $C$ in $X_{\CM}$, one also has  
$$
\sext^{n-a}_X(\sO_C,\omega_X)=
\left\{ \begin{array}{ll}
\omega_C & \mbox{ \ for \ } a=1\\
0 & \mbox{ \ for \ } a\neq 1,
\end{array} \right.
$$
Hence there is a Gysin map given as the composite
\begin{multline*}
H^1(C,\omega_C)=H^1(X,\sext^{n-1}_X(\sO_C,\omega_X)) \> \epsilon
>> {\rm Ext}^n_X(\sO_C,\omega_X)\\
\>>> H^n_C(X,\omega_X) \>>> H^n(X,\omega_X) 
\end{multline*}
where $\epsilon$ is the isomorphism resulting from the (degenerate)
spectral sequence 
$$
E_2^{a,b-a}=H^a(X,\sext^{b-a}_X(\sO_C,\omega_X)) \Longrightarrow
{\rm Ext}^b_X(\sO_C,\omega_{X}).
$$
The trace map ${\rm Tr}_C:H^1(C,\omega_C) \to \C$ (of Serre duality on
$C$) factors as 
$$
{\rm Tr}_C:H^1(C,\omega_C) \>{\rm Gysin}>> H^n(X,\omega_X) \> {\rm Tr}_X >>
\C 
$$
(one way to verify this is to show that the composite ${\rm Tr}_X\circ{\rm
Gysin}$ has the universal property of ${\rm Tr}_C$).

Now the claim~\ref{duality} amounts to the assertion that the following
diagram commutes: 
\[\begin{CD}
H^1(C,\sO_C) \>{\rm Gysin}>> H^n(X,\Omega^{n-1}_{X/\C})\\
\V{\psi_C(\varphi)} V V \V{\varphi} V V\\
H^1(C,\omega_C) \>{\rm Gysin}>> H^n(X,\omega_X)
\end{CD}\]
{From} Remark~\ref{gys_rem}, this will follow if we prove the commutativity
of the diagram of $\sO_X$-linear maps
$$
\begin{CD}
\sO_C \>>> \sext^{n-1}_X(\sO_C,\Omega^{n-1}_{X/\C})\\
\V \psi_C(\varphi) V V \V V \varphi V\\
\omega_C \> \cong >> \sext^{n-1}_X(\sO_C,\omega_X).
\end{CD}
$$
As $\omega_C$ is torsion-free, it is enough to check this commutativity on
a suitable open subset of the regular locus of $C$, where it is easily
verified. 

We now show the ``if'' part of (b) (note that the other direction
follows directly from (a)). By lemma~\ref{basic1}, it suffices to prove
that for each finite, linear projection\\
$\pi:X\to\P^n_{\C}$, the meromorphic 1-form $\omega$ determines a
section of 
$$\shom_X(\pi^*\Omega^{n-1}_{\P^n_{\C}/\C},\omega_X).$$
Since $\pi$ is a finite morphism, 
\[\pi_*\omega_X=\shom_{\P^n_{\C}}(\pi_*\sO_X,\omega_{\P^n_{\C}}),\]
and we have a sequence of natural identifications of sheaves
\begin{align*}
&\pi_*\shom_X(\pi^*\Omega^{n-1}_{\P^n_{\C}/\C},\omega_X)\cong
\shom_{\P^n_{\C}}(\Omega^{n-1}_{\P^n_{\C}/\C}\tensor\pi_*\sO_X,
\omega_{\P^n_{\C}})\\
&\cong \shom_{\P^n_{\C}}(\pi_*\sO_X,
\shom_{\P^n_{\C}}(\Omega^{n-1}_{\P^n_{\C}/\C},\omega_{\P^n_{\C}}))
\cong \sF\tensor\Omega^1_{\P^n_{\C}/\C},
\end{align*}
where $\sF=\shom_{\P^n_{\C}}(\pi_*\sO_X,\sO_{\P^n_{\C}})$
is a (non-zero) coherent reflexive sheaf on $\P^n_{\C}$. 

Let $W\subset\P^n_{\C}$ be a dense open subset such that
$\pi^{-1}(W)\subset U$. Then $\omega$ determines a section of
$\sF\tensor\Omega^1_{\P^n_{\C}/\C}$ on $W$, and we want to show it
extends to a global section of this sheaf. We do this by verifying that
the hypotheses of lemma~\ref{basic2} are satisfied.

Let $L$ be a line in $\P^n_{\C}$, disjoint from the
non-flat locus of $\pi:X^{(n)}\to\P^n_{\C}$ (which is a subset of $\P^n_{\C}$
of codimension $\geq 2$, since $X^{(n)}$ is reduced and purely of dimension
$n$). Then the scheme-theoretic inverse image of $L$ in $X^{(n)}$ is a closed,
local complete intersection subscheme of $X^{(n)}$, purely of dimension 1, and
which is contained in the Cohen-Macaulay locus of $X^{(n)}$ (since
$X^{(n)}$ is Cohen-Macaulay precisely at all points $x\in X$ where $\pi$
is flat).  If further $L$ is not contained in the branch locus of $\pi$ on
$X^{(n)}$ (\ie $\pi$ is \'etale over all but finitely many points of $L$),
then $\pi^{-1}(L)=D$ is non-singular outside a finite set. Thus $D$ is a
reduced, complete intersection curve in $X^{(n)}$. Further, if $D\cap
X^{<n}=\emptyset$, then $D$ is a reduced local complete intersection curve
in $X$, whose non-singular locus is contained in $X_{\reg}$. In particular
$D$ is a reduced Cartier curve in $X$. Finally,  if $L$ is not contained
in the image of $X-U$, then $D$ has finite intersection with $X-U$, and
hence $D\cap U$ is dense in $D$. Clearly the set of all such lines $L$
contains a non-empty open subset of the Grassmannian of lines. 

For a line $L$ as above, we have
\begin{align*}
\pi_*\omega_D&\cong \shom_{L}(\pi_*\sO_D,\omega_L)\cong
\shom_L(\pi_*\sO_X\tensor\sO_L,\omega_L)\\
&\cong \shom_L(\pi_*\sO_X\tensor\sO_L,\Omega^1_{L/\C})
\cong \sF\tensor\Omega^1_{L/\C},
\end{align*}
since $\sF\tensor\sO_L\cong \shom_L(\pi_*\sO_D,\sO_L)$ (as $\pi$ is flat
over $L$). Since we are given that the image of $\omega$ in
$\Omega^1_{D\cap U/\C}$ extends to a global section of $\omega_D$, 
it follows that the corresponding section of
$\sF\tensor\Omega^1_{L/\C}\mid_{L\cap W}$ extends to a global section of
$\sF\tensor\Omega^1_{L/\C}$. Thus we have verified the hypotheses of
lemma~\ref{basic2}.
\end{proof}

\begin{rmk}\label{example}
Two properties of $A^n(Y)$, which are true for smooth projective varieties
$Y$, do not carry over to the general case: the compatibility with
products, and its dimension being constant in a flat family. We give
examples to illustrate these pathologies.

Let $X$ and $Y$ be projective varieties of dimension $n$ and
$m$, respectively, and let $r(X)$ and $r(Y)$ denote the number
of irreducible components of dimensions $n$ and $m$ respectively. By 
\cite{D} the K\"unneth decomposition $$
H^{2(n+m)-1}(X\times Y, \Z)/_{({\rm torsion})} = 
\left[H^{2n-1}(X,\Z)^{r(Y)} \oplus H^{2m-1}(Y,\Z)^{r(X)}\right]
/_{({\rm torsion})}
$$
is compatible with the Hodge structure. Thus 
\begin{equation}\label{product1}
J^{n+m}(X\times Y) = J^n(X)^{r(Y)} \times J^m(Y)^{
r(X)} .
\end{equation}
For $A^{n+m}(X \times Y)$ the picture is wilder.
By (\ref{lie_descr}) in the proof of \ref{purity}, we have 
\begin{gather}\label{product2}
\Lie(A^{n+m}(X\times Y))=\\
\frac{H^n(\Omega^{n-1}_{X/\C})\otimes
H^m(\Omega^{m}_{Y/\C})\oplus H^n(\Omega^{n}_{X/\C})\otimes
H^m(\Omega^{m-1}_{Y/\C})}{H^n(\Omega^{n-2}_{X/\C})\otimes
H^m(\Omega^{m}_{Y/\C})\oplus H^n(\Omega^{n-1}_{X/\C})\otimes
H^m(\Omega^{m-1}_{Y/\C})\oplus H^n(\Omega^{n}_{X/\C})\otimes
H^m(\Omega^{m-2}_{Y/\C})} \notag
\end{gather}
where the maps from the denominator are 
$$
d_X\otimes {\rm id}_Y, \ \ d_X\otimes {\rm id}_Y + (-1)^{n-1}{\rm id}_X
\otimes d_Y\mbox{ \ \ and \ \ }{\rm id}_X \otimes d_Y.
$$

Consider an elliptic curve $E$, the rational curve
$\Gamma=(x^3-y^2z) \subset \P^2_\C$, with a cusp, and the union
of three rational curves $ C = (xyz) \subset \P^2_\C$. They all
are fibres of the family $\sC \to \P=\P(H^0(\P^2,\sO_{\P^2}(3)))$
of curves of degree three in $\P^2$.

Hodge theory implies that 
$$
d:H^1(E,\sO_E) \>>> H^1(E,\Omega_E^1)\mbox{ \ \  and \ \ }
d:H^1(C,\sO_C) \>>> H^1(C,\Omega_C^1)\cong \C^3
$$
are both zero. Using (\ref{product2}) this shows
\begin{align*}
&A^2(C \times E) = J^2(C \times E) = 
A^1(C)\times A^1(E)^{3}= \G_m 
\times E^{3}\\
&A^2(C \times C) = J^2(C \times C) = 
A^1(C)^{3}\times A^1(C)^{3}= \G_m^6.
\end{align*}
On the other hand, $\Gamma - (0:1:0)= \Spec(\C[t^2,t^3])$ and,
if $\pi: \tilde{\Gamma} \to \Gamma$ denotes the normalization, one has
exact sequences
\begin{gather*}
0\>>> \sO_{\Gamma} \>>> \pi_*(\sO_{\tilde{\Gamma}}) \>>> \C t \>>> 0
\mbox{ \ \ and}\\
0\>>> \Omega^1_{\Gamma/\C} \>>> \pi_*(\Omega^1_{\tilde{\Gamma}/\C}) \>>>
\C dt \>>> 0. \hspace*{1cm}
\end{gather*}
Thus $\C t=H^1(\Gamma,\sO_\Gamma) \> d >\cong > \C dt \> > \subset >
H^1(\Gamma,\Omega^1_{\Gamma/\C}) \cong \C^2$
and one obtains by \ref{product2}
\begin{align*}
&A^2(\Gamma \times E) = \G_a \times E = A^1(\Gamma)\times A^1(E)\\
&A^2(\Gamma \times \Gamma) = \frac{\C^2 \times \C^2}{\C} =
\G_a^3 \mbox{ \ \ whereas \ \ } A^1(\Gamma)\times A^1(\Gamma) =
\G_a \times \G_a. 
\end{align*}
In particular, a product formula as (\ref{product1}) fails for $A^n$
instead of $J^n$, and the dimension of $J^n$ and $A^n$ are not
constant for the fibres $\sC \times \sC \to \P \times \P$.

It is amusing to write down the cycle map for the last example.
Writing
$$
\Gamma_\reg\times \Gamma_\reg = (\Gamma -(0:0:1))\times
(\Gamma-(0:0:1)) = \Spec(\C[u]\otimes_\C \C[v]),
$$
$\Omega(\Gamma\times\Gamma)={\rm
Hom}_{\Gamma\times\Gamma}(\Omega^1_{\Gamma\times\Gamma/\C}, 
\omega_{\Gamma\times\Gamma})_{\rm cl}$ decomposes as
\begin{align*}
&(H^0(\Gamma,\shom(\Omega^1_{\Gamma/|C},\omega_\Gamma)) \otimes
H^0(\Gamma,\omega_\Gamma) \oplus H^0(\Gamma,\omega_\Gamma)
\otimes H^0(\Gamma,\shom(\Omega^1_{\Gamma/\C},\omega_\Gamma)))_{\rm
cl}\\
&= (\C dv \oplus \C udv \oplus \C du \oplus \C vdu)_{\rm cl} =
\C dv \oplus \C du \oplus \C (udv+vdu).
\end{align*}
The cycle map is
\begin{align*}
\Pi_{\Gamma_\reg \times \Gamma_\reg}=\A_\C^2 \times \A_\C^2
& \>>> \G_a^3\\
((x_1,x_2),(y_1,y_2)) & \longmapsto \left\{
\begin{array}{rll}
du &\mapsto &y_1-x_1\\
dv &\mapsto &y_2-x_2\\
udv+vdu & \mapsto & y_1y_2-x_1x_2.
\end{array}\right.
\end{align*}
\end{rmk}
\section{The universal property over $\C$}     
Let $U_1,\ldots,U_r$ be the connected components of $X_{\reg}$,
and for each $i$, let $p_i\in U_i$ be a base point. 

Let $G$ be a commutative algebraic group.
By \ref{def-reg} and \ref{equ-reg} a homomorphism (of abstract groups)
$\phi:CH^n(X)_{\deg 0}\to G$ is regular, if and only if 
$\phi\circ\gamma_m:S^m(X_\reg)\to G$ is a morphism of varieties, for
some $m>0$. 

\begin{thm}\label{regular} \ \
\begin{points}
\item The homomorphism $\varphi:CH^n(X)_{\deg 0}\to A^n(X)$ constructed in
lemma~\ref{equivalence} is regular and surjective. 
\item The cokernel of the map $H_1(X_\reg,\Z) \to \Lie(A^n(X))$, defined by
integration of 1-forms over homology classes, is naturally
isomorphic to $A^n(X)$ and the composite 
$(\varphi\circ\gamma)^{(-)}: \Pi_{X_\reg} \to CH^n(X)_{\deg
0}\to A^n(X)$ is given by 
$$
(x,y) \longmapsto \bigr\{ \omega \mapsto \int^y_x \omega \bigr\}
$$
\item {\em (Universality)} $\varphi$ satisfies the following universal
property: for any regular homomorphism $\phi:CH^n(X)_{\deg 0}\to
G$ to a commutative algebraic group there exists a unique
homomorphism $h:A^n(X)\to G$ of algebraic groups with
$\phi=h\circ\varphi$.  
\end{points}
\end{thm}
\begin{proof}[Proof of (i)] It suffices to prove that
$\varphi\circ\gamma_1:U=X_\reg\to A^n(X)$ is a morphism. Note
that, from the definition, it is clearly analytic. Further, we
have the following. 
\begin{enumerate}
\item[(a)] the composition $U\to A^n(X)\to \Alb(\tilde{X})$ is a
morphism, where $\tilde{X}$ is a resolution of singularities of $X^{(n)}$,
since we may then regard $U$ as an open subset of $\tilde{X}$, and the map
$U\to \Alb(\tilde{X})$ is the restriction of the Albanese mapping for
$\tilde{X}$, with appropriate base-points. Here $\Alb(\tilde{X})$ is the
product of the Albanese varities of the connected components of
$\tilde{X}$, and $A^n(X)$ is an extension of $\Alb(\tilde{X})$ by a group
$\G_a^r\times\G_m^s$, so that in particular $A^n(X)\to\Alb(\tilde{X})$ is
a Zariski locally trivial fibre bundle.   
\item[(b)] For each reduced Cartier curve $C\subset X$, the composite
$$C_{\reg}\>>> U\>>> A^n(X)$$
is a morphism. Indeed, for each component $B_0$ of
$C_\reg$, the composition 
\[B_0\>>> U\>{\gamma_1}>>CH^n(X)_{\deg 0}\>{\varphi}>> A^n(X)\]
agrees with
\[B_0\>>> \Pic^0(C)\>>> CH^n(X)_{\deg 0}\>{\varphi}>>A^n(X)\]
up to a translation, and by corollary~\ref{algcurve}, the latter is
algebraic.  
\end{enumerate}
Now we may argue as in \cite{BiS}: we are reduced to proving that if $V$ is a
non-singular affine variety, a holomorphic function on $V$ which is
algebraic when restricted to ``almost all'' algebraic curves in $V$, is in
fact an algebraic regular function. This may be proved using Noether
normalization and power series expansions for holomorphic functions on
$\C^n$, or deduced from \cite{Si}, (1.1).

Since $\Omega(A^n(X))$ is a finite dimensional subspace of
1-forms on $U$, there exist reduced local complete intersection Cartier
curves $C_i \subset X$, for $i=1,\ldots,s$, such that 
$$
\Omega(A^n(X)) \>>> \bigoplus_{i=1}^s H^0(C_i,\omega_{C_i})
$$
is injective. Hence 
$$
\bigoplus_{i=1}^s \Pic^0(C_i) \> \oplus \psi_i >> A^n(X)
$$
is surjective.

\noindent
{\it Proof of (ii) and (iii):}\quad Let $\phi:CH^n(X)_{\deg
0}\to G$ be a regular homomorphism to a commutative algebraic group $G$. 
By lemma \ref{pic} the image of $\phi$ is contained in the connected
component of the identity of $G$. Hence we may assume without loss of
generality that $G$ is connected. 

Now $\Omega(G)$ consists of closed, translation-invariant 1-forms. Thus
if 
$$h=\phi\circ\gamma_1:U\>>> G,$$
then the image of $h^*:\Omega(G)\to
\Gamma(U,\Omega^1_{X/\C})$ is contained in the subspace of closed
1-forms. We claim that in fact 
$h^*(\Omega(G))\subset \Omega(A^n(X)).$

This is deduced from the criterion of proposition~\ref{basic},
(b), since we know that for any reduced Cartier curve $C$ in $X$, the
composition 
$$\Pic^0(C)\>>> CH^n(X)_{\deg 0}\>{\phi}>> G$$
is a homomorphism of algebraic groups. Now we observe that if
$B_0$ is any component of $C_{\reg}$, then 
$$
B_0 \>>> \Pic^0(C)=\Lie(A^1(C))/{\rm image}\,H_1(C_\reg,\Z)
$$
is given by integration of 1-forms in $H^0(C,\omega_C)$.
Moreover the composite 
$$
B_0\>>> U\>{h}>> G
$$
agrees with 
$$
B_0\>>>\Pic^0(C)\>>>G,
$$
up to a translation by an element of $G$ (and elements of
$\Omega(G)$ are translation invariant).  
Dualizing the above inclusion on 1-forms, we thus obtain a map
on Lie algebras 
$\Lie(A^n(X))\to\Lie(G).$
This fits into a commutative diagram
\[\begin{CD}
H_1(U,\Z)\>>> \Lie (A^n(X))\\
\V V V \V V V \\
H_1(G,\Z)\>>> \Lie(G)
\end{CD}\]
where the horizontal arrows are given by integration of 1-forms over
homology classes. Further, there is a commutative diagram
\begin{equation}\label{diag}
\begin{CD}
U \>>> \Lie(A^n(X))/{\rm image}\,H_1(U,\Z)\\ 
\V \gamma_1 V V \V V \tilde{\phi}V \\
CH^n(X)_{\deg 0}\> {\phi}>> G=\Lie(G)/{\rm image}\,H_1(G,\Z)
\end{CD}
\end{equation}
where $\tilde{\phi}$ is a homomorphism of analytic groups, and
where the upper horizontal arrow is given by integration of
1-forms in $\Omega(A^n(X))$.

We claim that the map $H_1(U,\Z)\to \Lie(A^n(X))=\Omega(A^n(X))^*$
factors through the (surjective) composition
\[H_1(U,\Z)\>{\cong}>> H^{2n-1}_c(U,\Z(n))\onto H^{2n-1}(X,\Z(n)),\]
where $H^*_c$ denotes compactly supported cohomology, and the isomorphism
is by Poincar\'e duality. Indeed, let $C\subset X^{(n)}$ be a sufficiently
general reduced complete intersection curve in $X^{(n)}$. Then 
\[C\cap X^{<n}=\emptyset,\;\; C_\sing=C\cap X^{(n)}_{sing}=C\cap 
X_{\sing},\] 
and one has a Gysin homomorphism 
\[H^1(C,\Z(1))\to H^{2n-1}(X,\Z(n))\cong H^{2n-1}(X^{(n)},\Z(n))\] 
which fits into a commutative diagram with exact rows 
\[\begin{CD}
H^0(C_{\sing},\Z(1)) \>>>  H^1_c(C\cap U,\Z(1)) \>>>
\hspace*{-.82cm}\to \hspace{.42cm} H^1(C,\Z(1))\hspace{.3cm}\\
\V V V \V V V \V V{\rm Gysin}V\\
H^{2n-2}(X^{(n)}_{\sing},\Z(n))\>>> H^{2n-1}_c(U,\Z(n))\>>>
\hspace*{-.7cm} \to \hspace{.2cm} H^{2n-1}(X^{(n)},\Z(n))
\end{CD}\]
The left hand vertical arrow is in fact surjective, since
$H^{2n-2}(X^{(n)}_{\sing},\Z(n-1))$ is the free abelian group on the
$(n-1)$-dimensional components of $X^{(n)}_{\sing}$, and (since $C$ is a
general complete intersection) $C_{\sing}$ has non-empty intersection
(which is supported at smooth points, and is transverse) with each such
component of $X^{(n)}_{\sing}$. Now we note that the composite $H_1(C\cap
U,\Z)\to H_1(U,\Z)\to \Lie(A^n(X))$ factors through the surjective
composite $H_1(U,\Z)\cong H^1_c(U,\Z(1))\to H^1(C,\Z(1))$, since $C\cap
U\to U\to A^n(X)$ is compatible with a homomorphism $\Pic^0(C)\to A^n(X)$
(here ``compatible'' means that for any component $B_0$ of $C\cap U$, the
composites $B_0\to U\to A^n(X)$ and $B_0\to\Pic^0(C)\to A^n(X)$ agree up
to translation by an element of $A^n(X)$). Now a diagram chase implies the
claim made at the beginning of the paragraph.

Thus in the diagram (\ref{diag}) we see that $\Lie(A^n(X))/{\rm
image}\,H_1(U,\Z)$ is identified with 
\[{\Lie(A^n(X))}/{{\rm image}\,H^{2n-1}(X,\Z(n))}\;=A^n(X).\]
Hence there is a homomorphism $\tilde{\phi}: A^n(X)\to G$, 
such that $\gamma_1 \circ \phi:U\to G$ factors through $A^n(X)$. 
Since $\gamma_1^*:\Omega(A^n(X))\to
\Gamma(U,\Omega^1_{U/\C})$ is injective, the induced map $A^n(X)\to G$
with this property is unique, since the corresponding map on Lie algebras
is uniquely determined. Since $\image\gamma_1$ generates $CH^n(X)_{\deg
0}$, the two homomorphisms 
\[\phi:CH^n(X)_{\deg 0}\>>> G,\;\; CH^n(X)_{\deg 0}\>{\varphi}>>
A^n(X)\>{\tilde{\phi}}>> G\]
must coincide. This proves the universal property of $\varphi$, except
that we need to note that $\tilde{\phi}$ is a morphism. 
By lemma \ref{equ-reg}, $\phi$ induces an algebraic group
homomorphism $\Pic^0(C) \to G$ for all admissible pairs $(C',
\iota)$, with $C = \iota(C')$. 
As above, we can choose reduced complete intersection curves $C_i$, $i=1,
\ldots, s$, such that 
$$
\bigoplus_{i=1}^s \Pic^0(C_i) \> \oplus \psi_i >> A^n(X)
$$
is surjective. As $\tilde{\phi} \circ (\oplus \psi_i)$ is an
algebraic group homomorphism, $\tilde{\phi}$ is an algebraic group
homomorphism as well. 
\end{proof}
\begin{rmk} Lemma~\ref{algebraicity}, combined with the Roitman Theorem
proved in [BiS], imply that $\varphi:CH^n(X)_{\deg 0}\to A^n(X)$ is an
isomorphism on torsion subgroups. In other words, the Roitman Theorem is
valid for $\varphi:CH^n(X)_{\deg 0}\to A^n(X)$, over $\C$. This is another
similarity with the Albanese mapping for a non-singular projective variety.
\end{rmk}
\begin{rmk}
The proof of theorem \ref{regular} is close in spirit to the
construction of a ``generalized Albanese variety'' in
\cite{FW}. There Faltings and W\"ustholz consider a finite
dimensional subspace $V \subset H^0(X_\reg,\Omega^1_{X_\reg})$,
containing the 1-forms with logarithmic poles on some
desingularization of $X^{(n)}$, and they construct a commutative algebraic
group $G_V$ together with a morphism $X_\reg \to G_V$, which is
universal among the morphisms $\tau: X_\reg \to H$ to
commutative algebraic groups $H$, with $\tau^*(\Omega(G))
\subset V$.
\end{rmk}
\section{Picard groups of Cartier curves}
In the next section, we give an algebraic construction of $A^n(X)$
for a reduced projective $n$-dimensional variety $X$, defined over an
algebraically closed field $k$. As in the analytic case, we will use the
Picard scheme for Cartier curves in $X$ and for families of such curves.
In this section, we discuss some properties of such families of curves, and
the corresponding Picard schemes. In particular, we establish the
technical results \ref{mu2} and \ref{mu3}, which are important steps in
the algebraic construction of $A^n(X)$.

Let $S$ be a non-singular variety, and let $f:\sC \to S$ be a flat
family of projective curves with reduced geometric fibres
$C_s=f^{-1}(s)$. Then
$$
g(C_s):= \dim_{k(s)} H^1(C_s,\sO_{C_s}) \mbox{ \ \ \ and \ \ \ }
\# C_s:= \dim_{k(s)} H^0(C_s,\sO_{C_s})
$$
are both constant on $S$. In fact, let 
$\sC \> \tilde{f} >> \tilde{S} \>\varkappa >> S$ 
be the Stein factorization of $f$. Since the fibres of $f$ are reduced,
$\varkappa:\tilde{S}\to S$ is a finite \'etale morphism and 
$\# C_s=\#\varkappa^{-1}(s)$ is constant, as well as 
$g(C_s)= \chi(C_s,\sO_{C_s}) - \# C_s$.

Let $S' \to S$ be a finite (possibly branched) covering such that
$f':\sC' = \sC \times_SS'\to S'$ is the disjoint union of families
of curves $\sC'_i \to S'$, for $i=1,\ldots ,s$ with connected
fibres. By \cite{BLR}, 8.3, theorem 1, the relative Picard
functors $\Pic_{\sC'_i/S'}$ are represented by an algebraic
space $\Pic(\sC'_i/S')$, and we define 
$$\Pic(\sC'/S')=\Pic(\sC'_i/S')\times_{S'}\cdots\times_{S"}\Pic(\sC'/S').
$$
For the smooth locus $\sC_\sm$ of $f$ consider the $g$-th
symmetric product 
$$
{f'}^g:S^g(\sC'_\sm /S')\>>> S'
$$
over $S'$. For any open subscheme $W'\subset S^g(\sC'_\sm /S')$
there is a natural map $\vartheta_{W'} : W' \to \Pic(\sC'/S')$.
By \cite{BLR}, 9.3, lemmas 5 and 6, one has the following
generalization of \ref{generators2}:
\begin{lemma}\label{relpic2} 
After replacing $S'$ by an \'etale covering, there exists
an open subscheme $W'\subset S^g(\sC'_\sm /S')$ with
geometrically connected fibres over $S'$,
such that $\vartheta_{W'} : W' \to \Pic(\sC'/S')$ is an open
embedding. 
\end{lemma}

Recall that $X^{(n)}$ denotes the union of the $n$-dimensional irreducible
components of $X$, and $X^{<n}$ is the union of the smaller dimensional
components. 
\begin{notations}\label{linearsystem}
For a very ample invertible sheaf $\sL$ on $X^{(n)}$ we write
$$
|\sL|^{n-1} = \P(H^0(X^{(n)},\sL))\times \cdots \times \P(H^0(X^{(n)},\sL))
\ \ \ \ (n-1)\mbox{-times}
$$
and $|\sL|^{n-1}_0$ for the open subscheme defined by $n-1$-tuples 
$D_1, \ldots , D_{n-1}$ of divisors such that 
\begin{enumerate}
\item[(i)] $C = D_1 \cap \cdots\cap D_{n-1}$ is a reduced complete 
intersection curve in $X^{(n)}$, 
\item[(ii)] $C\cap X^{<n}=\emptyset$, and 
\item[(iii)] $X_\reg \cap C$ is non-singular and dense in $C$. 
\end{enumerate}
Note that by (ii), $C$ is a reduced Cartier curve in $X$ which is
a local complete intersection. By abuse of notation we 
will sometimes 
write $C \in |\sL|^{n-1}$ instead of $(D_1, \ldots , D_{n-1}) \in
|\sL|^{n-1}$. 
\end{notations}

The normalization $\pi :\tilde{C} \to C$ induces a surjection
$\pi^*:\Pic^0(C) \to \Pic^0(\tilde{C})$. By \cite{BLR}, 9.2, the
kernel of $\pi^*$ is the largest linear subgroup $H(C)$ of
$\Pic(C)$. 
One has
\begin{align}\label{chieq}
\dim(H(C)) & = \dim(\Pic^0(C)) - \dim(\Pic^0(\tilde{C}))\notag\\
& = \dim_k(H^1(C,\sO_C)) - \dim_k(H^1(\tilde{C},\sO_{\tilde{C}}))\\ 
& = \chi(\tilde{C},\sO_{\tilde{C}}) - \chi({C},\sO_{{C}})
- ( \# \tilde{C} - \# C ),\notag
\end{align}
where again $\# C$ and $\# \tilde{C}$ denote the number of connected
components of $C$ and $\tilde{C}$, respectively.

Given a flat family of projective curves $f:\sC \to S$ over an irreducible
variety $S$ with reduced geometric fibres $C_s$, there exists a finite
(possibly branched) covering $S'\to S$ and an open dense subscheme
$S'_0\subset S'$ such that the normalization of $\sC\times_S S'_0$ is
smooth over $S'_0$. Hence $\# \tilde{C}_s$, and the dimension of the
linear part $H(C_s)$ of $\Pic^0(C_s)$, are both constant on the image of
$S'_0$. 
\begin{defn}\label{mu}
For a reduced projective curve $C$ we define $r(C)$ to be the
number of irreducible components of $C$, and $\mu(C)$ to be the
dimension of the largest linear subgroup of $\Pic^0(C)$.

By \cite{F}, Satz 5.2, for a very ample invertible sheaf $\sL$
the open subscheme $|\sL|^{n-1}_0$ is not empty. Then $r(\sL)$ and
$\mu(\sL)$ denote the values of $r(C)$ and of $\mu(C)$ for
$C \in |\sL|^{n-1}_0$ in general position.
\end{defn}
By the equality (\ref{chieq}) one has:
\begin{equation}\label{chieq2}
\chi(\tilde{C},\sO_{\tilde{C}}) - \chi({C},\sO_{{C}}) \geq
\mu(C) \geq \chi(\tilde{C},\sO_{\tilde{C}}) -
\chi({C},\sO_{{C}}) - r(C) + 1.
\end{equation}
\begin{lemma}\label{muinequality}
For a very ample invertible sheaf $\sL$ and for a positive integer $N$,
$$\mu(\sL^N) \leq N^{n-1}\cdot (\mu(\sL) + r(\sL) -1).$$
\end{lemma}
\begin{proof}
Given $D^{(i)}_j \in |\sL|$, for $i=1, \ldots , N$ and $j=1,
\ldots , n-1$, we write
\begin{align*}
I  &= \{1, \ldots ,N\}^{n-1}\\
C^{(\underline{i})} &= D^{(i_1)}_1 \cap \cdots \cap
D^{(i_{n-1})}_{n-1}\mbox{ \ \ for \ \ } \underline{i} = (i_1, \ldots ,
i_{n-1}) \in I\\
\mbox{and \ \ \ \ } C &= \bigcup_{\underline{i} \in I}
C^{(\underline{i})} = \bigcap_{j=1}^{n-1} (D_j^{(1)} \cup \cdots
\cup D_j^{(N)}). 
\end{align*}
\begin{claim}\label{general}
There exists a choice of the divisors $D^{(i)}_j \in |\sL|$ such
that
\begin{enumerate}
\item[(a)] $C^{(\underline{i})} \in |\sL|^{n-1}_0$, \ 
$\mu(C^{(\underline{i})})= \mu(\sL)$ \ and \
$r(C^{(\underline{i})})= r(\sL)$  \vspace{.05cm}
\item[(b)] $C^{(\underline{i})} \cap C^{(\underline{i}')} \cap
X_\sing = \emptyset$ for $\underline{i} \neq \underline{i}'$\vspace{.05cm}
\item[(c)] each point $x \in C_\sing \cap X_\reg$ lies on 
exactly two components $C^{(\underline{i})}$ and
$C^{(\underline{i}')}$. In this case, there exists one $\nu$ with
$i_j = i'_j$ for all $j\neq \nu$. Locally in $x$ the 
surface 
$$Y = D^{(i_1)}_1 \cap \cdots \cap \widehat{D^{(i_{\nu})}_\nu} \cap
\cdots \cap D^{(i_{n-1})}_{n-1}\cap X_\reg
$$
is nonsingular and contains $C^{(\underline{i})}$ and
$C^{(\underline{i}')}$ as two smooth divisors intersecting
transversally. \vspace{.05cm}
\item[(d)] $C$ is a reduced complete intersection curve in $|\sL^N|^{n-1}$.
\end{enumerate}
\end{claim}
\begin{proof} (d) follows from (a), (b) and (c).
Since $|\sL|^{n-1}_0$ is open and dense in $|\sL|^{n-1}$
(a) holds true for sufficiently general divisors. Counting
dimensions one finds that for $\underline{i}\neq \underline{i}'$
the intersection $C^{(\underline{i})} \cap C^{(\underline{i}')}$
is either empty or consists of finitely many points. The latter
can only happen, if all but one entry in $\underline{i}$ and
$\underline{i}'$ are the same, and obviously one may assume that
the intersection points all avoid $X_\sing$. Moreover
$$
C^{(\underline{i})} \cap C^{(\underline{i}')} \cap
C^{(\underline{i}'')}= \emptyset
$$
for pairwise different $\underline{i},\ \underline{i}'\
\underline{i}''\in I$. Now (c) follows from the Bertini theorem
\cite{F}, Satz 5.2, saying that for sufficiently general divisors
$D_j^{(i)}$ 
\begin{gather*}
Y = D^{(i_1)}_1 \cap \cdots \cap \widehat{D^{(i_{\nu})}_\nu} \cap
\cdots \cap D^{(i_{n-1})}_{n-1}\cap X_\reg\\
C^{(\underline{i})}=Y \cap D_\nu^{(i_\nu)}
\mbox{ \ \ \ and \ \ \ }C^{(\underline{i}')}=Y \cap D_\nu^{(i'_\nu)}
\end{gather*}
are non-singular and that $C^{(\underline{i})}$ and
$C^{(\underline{i}')}$ meet tranversally on $Y$.
\end{proof}
Let $\A^{M+1} \subset |\sL^N|^{n-1}$ be an affine open subspace
containing the point $s_0$ which corresponds to the tuple
$\{D_j^{(1)}\cup \cdots \cup D_j^{(N)}\}_{j=1,\ldots,n-1}$,
and let $\P^M$ be the projective space parametrizing lines in 
$\A^{M+1}$, passing through $s_0$. There is a line $S \in \P^M$ such that 
\begin{points}
\item
the total space $\sC$ of the restriction 
$$
\begin{CD}
\sC \>\tau>> X\\
\V f V V\\
S \> \subset >> |\sL^N|^{n-1}.
\end{CD}
$$
of the universal family to $S$ is non-singular in a
neighbourhood of each point $x\in C_\sing \cap X_\reg$
\item
the intersection of $X_\reg$ with the
general fibre of $f:\sC \to S$ is non-singular. 
\end{points}
In fact, using the notation from \ref{general} (c), we can
choose for a point $x\in C_\sing \cap X_\reg$ a line $S$
connecting $s_0$ with a point $(D'_1, 
\ldots ,D'_{n-1})\in \A^{M+1}$, where
$$
D'_j = D_j^{(1)}\cup \cdots \cup D_j^{(N)} \mbox{ \ \ for \ \ }
j \neq \nu,
$$
where $D'_\nu \cap Y_\reg$ is non singular, and where $x \not\in
D'_\nu$. By this choice, in a neigbourhood of $x$ the
restriction of the universal family $\sC$ to $S$ is just a
fibering of $Y$ over $S$. Hence the condition (i) is valid
for the chosen point $x$.

However,  for each point $x \in C_\sing \cap
X_\reg$, the condition (i) is an open condition in $\P^M$, and hence for a
general line $S$, (i) hold true for all points in $C_\sing \cap X_\reg$;
clearly the second condition (ii) holds as well.
The family $f:\sC \to S$ has only finitely many non-reduced
fibres and outside of them $U=\tau^{-1}(X_\reg)$ contains
only finitely many points, which are singularities of the fibres.

Replacing $S$ by an open neighbourhood of $s_0$, we may assume
thereby, that for $s\neq s_0$ the fibre $C_s=f^{-1}(s)$ is
reduced, that $C_s\cap X_\reg$ is non-singular and dense in $C_s$ and that
$\mu(C_s)=\mu(\sL^N)$. In particular $U$ is non-singular outside
of the points $C_\sing\cap X_\reg$, and by condition (i) $U$ is
non singular. Moreover, $f|_U : U \to
S$ is semi-stable; hence $f|_U$ is a local complete intersection morphism,
smooth outside a finite subset of $U$. Let $L$ be a finite extension of the
function field $k(S)$ such that the normalization of
$\sC\times_S \Spec(L)$ is smooth over $L$, and let $S'$ be the
normalization of $S$ in $L$. Consider 
$$
\begin{CDS}
\sC' \> \sigma >> \sC\times_S S' \> \eta' >> \sC\\
\novarr \SE E f' E \V pr_2 V V \novarr \V V f V \\
\noharr S' \> \eta >> S 
\end{CDS}
$$
where $\sigma$ denotes the normalization. Since
$U\times_S S' \to S'$ is a local complete intersection morphism, smooth
outside a finite subset of the domain, $U\times_S S'$ is normal and
$\sigma$ restricted to $U'=\sigma^{-1}(U\times_S S')$ is an
isomorphism. By construction the general fibre of $f'$ is smooth and
$\sC'$ is normal. Since for all $s' \in S'$ the fibres
$C'_{s'}={f'}^{-1}(s')$ of $f'$ are reduced on the open dense
subvariety $U'$, they are reduced everywhere. Note also that $\sC-U\to S$
is finite, and hence so is $\sC'-U'\to S'$.

Let $s', s'_0 \in S'$ be points, with $s_0=\eta(s'_0)$, and with
$s=\eta(s')$ in general position. The inequality (\ref{chieq2})
implies that
\begin{align*}
\mu(\sL^N) = \mu(C_{s}) &\leq \chi(C'_{s'},\sO_{C'_{s'}})-
\chi(C_s,\sO_{C_s})
\\ 
&= \chi(C'_{s'_0},\sO_{C'_{s'_0}})-\chi(C_{s_0},\sO_{C_{s_0}}).
\end{align*}
Since $C'_{s'_0}\cap U'$ is isomorphic to
$C_{s_0}\cap U$ the curve $C'_{s'_o}$ is
finite over and birational to $C=C_{s_0}$.
Moreover, the fibres $C'_{s'_0}\cap U'$ and $C\cap U$ have the same number
$\delta$ of double points. 
Writing ${C'}^{(\underline{i})}$ for the preimage of
${C}^{(\underline{i})}$ in $C'_{s'_0}$ one obtains
\begin{align*}
& \chi(C,\sO_{C}) + \delta = \sum_{\underline{i}\in I}
\chi({C}^{(\underline{i})},\sO_{{C}^{(\underline{i})}}),\\
& \chi(C'_{s'_0},\sO_{C'_{s'_0}})+ \delta = \sum_{\underline{i}\in I}
\chi({C'}^{(\underline{i})},\sO_{{C'}^{(\underline{i})}})\\
\mbox{and \ \ \ }
&\mu(\sL^N)\leq \sum_{\underline{i}\in
I}(\chi({C'}^{(\underline{i})},\sO_{{C'}^{(\underline{i})}})-
\chi({C}^{(\underline{i})},\sO_{{C}^{(\underline{i})}}) 
).
\end{align*}
Finally, ${C'}^{(\underline{i})}$ is finite over and birational
to ${C}^{(\underline{i})}$, thus it is dominated by the
normalization of ${C}^{(\underline{i})}$, and (\ref{chieq2})
implies 
\begin{align*}
\sum_{\underline{i}\in
I}(\chi({C'}^{(\underline{i})},\sO_{{C'}^{(\underline{i})}})-
\chi({C}^{(\underline{i})},\sO_{{C}^{(\underline{i})}}))
&\leq \sum_{\underline{i}\in I}(\mu({C}^{(\underline{i})})
+ r({C}^{(\underline{i})}) -1)\\
&=N^{n-1}\cdot (\mu(\sL) +r(\sL) -1).
\end{align*}
\end{proof}
Replacing $\sL$ by its $N$-th power one obtains by lemma
\ref{muinequality} ample invertible sheaves on $X^{(n)}$ with many more
linearly independent sections than $\mu(\sL)$. For example, if  
$X_1, \ldots ,X_r$ are the irreducible components of $X^{(n)}$, then
$$
\image( H^0(X^{(n)},\sL^N) \to H^0(X_i,\sL^N|_{X_i}))=
H^0(X_i,\sL^N|_{X_i}),
$$
for sufficiently large $N$, and its dimension i bounded below by a
non-zero multiple of $N^n$, whereas by \ref{muinequality}, $\mu(\sL^N)$
is bounded above by $(\mu(\sL)+r(\sL)-1)\cdot N^{n-1}$. One obtains:
\begin{cor}\label{mu2}
There exists a very ample sheaf $\sL$ on $X^{(n)}$ with
$$
\dim_k (\image( H^0(X^{(n)},\sL) \to H^0(X_i,\sL|_{X_i}))) \geq
2\cdot \mu(\sL)+r+2,
$$
for $i=1, \ldots, r$.
\end{cor}
Over a field $k$ of positive characteristic we will need
a stronger technical condition. Recall that
$\Pi_{X_\reg} = \bigcup_{i=1}^r (U_i\times U_i)$, where
$U_i=X_i \cap X_\reg$ are the irreducible components of
$X_\reg$. 
\begin{ass}\label{dominant}
Let $Z \subset S^d(\Pi_{X_\reg})\times
S^d(\Pi_{X_\reg})\times |\sL|^{n-1}_0$ be the incidence variety
of points
\begin{multline*}
(((x_1,x'_1), \ldots ,(x_d,x'_d)),((x_{d+1},x'_{d+1}),
\ldots ,(x_{2d},x'_{2d})), (D_1, \ldots ,D_{n-1}))\\
\in S^d(\Pi_{X_\reg})\times S^d(\Pi_{X_\reg})\times |\sL|^{n-1}_0
\end{multline*}
with 
$x_1, \ldots, x_{2d},x'_1, \ldots, x'_{2d} \in C=D_1 \cap \cdots
\cap D_{n-1}.$
Then the projection
$$
pr_{12}'=pr_{12}|_{Z}: Z \>>> S^d(\Pi_{X_\reg})\times
S^d(\Pi_{X_\reg})
$$
is dominant. 
\end{ass}
\begin{prop}\label{mu3}
There exists a very ample sheaf $\sL$ on $X^{(n)}$ which satisfies the
assumption \ref{dominant}, for all $d \leq \mu(\sL)$.
\end{prop}
\begin{proof}
Let $\sI_i$ be the ideal sheaf of $\bigcup_{j\neq i} X_j$ on
$X^{(n)}$. In particular $\sI_i|_{X_j}$ is zero, for $j\neq i$.
Hence if $\sF$ is a torsion-free coherent sheaf on $X^{(n)}$
then
\begin{align*}
&H^0(X_i, \sF \otimes \sI_i|_{X_i}/_{\rm
torsion})=H^0(X^{(n)},\sF\otimes \sI_i/_{\rm torsion})\\
\mbox{and \ \ \ }&\bigoplus_{i=1}^r H^0(X_i, \sF \otimes
\sI_i|_{X_i}/_{\rm torsion}) \subset H^0(X^{(n)},\sF). 
\end{align*}
\begin{claim}\label{mu4}
There exists a very ample invertible sheaf $\sL$ on $X^{(n)}$ such
that 
\begin{equation}\label{mueq}
\dim_k(\image(H^0(X_i, \sL \otimes \sI_i|_{X_i}) \to
H^0(C, \sL|_C))) \geq 4\cdot \mu(\sL),
\end{equation}
for $i = 1, \ldots , r$, and for all $C \in |\sL|^{n-1}_0$.
\end{claim}
\begin{proof}
Given a very ample invertible sheaf $\sL$ it suffices to find a lower
bound for the dimension of the image of the composite map
$$
\tau:
H^0(X_i, \sL^N \otimes \sI_i|_{X_i}) \>>> H^0(X^{(n)},\sL^N) \>>> H^0(C,
\sL^N|_C) $$
which is independent of $C\in |\sL^N|^{n-1}_0$ and grows like
$N^n$. If $\sJ_C$ denotes the ideal sheaf of $C$ on $X^{(n)}$, then
$$
\ker(\tau)= H^0\left(X_i, \sL^N \otimes \sJ_C \otimes
\sI_i|_{X_i}/_{({\rm torsion})}\right)
\subset H^0(X^{(n)}, \sL^N \otimes \sJ_C).
$$
Hence it is sufficient to give an upper bound for
$\dim(H^0(X^{(n)}, \sL^N \otimes \sJ_C))$ by some polynomial in $N$ of
degree $n-1$, independent of $C$. 

For $j < n$ the dimension of 
$H^{j}(X^{(n)}, \sL^{-N})$ is bounded by a polynomial of
degree $n-2$. In fact, $X^{(n)}$ is a subscheme of $\P^M=\P(H^0(X^{(n)},
\sL))$ and by \cite{Ha}, III.7.1 and III.6.9 one has, for $N$
sufficiently large,
\begin{align*}
H^j(X^{(n)}, \sL^{- N}) &\cong {\rm Ext}^{M-j}(\sO_{X^{(n)}}\otimes
\sL^{-N},\omega_{\P^M}) \\
&\cong H^0(\P^M,\sext^{M-j}(\sO_{X^{(n)}},\omega_{\P^M}\otimes
\sO_{\P^M}(N)))\\ 
&\cong H^0(\P^M,\sext^{M-j}(\sO_{X^{(n)}},\omega_{\P^M})\otimes \sO_{\P^M}(N)).
\end{align*}
Since $X$ is Cohen-Macaulay outside of a subscheme $T$ of
codimension $2$, the support of $\sext^{M-j}(\sO_{X^{(n)}},\omega_{\P^M})$
lies in $T$ for $M-j >M-n$.

The curve $C$ being a complete intersection of divisors in
$|\sL^N|$, a resolution of the ideal sheaf $\sJ_C$ on $X^{(n)}$ is given by
the Koszul complex 
$$
0 \to \sL^{-(n-1)N}=\bigwedge^{n-1}\bigl( \soplus{n-1} \sL^{-N} \bigr) \to
\ \ \cdots \ \ \to 
\bigwedge^{2}\bigl( \soplus{n-1} \sL^{-N} \bigr) \to
\soplus{n-1} \sL^{-N} \to \sJ_C \to 0.
$$
Therefore $\dim(H^0(X^{(n)}, \sL^N \otimes \sJ_C))$ is bounded from above
by 
$$\sum_{j=0}^{n-2}
\dim_k \bigl( H^j\bigl(X^{(n)},\sL^N\otimes \bigwedge^{j+1}\bigl( \soplus{n-1}
\sL^{-N} \bigr)\bigr)=\sum_{j=0}^{n-2}
\dim_k \bigl( H^j\bigl(X^{(n)},\sL^{-jN}\bigr)\bigr)
{\binom{n-1}{j+1}}. 
$$
\end{proof}
Let $\sL$ be a very ample invertible sheaf on $X^{(n)}$ which
satisfies the inequality (\ref{mueq}) in \ref{mu4}. We fix some
curve $C\in |\sL|^{n-1}_0$ and some natural number $d \leq
\mu(\sL)$. 

Each irreducible component of $S^d(\Pi_{X_\reg})\times
S^d(\Pi_{X_\reg})$ is of the form
$$
S_{\underline{d}}=(S^{d_1}(U_1 \times U_1) \times \cdots \times
S^{d_r}(U_r \times U_r)) \times (S^{d_{r+1}}(U_1 \times U_1)
\times \cdots \times S^{d_{2r}}(U_r \times U_r)),
$$
for some tuple $\underline{d}=(d_1, \ldots d_{2r})$ of non-negative
integers with 
$$
d_1 + \cdots + d_r = d_{r+1} + \cdots + d_{2r} =d.
$$
Given such a tuple ${\underline{d}}$, we claim that there are
(pairwise distinct) points $x_1, \ldots x_{2d},$ 
$x'_1 \ldots x'_{2d}$, with $x_\nu,x'_\nu \in C \cap U_i$, for
$$
\sum_{j=1}^{i-1} d_j < \nu \leq \sum_{j=1}^i d_j
\mbox{ \ \ \ and for \ \ \ }
\sum_{j=1}^{r+i-1} d_j < \nu \leq \sum_{j=1}^{r+i} d_j,
$$
such that the restriction map
\begin{equation}\label{restr}
H^0(X^{(n)},\sL) \>>> \bigoplus_{i=1}^{2d} k_{x_i}\oplus k_{x'_i}
\end{equation}
is surjective. In fact, by the inequality (\ref{mueq}) the
dimension of the image of 
$$
H^0(X_i, \sL \otimes \sI_i|_{X_i}) \>>>
H^0(C\cap X_i, \sL \otimes \sI_i|_{C\cap X_i}) \>>>
H^0(C, \sL|_C)))
$$
is at least  $4\cdot \mu(\sL) \geq 4 \cdot d \geq 2
\cdot (d_i + d_{r+i})$
and for sufficiently general points $x_1, \ldots x_{2d}, x'_1
\ldots x'_{2d} \in C$ the composite
$$
\bigoplus_{i=1}^rH^0(X_i,\sL\otimes \sI_i|_{X_i}) \subset
H^0(X^{(n)},\sL) \>>> \bigoplus_{i=1}^{2d} k_{x_i} \oplus k_{x'_i}
$$
is surjective.

By construction 
$$
w: =(((x_1,x'_1), \ldots ,(x_d,x'_d)),((x_{d+1},x'_{d+1}),
\ldots ,(x_{2d},x'_{2d}))) \in S_{\underline{d}}.
$$
Let $V$ denote the subspace of divisors $D \in |\sL|$
with
$$
x_1, \ldots ,x_{2d},x'_1,\ldots,x'_{2d} \in D.
$$ 
The fibre ${pr'}^{-1}_{12}(w)$ of the morphism
$pr_{12}': Z \to S^d(\Pi_{X_\reg})\times S^d(\Pi_{X_\reg})$
is the intersection of $V^{n-1}$ with $|\sL|^{n-1}_0$. 
In particular, since $(w,C) \in Z$, this intersection is
non-empty. 

If $\delta = \dim(|\sL|)$, the surjectivity of the restriction
map (\ref{restr}) implies that $\dim(V)= \delta - 4 \cdot d$
and $\dim({pr'}_{12}^{-1}(w)) = (n-1)\cdot (\delta - 4 \cdot d).$
The fibres of $pr_3':Z \to |\sL|^{n-1}_0$ are equidimensional of
dimension $4\cdot d$ and hence $Z$ is equidimensional of dimension
$(n-1)\cdot \delta + 4\cdot d$. Therefore the dimension of 
$pr_{12}'(Z) \cap S_{\underline{d}}$ can not be smaller than
$$
(n-1)\cdot \delta + 4\cdot d - (n-1)\cdot (\delta - 4 \cdot d)=
n \cdot 4 \cdot d = \dim(S_{\underline{d}}).
$$
\end{proof}
\section{The algebraic construction of $A^n(X)$}
Let $X$ be a projective variety of dimension $n$, defined over
an algebraically closed field $k$. 
As a first step towards the construction of $A^n(X)$ we need 
to bound the dimension of the image of a regular homomorphism 
$$\phi:CH^n(X)_{\deg 0} \>>> G
$$
to a smooth connected commutative algebraic group $G$.

By the theorem of Chevalley and Rosenlicht (theorems 1 and 2 in
\cite{BLR}, 9.2) there exists a unique smooth linear subgroup
$L$ of $G$ such that $G/L=A$ is an abelian variety. In addition,
$L$ is canonically isomorphic to a product of a unipotent group
and a torus. Let us write
$$
0\>>> L \>>> G \> \delta >> A \>>> 0
$$
for the extension.
\begin{lemma}\label{subgroups}
There exists a unique smooth connected algebraic subgroup $H$ of $G$, with
$\delta(H) = A$, such that every smooth connected algebraic subgroup $J$
of $G$ with $\delta(J)=A$ contains $H$. Moreover, the quotient
group $G/H$ is linear.
\end{lemma}
\begin{proof}
Given a smooth algebraic subgroup $J$ of $G$, one has the
commutative diagram of exact sequences
$$
\begin{CD}
\noarr 0 \noarr 0 \noarr 0 \\
\noarr \V V V \V V V \V V V\\
0 \>>> L\cap J \>>> J \>>> \delta(J) \>>> 0\\
\noarr \V V V \V V V \V V V\\
0 \>>> L \>>> G \>>> A \>>> 0\\
\noarr \V V V \V V V \V V V\\
0 \>>> L/(L\cap J) \>>> G/J \>>> A/\delta(J) \>>> 0\\
\noarr \V V V \V V V \V V V\\
\noarr 0 \noarr 0 \noarr 0 
\end{CD}
$$
Since $A/\delta(J)$ is an abelian variety and $L/(L\cap J)$
a linear algebraic group, $\delta(J)=A$ if and only if
$G/J$ is linear. Observe further, that $\delta(J) = A$ if and
only if $\delta(J') = A$ for the connected component $J'$ of $J$
containing the identity. 

Choose $H$ to be any smooth connected algebraic subgroup of $G$ with 
$\delta(H)=A$ and such that $\delta(H') \neq A$ for all proper
algebraic subgroups $H'$ of $H$. For $J$ as in \ref{subgroups}
consider the commutative diagram of exact sequences
$$
\begin{CD}
0 \>>> J\cap H \>>> G \>>> G/(J\cap H) \>>> 0\\
\noarr \V V V \V \Delta V V \V \iota V V\\
0 \>>> J\oplus H \>>> G\oplus G \>>> G/J \oplus G/H \>>> 0
\end{CD}
$$
where $\Delta$ is the diagonal embedding. Since
$J\cap H = \Delta^{-1}(J\oplus H)$ the morphism $\iota$
is injective on closed points, and hence $G/(J\cap H)$ is a linear
algebraic group. By the choice of $H$ one obtains $J\cap H = H$. 
\end{proof}
Recall that $X$ has $n$-dimensional irreducible components $X_1, \ldots
,X_r$, whose union is denoted $X^{(n)}$, and $U_i=X_{reg}\cap X_i$.  Also
$X^{<n}$ is the union of the lower dimensional components of $X$. 
\begin{prop}\label{bound}
Let $\sL$ be a very ample invertible sheaf on $X^{(n)}$ which
satisfies the assumption \ref{dominant}. Let
$g=\dim_k(H^1(C,\sO_C))$, for $C \in |\sL|^{n-1}_0$.

Let $\phi:CH^n(X)_{\deg 0} \to G$ be a surjective regular
homomorphism to a smooth connected commutative algebraic group $G$.
Then the induced morphism (see \ref{difference})
$$
\pi^{(-)}: S^{g+\nu\cdot \mu(\sL)}(\Pi_{X_\reg}) \>>>
CH^n(X)_{\deg 0} \>\phi >> G 
$$
is dominant, for $\nu > 0$, and surjective, for $\nu > 1$. 
In particular the dimension of $G$ is bounded by $2\cdot n \cdot
(g+\mu(\sL))$. 
\end{prop}
Probably the bound for the dimension of $G$ is far from being
optimal. We will indicate in \ref{bound2} how to obtain 
$\dim (G) \leq g$ in characteristic zero, under a weaker
assumption on $\sL$.
\begin{proof}[Proof of \ref{bound}]
Let again $L$ be the largest smooth linear algebraic subgroup
and $\delta:G \onto A=G/L$ the projective quotient group.
Recall that $|\sL|^{n-1}_0$ denotes the set of tuples
$(D_1, \ldots , D_{n-1})$ of divisors in the linear system
$|\sL|$ for which $C=D_1\cap \cdots \cap D_{n-1}$ is a reduced
complete intersection curve (in $X^{(n)}$), $C\cap X^{<n}=\emptyset$, and
$C\cap X_\reg$ non-singular and dense in $C$. 
\begin{claim}\label{surjective}
There exists an open dense subscheme $S \subset |\sL|^{n-1}_0$
such that 
$$
\Pic^0(C) \>\psi >> G \> \delta >> A
$$
is surjective and such that the dimension of $\image (\psi:
\Pic^0(C) \to G)$ is constant, for $C \in S$.
\end{claim}
\begin{proof}
Returning to the notation introduced in \ref{linearsystem}
let $S\subset |\sL|^{n-1}_0$ be an open subvariety, and let
$$
\begin{CD}
\sC \> \sigma >> X\\
\V f V V\\
S
\end{CD}
$$
denote the restriction of the universal complete intersection to
$S$. The smooth locus of $f$ is $\sC_{\sm}=\sigma^{-1}(X_\reg)$
and $\sC_\sm$ is dominant over $X_\reg$. 
Let $S' \to S$ be the finite morphism, and let
$W' \subset S^g(\sC_\sm\times_S S' /S')$
be the open connected subscheeme considered in lemma
\ref{relpic2}. Replacing $S$ and $S'$ by some open subschemes, and $S'$ by
a branched cover if necessary, we may assume that there exists a section
$\epsilon'$ of $W' \to S'$. By \ref{relpic2}
\begin{align*}
W' &\>>> \Pic(\sC/S)\times_S S'= \Pic(\sC\times_S S'/S')\\
w' &\longmapsto \vartheta_{W'}(w') - \vartheta_{W'}(\epsilon'(pr_2(w')))
\end{align*}
is an open embedding. On the other hand, one has a morphism
of schemes
$$
h: W' \>>> S^g(\sC_\sm\times_S S' /S') \>>> S^g(X_\reg),
$$
and the image of the connected scheme $W'$ lies in some
connected component, say
$S_{\underline{g}}=S^{g_1}(U_1) \times \cdots \times S^{g_r}(U_r).$
Since $\phi: CH^n(X)_{\deg 0} \to G$ is regular, the composite
$$
h^{(-)}:W'\times_{S'} W' \>>> S_{\underline{g}}\times S_{\underline{g}} \>
\theta >> G
$$
is a morphism, where
$$
\theta(\underline{x},\underline{x}')= \phi\bigl(
\sum_{i=1}^g \gamma(x_i) - \sum_{i=1}^g \gamma(x'_i)\bigr) .
$$ 
The morphism $h^{(-)}$ induces $S'$ morphisms
$$
h_{S'}^{(-)}: W'\times _{S'} W' \>>> G\times S' \mbox{ \ \ and \ \ }
h_{S'}^{(-)}\circ\delta:W'\times_{S'} W' \>>> A \times S'. 
$$
Let $W_G$ and $W_A$ be locally closed subschemes of the images
$$
h_{S'}^{(-)}(W'\times_{S'} W')\mbox{ \ \ and \ \ }
h_{S'}^{(-)}\circ\delta(W'\times_{S'} W')
$$
respectively, dense in the closure of the images.
Choosing $S'$ and $S$ small enough, one may assume that 
$S' \to S$ is surjective and that $W_G$ and $W_A$ 
are both equidimensional over $S'$. For $C\in S$ choose a point
$s' \in S'$ mapping to 
$C \in S$ and let $W'_{s'}$ denote the fibre of $W'$ over $s'$.
Then the image of $W'_{s'} \times W'_{s'}$ in $\Pic^0(C)$ is dense and
thereby $\dim( \psi(\Pic^0(C)))$ and
$\dim(\delta(\psi(\Pic^0(C)))) = d'$ are both constant on $S$. 

Asume that $d' < \dim(A)$. The closure of
$\delta(h^{(-)}(W'_C\times W'_C))$ is the image of $\Pic^0(C)$, 
hence $\delta(h^{(-)}(W'_C\times W'_C))$ lies in some abelian
subvariety $B$ of $A$ of dimension $d'< \dim(A)$. Since $S'$ and $W'$ are
connected, and since an abelian variety $A$ does not contain
non-trivial families of abelian subvarieties, $B$ is independent
of the curve $C$ chosen. 

$\sC_\sm$ being dominant over $X_\reg$ this implies
that the image $\delta\phi(CH^n(X)_{\deg 0})$ lies in $B$,
contradicting the assumptions made.
\end{proof}
In general, a commutative algebraic group $G$ can contain non-trivial
families of subgroups and the argument used above does not
extend to $G$ instead of $A$.

Let $H \subset G$ be the smallest connected algebraic subgroup
with $\delta(H) = A$, 
as constructed in \ref{subgroups}. By \ref{surjective} and by
the universal property in \ref{subgroups}, for $C \in S$ the
image of $\psi(\Pic^0(C))$ contains $H$. 

By \ref{generators} the image
of $S^g(\Pi_{C_\reg})$ in $G$ is $\psi(\Pic^0(C))$ and hence $H$
is contained in the image of $S^g(\Pi_{X_\reg})$. In order to
show that 
$$
\pi^{(-)}: S^{g+\mu(\sL)}(\Pi_{X_\reg}) \>>>
CH^n(X)_{\deg 0} \>\phi >> G 
$$
is dominant, it suffices to verify that the
image $Y_0$ of the composite 
$$
\tau^{(-)}: S^{d}(\Pi_{X_\reg}) \>>>
CH^n(X)_{\deg 0} \>\phi >> G \>>> G/H
$$
is dense, for some $d \leq \mu(\sL)$. Applying claim
\ref{surjective} to $G/H$ instead of $G$ one finds a non-empty
open subscheme $S\subset |\sL|^{n-1}_0$ such that the dimension $d$
of $\psi'(\Pic^0(C))$ is constant on $S$, where $\psi':\Pic^0(C) \to
G/H$ is the natural map (see \ref{pic}). Since $G/H$ is a linear
algebraic group, we must have $d\leq \mu(\sL)$, and
choosing $S$ small enough, we may assume that
\begin{equation}\label{bd}
d=\dim(\psi'(\Pic^0(C))) \leq \mu(C) = \mu(\sL), \mbox{ \ \ \ \
for all \ $C \in S$.} 
\end{equation}

Since $Y_0$ generates the group $G/H$, it is dense in $G/H$ 
if and only if its closure $Y$ is a group. By assumption
the image of the incidence variety 
$$
Z \> pr'_{12} >>S^{d}(\Pi_{X_\reg})\times S^{d}(\Pi_{X_\reg})
\>\tau^{(-)}\times \tau^{(-)}>> Y \times Y
$$
defined in \ref{dominant} contains some open dense subscheme
$T$. By definition, for each $t \in T$ there exist divisors
$D_1, \ldots, D_{n-1}$ with 
$C=D_1\cap \cdots \cap D_{n-1} \in S$ and with
$$
t \in \image( S^{d}(\Pi_{B})\times S^{d}(\Pi_{B})
\>{\vartheta}^{(-)}\times {\vartheta}^{(-)}>> Y \times Y)
$$
for $B = C\cap X_\reg$ and for the induced map $\vartheta^{(-)}$
from $\Pi_B$ to $G/H$.
By \ref{generators}
$\psi'(\Pic^0(C))={\vartheta}^{(-)}(S^d(\Pi_B))\subset Y$. Since
$\psi'(\Pic^0(C))$ is an algebraic subgroup of $G/H$, the image of
$t$ under the morphism 
$$
{\rm diff}:G/H \times G/H \>>> G/H \mbox{ \ \ \ with \ \ \ } (g,g')
\longmapsto g-g'
$$
is contained in $\psi'(\Pic^0(C))$, hence in $Y$.

Thereby $T$ is a subset of ${\rm diff}^{-1}(Y)$, and the same is
true of its closure $Y \times Y$. One obtains that
${\rm diff} (Y \times Y) \subset Y$ and $Y$ is a subgroup of
$G/H$. 

Since $Y_0$ is dense in $G/H$, by lemma \ref{generators} (ii)
the image of $S^{2d}(\Pi_{X_\reg})$ is $G/H$.
\end{proof}
As indicated already, the proposition \ref{bound} can be
improved in characteristic zero.
\begin{variant}\label{bound2} Assume that $\Char (k) =0$.
Let $\sL$ be a very ample invertible sheaf on $X^{(n)}$ with
\begin{equation}\label{eqmu2}
\dim_k (\image( H^0(X^{(n)},\sL) \to H^0(X_i,\sL|_{X_i}))) \geq
2\cdot \mu(\sL) + r + 2,
\end{equation}
for $i=1,\ldots ,r$. Let $G$ be a smooth connected commutative
algebraic group, and let $\phi:CH^n(X)_{\deg 0} \to G$ be a
surjective regular homomorphism. Then there exists
an open dense subvariety $S \subset |\sL|^{n-1}_0$ such that for
each $C \in S$ the induced homomorphism (see \ref{pic})
$$
\psi: \Pic^0(C) \>>> CH^n(X)_{\deg 0} \>\phi >> G
$$
is surjective. In particular the dimension of $G$ is bounded by 
$g=\dim_k H^1(C,\sO_C)$.
\end{variant}
\begin{proof}
The first part of the proof is the same as the one for
\ref{bound}. In particular we may assume claim \ref{surjective}
to hold true. 

Let $H \subset G$ be the smallest subgroup with $\delta(H) = A$,
as constructed in \ref{subgroups}. By \ref{surjective} and by
the universal property in \ref{subgroups}, for all $C \in S$ the
image of $\psi(\Pic^0(C))$ contains $H$. Hence $\psi:\Pic^0(C)
\>>> G$ is surjective if and only if
$$
\Pic^0(C) \> \psi >> G \>>> G/H
$$
is surjective. In order to prove \ref{bound2} we may assume
thereby that $G$ is linear and $A=0$.

For $C \in S$, let $\gamma_B : B=C\cap X_\reg \to CH^n(X)$
denote the natural map and let $\Gamma_B$ be the image of the
composite 
$$
\vartheta^{(-)}:\Pi_B \> \gamma_B^{(-)} >> CH^n(X)_{\deg 0} \> \phi >> G.
$$
For any subset $M \subset G$ we will denote by $G(M)$ the smallest
algebraic subgroup of $G$ which contains $M$. If $M$ contains a
point of infinite order, then $\dim(G(M)) > 0$. In
characteristic zero the converse holds true, as well.
In fact, if $\dim(G(M)) > 0$ then $G(M)$ contains a subgroup
isomorphic either to $\G_a$ or to $\G_m$. In characteristic
zero, both contain points of infinite order.

Hence if the dimension of $G(\Gamma_B) = \psi(\Pic^0(C))$
is larger than zero, the constructible set
$\Gamma_B$ contains a point $\alpha_1$ of infinite order and
$\dim(G(\alpha_1)) > 0$. Repeating 
this for $G/G(\alpha_1,\ldots,\alpha_\nu)$ instead of
$G$, we find recursively points $\alpha_1, \ldots, \alpha_d\in
\Gamma_B$ with $G(\Gamma_B)=G(\alpha_1,\ldots, \alpha_d).$

Let us choose points $x_1, \ldots, x_d,x'_1, \ldots, x'_d \in B$
with $\alpha_j=\vartheta^{(-)}((x_j,x'_j))$, and moreover,
for each component $X_i$ of $X^{(n)}$, choose a base point $q_i \in B\cap 
X_i$. 
\begin{claim}\label{dominant2}
There exists a closed suscheme $Z\subset S$ such that the restriction
$$
\begin{CD}
\sC'=\sC\times_S Z \> \sigma'=\sigma|_{\sC'} >> X^{(n)}\>\subset>> X\\
\V f'=f|_{\sC'} V V\\
Z \>\subset >> S\> \subset >>|\sL|^{n-1}_0
\end{CD}
$$
of the universal family satisfies:
\begin{enumerate}
\item[(a)] For each point $z \in Z$ the curve $C_z={f'}^{-1}(z)$
contains the points 
$$x_1, \ldots, x_d,x'_1, \ldots, x'_d, q_1,\ldots , q_r.$$
\item[(b)] $\sigma': \sC' \to X^{(n)}$ is dominant.
\end{enumerate}
\end{claim}
\begin{proof}
For 
$$
V_i=(\image(H^0(X^{(n)},\sL) \to H^0(X_i,\sL|_{X_i}))-0)/k^* \subset
|(\sL|_{X_i})|
$$
consider the rational map $\tilde{p}_i:|\sL|^{n-1} \to
V_i^{n-1}$. Since each $C\in S$ is a complete intersection
curve, the restriction $p_i:S\to V_i^{n-1}$ of $\tilde{p}_i$ is
a morphism. 

For $x \in X_i\cap \sigma(\sC)$ the condition ``$x
\in C_s$'' defines a multilinear subspace $\Delta_x^i$ of
$V_i^{n-1}$ of codimension $n-1$. 
Let $I_i\subset\{ 1, \ldots,d\}$ denote the set of all the $\nu$
with $x_\nu, x'_\nu \in X_i$. Then the codimension of
$$
\Delta^i = \Delta^i_{q_i} \cap \bigcap_{\nu \in
I_i}(\Delta^i_{x_\nu} \cap \Delta^i_{x'_\nu})
$$
is at most $(n-1)\cdot(2\cdot \# I_i +1)$.

Let $\sC^i \to \Delta^i$ be the intersection on $X_i$ of the
divisors in $\Delta^i \subset V_i^{n-1}$. Then the 
general fibre of $\sC^i \to X_i$ has dimension at least 
$$\dim\Delta^i+1-\dim X=\dim V_i^{n-1}+1-n-\codim\Delta^i$$
$$\geq
\dim (V_i^{n-1}) +1 -n - (n-1)\cdot(2\cdot \# I_i +1) \geq
(n-1)\cdot(2\cdot(\mu(\sL) - \# I_i)+r-1).
$$
Since some $C \in S$ contains all the points $x_j, x'_j$ and
$q_i$, the intersection
$$
Z = S\cap \bigcap_{i=1}^r p_i^{-1}(\Delta^i)
$$
is non-empty. For the restriction $\sC'$ of the universal curve
$\sC$ to $Z$ the dimension of the general fibre of $\sigma':\sC'
\to X$ over $X_i$ has dimension larger than or equal to
\begin{align*}
& (n-1)\cdot(2\cdot(\mu(\sL) - \# I_i)+r-1)-\sum_{j\neq i}
(n-1)\cdot(2\cdot \# I_j +1)\\
&= (n-1)\cdot 2\cdot(\mu(\sL) - \sum_{j=1}^r \# I_i)
=(n-1)\cdot 2\cdot(\mu(\sL) - d).
\end{align*}
By the inequality (\ref{bd}) the last expression is larger than or
equal to $0$ and $\sigma'$ is dominant.
\end{proof}
Let $G(C_z)$ denote the image of $\Pic^0(C_z)$ in $G$.
By the choice of $Z$ the intersection $B_z=C_z\cap X_\reg$ is
non-singular and the dimension of $G(C_z)=G(\Gamma_{B_z})$ is
the same as the dimension of $G(\Gamma_B)=G(\alpha_1,\ldots, \alpha_d)$.

By \ref{dominant2} the points $\alpha_i =
\phi(\sigma'(x_i)-\sigma'(x'_i))$ are contained in
$\Gamma_{B_z}$, hence 
$$
G(C_z)=G(\Gamma_{B_z}) = G(\Gamma_B)=G(C)
$$
for all $z \in Z$.

As $\sigma'$ is dominant, $\sigma'(\sC')$ contains some $V$,
open and dense in $X_{reg}$ (and hence in $X^{(n)}$). For $q\in V\cap X_i$
one finds some $z \in Z$ with $q \in C_z$. By \ref{dominant} $C_z$
contains the chosen base point $q_i$ and 
$$
\phi(\gamma(q_i) - \gamma(q)) \in G(C_z)=G(C).
$$
By \ref{equidim} (i), the points $\gamma(q_i) - \gamma(q)$ (for $q\in V$)
generate $CH^n(X)_{\deg 0}$. Since $\phi$ was assumed to be surjective,
we obtain $G=G(C)$, as claimed. 
\end{proof}
Using proposition \ref{bound} or its variant \ref{bound2}
the construction of $A^n(X)$
proceeds now along the lines of Lang's construction in 
\cite{La} of the Albanese variety of a smooth projective variety.
\begin{thm}\label{existence}
There exists a smooth connected commutative algebraic group $A^n(X)$
and a surjective regular homomorphism $\varphi: CH^n(X)_{\deg 0} \to
A^n(X)$  satisfying the following universal property:
For any smooth commutative algebraic group $G$ and for any
regular homomorphism $\phi:CH^n(X)_{\deg 0}\to G$ there
exists a unique homomorphism $h:A^n(X)\to G$ of algebraic groups
with $\phi=h\circ\varphi$. 

Moreover, if $k \subset K$ is an extension of algebraically
closed fields, then 
$$
A^n(X\times_kK)=A^n(X)\times_kK.
$$
\end{thm}
\begin{proof}
By lemma \ref{reg_sur} it is sufficient to consider connected
groups $G$, and surjective regular homomorphisms $\phi$.
 
By \ref{mu3} there exists a very ample invertible
sheaf $\sL$ which satisfies the assumption
\ref{dominant} and we can apply \ref{bound}.
(As we have seen in \ref{mu2} the inequality (\ref{eqmu2}) in
\ref{bound2} holds true for some $\sL$, and if $\Char(k)=0$
we can use the variant \ref{bound2}, as well.)

Let $g=\dim_k(H^1(C,\sO_C))$, for some curve $C\in |\sL|^{n-1}_0$ in
general position. Then for all regular homomorphisms
$\phi:CH^n(X)_{\deg 0}\to G$ to smooth connected commutative
algebraic groups $G$ the induced morphism
$$
\pi^{(-)}: S^{g+\mu(\sL)}(\Pi_{X_\reg})
\>\gamma^{(-)}_{g+\mu(\sL)}>> 
CH^n(X)_{\deg 0} \>\phi >> G 
$$
has a dense image in $G$. Hence for the product $\Pi$ of all the
different connected components of $S^{g+\mu(\sL)}(\Pi_{X_\reg})$
the induced morphism $\pi' : \Pi \to G$ is dominant and $\pi'$
induces a unique embedding of function fields $k(G) \subset k(\Pi)$.

If $\phi_\nu : CH^n(X)_{\deg 0} \to G_\nu$, for $\nu = 1,2$  are
two surjective regular homomorphisms to smooth connected commutative algebraic
groups, then 
$$
\phi_3:CH^n(X)_{\deg 0} \>>> G_1\times G_2
$$
is regular. Let $G_3$ be the image of $\phi_3$.
Then $\phi_\nu$ factors through the regular homomorphism
$\phi_3:CH^n(X)_{\deg 0} \to G_3$ and 
$k(G_\nu) \subset k(G_3) \subset k(\Pi)$, for $\nu = 1,2$.

Hence among the smooth connected commutative algebraic groups
$G$ with a regular surjective homomorphisms from
$\phi:CH^n(X)_{\deg 0}\to G$, there is one, $A^n(X)$, for which the subfield
$k(A^n(X))$ is maximal in $k(\Pi)$ and $A^n(X)$ dominates all
the other $G$ in a unique way.

It remains to show that $A^n(X)$ satisfies base-change for
algebraically closed fields. Let us write $Z_K = Z\times_kK$,
for a variety $Z$ defined over $k$. We first show:
\begin{claim}\label{basechange1} 
Let $K\supset k$ be an algebraically closed extension field of $k$.
The cycle map $\varkappa^{(-)}_K:(\Pi_{X_\reg})_K \to A^n(X)_K$
factors through a surjective homomorphism
$u_{K,k}:A^n(X_K) \to A^n(X)_K$ of algebraic groups.
\end{claim}
\begin{proof} Let $U=X_{reg}$, and let $(C', \iota)$ be an admissible pair
defined over $K$, with $B =\iota^{-1}(U_K)_{\reg}$.  
Choose a rational function $f \in R(C',X_K)$ such that 
$$ {\rm div} f = \sum a_i - \sum b_i $$ 
for $p =(a_1,b_1,\ldots,a_m,b_m)\in S^m(\Pi_B)(K)$.  Choose a smooth $k$
variety $S$ with  
$k(S)\subset K$, such that $C'$, $B$, $p$, $a_i$, $b_i$, $f$ come by
base-change from $k(S)$ to $K$ from 
$$\sC' \to S, \ \ \sB \to S, \ \ \pi: S \to S^m(\Pi_{\sB/S}), \ \ \alpha_i,
\beta_i: S \to \sB, \ \ \varphi \in k(\sB)^{\times} $$
with ${\rm div} \varphi = \sum \alpha_i - \sum \beta_i$. Since $f\in
R(C',X_K)$, we can replace $S$ by a dense open subscheme, so that we can
arrange that for each $s\in S(k)$, if we specialize to $\sC'_s = \sC _S
\times s$, then $\pi(s)$ maps to zero in $\Pic^0(\sC'_s)$. As 
\[S^m(\Pi_{\sB_s}) \>>> S^m (\Pi_U) \times s \>>> A^n(X) \times s\]
factors through $\Pic^0(\sC'_s)$, the composite morphism 
$$
S \> \pi >> S^m(\Pi_{\sB/S}) \>>> S^m(\Pi_U)\times S
\>>> A^n(X)\times S
$$
maps all $k$-points of $S$ to the zero section. Thus it is the zero
section, and therefore  
$S^m(\Pi_U)_K \to A^n(X)_K$
factors through $CH^n(X_K)$, inducing $u_{K,k}$ by lemma \ref{pic}.
\end{proof}

 Since $d_K:= \dim A^n(X_K)$ is bounded by $2n(g + \mu(\sL))$
(proposition~\ref{bound}), there is an algebraically closed field $K_1$
with $d_{K_1}=d_L$ for all $L\supset K_1$ algebraically closed. Since for
any ascending chain $K_i\subset K_{i+1}$ of algebraically closed fields
with $K_i \supset K_1$ one has 
$$
{\rm deg} \ u_{K_i,K_1} \leq {\rm deg} \ u_{K_{i+1},K_1} \leq
{\rm deg} \ u_{\cup K_i, K_1},
$$
one concludes that there is an algebraically closed field
$E\supset K_1$ such that $u_{E,L}$ is an isomorphism for all
algebraically closed fields $L\supset E$.
 
We will make use of the following lemma. 
\begin{lemma}\label{descent}
Let $K$ be a field, $W$, $Y$, $Z$ be geometrically integral $K$-varieties,
such that there are $K$-morphisms $\alpha:W\to Y$,
$\beta:W\to Z$, such that $\alpha$ has dense image. Then:
\begin{enumerate}
\item[(i)] there is at most one $K$-morphism $f:Y\to Z$ such that
$\beta=f\circ\alpha$
\item[(ii)] suppose that for some extension field $L$ of $K$, there
is an $L$-morphism $h:Y_L\to Z_L$ such that $\beta_L=h\circ
\gamma_L:W_L\to Z_L$;  then there is a $K$-morphism $f:Y\to Z$ as in
(i), and we have $h=f_L$. 
\end{enumerate}
\end{lemma}
\begin{proof} Let $\Gamma\subset W\times_K Z$ be the graph of $\beta$, and
let $\bar{\Gamma}\subset Y\times_K Z$ be the closure of $\alpha\times
1_Z(\Gamma)$. The projection $\bar{\Gamma}\to Y$ has dense image. If there
is a $K$-morphism $f:Y\to Z$ as in (i), then $\bar{\Gamma}$ must be its
graph, and so there is at most one such morphism, which exists precisely
when $\bar{\Gamma}\to Y$ is an isomorphism. Clearly if this is an 
isomorphism after base change to $L$, it is an isomorphism to begin with. 
\end{proof}

There is a smooth $k$ variety $S$, with $k(S) \subset E$, together
with a smooth commutative $S$-group scheme $\sA \to S$ with connected
fibers, such that $A^n(X_E) = \sA \times_S \Spec E$. Choosing $S$ small
enough one also has natural surjective $S$-morphisms $\Pi\times S \to \sA$
and $u_{S,k}:\sA \to A^n(X)\times S$, where $\Pi$ is the irreducible
variety constructed in the first part of the proof.

Let $F$  be an algebraic closure of the quotient field of $E\tensor_k E$,
$p: k(S\times_k S) \hookrightarrow F$ the natural inclusion, and let 
$p_i^*: k(S)\hookrightarrow k(S\times_k S)$, $i=1,2$  be the inclusions
defined by the two projections $p_i:S\times_k S\to S$. Set $q_i= p \circ
p_i^*$, and for any $S$-scheme $T$, let $q_i^*T$ be the $F$-scheme
obtained by the base change to $F$ determined by $q_i$.

The surjective $S$-morphism $\Pi\times S \to \sA$ gives rise  to the
surjections 
\[\alpha_i':q_i^* (\Pi\times_kS) = \Pi_{F} \>>> q_i^*\sA.\] 
By the assumption on $E$ the two $F$-varieties $q_i^*\sA$ are isomorphic 
via  
$$ u'= u_{F, E_1} \circ u^{-1}_{F,E_2}: q_1^*A^n(X_E) \>>> q_2^*A^n(X_E).
$$ 
where $E_i\subset F$ are the images of the two embeddings $E\into
E\tensor_k E\into F$, $x\mapsto x\tensor 1$ and $x\mapsto 1\tensor x$.
By construction, $u'$ satisfies $\alpha_2'=u'\circ\alpha_1'$. 

Hence by lemma~\ref{descent}, applied to the extension of fields
$k(S\times_kS)\into F$, the isomorphism $u'$ comes from an isomorphism 
$$u: (p_1^*\sA)_{k(S\times_k S)} \>>> (p_2^*\sA)_{k(S\times_k S)},$$
Then $u$ in fact extends uniquely to an isomorphism of groups schemes
(again denoted $u$)
\[u: (p_1^*\sA)_U \>>> (p_2^*\sA)_U\]
over an open dense subset $U \subset S\times_k S$.
Replacing $S$ by some open  dense subscheme, we may assume that $p_i:U\to S$
is surjective, for $i=1,2$. Further, if $\alpha_i:\Pi\times U\to
(p_i^*\sA)_U$, $i=1,2$ are the natural surjections, then $\alpha_2=u\circ
\alpha_1$. 

The uniqueness statement in lemma~\ref{descent}~(i) similarly implies
that $u$ satisfies the ``cocycle condition''
\[u_{23}\circ u_{12}=u_{13}:\pi_1^*\sA\>>>\pi_3^*\sA\]
on the fibers over the generic point of $S\times_k S\times_k S$, and hence
(by continuity) over the open dense subset  
\[\pi_{12}^{-1}(U)\cap\pi_{23}^{-1}(U)\cap \pi_{13}^{-1}(U)\subset
S\times_k S\times_k S.\]  
 Here $\pi_j:S\times_k S\times_k S\to S$ are the 3 projections, and
$u_{ij}=\pi_{ij}^*{u}$, for the 3 projections 
$\pi_{ij}:S\times_k S\times_k S\to S\times_k S$. 

Given two points $s_i \in S(k)$, one finds a third one $s 
\in S(k)$ such that $(s_1,s) \in U(k)$ and $(s, s_2) \in U(k)$. The cocycle 
condition implies that the induced composite isomorphism
$$ \theta_{s_1s_2}: \sA|_{s_1} \> u|_{(s_1,s)} >> \sA|_s \> u|_{(s,s_2)} 
>> \sA|_{s_2}$$
does not depend on the point $s\in S(k)$ chosen. Also
$u$ is compatible with the surjective morphisms $\Pi \times U \to
p_i^*\sA$. 

We claim that for each closed point $s\in S(k)$, the morphism $\Pi\times
s\to \sA|_s$ induces a regular homomorphism $CH^n(X)_{\deg 0}\to
\sA|_s$. Let $(C',\iota)$ be an admissible pair on $X$, defined over
$k$, with $B = \iota^{-1}(U)_{\reg}$. The morphism $(\Pi_B)_E\to
A^n(X_E)=\sA_E$, and the resulting morphism $S^g(\Pi_B)_E\to \sA_E$ (with
$g:=\dim\Pic^0(C')$) induces a homomorphism $\Pic^0(C')_E\to \sA_E$, since
by the defining property of $A^n(X_E)$, we have a factorization through
$CH^n(X_E)_{\deg 0)}$. Since $S^g(\Pi_B)\onto\Pic^0(C')$,
lemma~\ref{descent} gives a map $\Pic^0(C') \times_k k(S) \to
\sA_{k(S)}$, compatible with the maps from $(\Pi_B)_{k(S)}$. This then
induces a map $\Pic^0(C')\times S^0\to\sA_{S^0}$ for some open dense
subscheme $S^0\subset S$.  Choosing a $k$-point $s_1\in S^0(k)$, we get
that the map $\Pi_B\times s_1\to \sA|_{s_1}$ is compatible with a
homomorphism $\Pic^0(C')\to \sA|_{s_1}$. The isomorphism $u$ is 
compatible with the morphisms $\Pi_B \times U \to p_i^*\sA$.
Hence, the isomorphisms $\theta_{ss_1}$ are compatible with the maps  
$\Pi_B\cong \Pi_B\times s\to \sA|_s$ and $\Pi_B\cong \Pi_B\times s_1\to
\sA|_{s_1}$, for all $s\in S(k)$.  We deduce that for any $s\in S(k)$, 
the map $\Pi_B\times s\to \sA|_s$ gives rise to a compatible morphism
$\Pic^0(C') \times _k s \to \sA|_s$. This implies that there is an induced
regular homomorphism $CH^n(X)_{\deg 0}\to \sA|_s$ for each $s\in S(k)$.

Hence, one obtains morphisms $v_s: A^n(X) \to \sA|_s$, verifying $ v_t =
\theta_{st}      
\circ v_s$ for all $s,\ t \in S(k)$. Choosing now $s \in S(k)$,
we set $G= \sA|_s$, and $ v= v_s$.
The surjective morphism $\Pi\times_k S \to \sA$ induces a surjection from
$\Pi \times_k s$ onto $G$, hence $v$ is surjective. Since the composite
$$
u_{S,k}\circ v : A^n(X) \>>> G \>>> A^n(X)
$$
is an isomorphism, $v$ is an isomorphism. Thus
$u_{S, k}: \sA \to A^n(X)\times S$ is an isomorphism when restricted to
each $t \in S(k)$, and is hence an isomorphism. By base change to $E$ one
finds that $u_{E,k}:A^n(X_E) \to A^n(X)_E$ is an isomorphism.  

Now if $K \supset k$ is any algebraically closed field, we
choose an algebraically closed field $F$ with 
\begin{gather*}
\hspace*{1.6cm}F\supset K \supset k\mbox{ \ \ \ and \ \ \ }F\supset
E\supset k, \mbox{ \ \ \ hence}\\
u_{K, k} \otimes {\rm id}_F \circ u_{F,K}=
u_{F, k} = u_{E, k} \otimes {\rm id}_F \circ u_{F,E},
\end{gather*}
and $u_{K,k}$ is an isomorphism as well. 
\end{proof}

\section{Finite dimensional Chow groups of zero cycles}
The definition of finite dimensionality for the Chow group
of 0-cycles is a natural generalization of the definition in the
non-singular (and normal) case (see \cite{M}, \cite{S}).
\begin{defn}\label{fin-dim} $CH^n(X)$ is said to be {\em finite
dimensional} if for some $m>0$, the map 
$$
\gamma_m : S^m(X_\reg) \>>> CH^n(X)_{\deg 0}
$$
(introduced in \ref{difference}) is surjective.
\end{defn}
One can see that this is also equivalent to the statement that for some
integer $m'>0$, depending only on $X$, any element
of $CH^n(X)_{\deg 0}$ is represented by a 0-cycle $\sum_{i=1}^r\delta_i$,
where for each $i$, the cycle $\delta_i$ is a difference of two effective
0-cycles of degree $m'$ supported in $X_i$.  

In the proof of the next theorem we will use the notion of a
{\em regular map} $f:Z\to CH^n(X)_{\deg 0}$ from a variety $Z$.
This is a map of sets such that 
\begin{points}
\item the composition $Z\to CH^n(X)_{\deg 0}\to A^n(X)$ is a
morphism \vspace{.05cm}
\item there is a surjective morphism $W\to Z$ such that
$$W\>>> Z\>{f}>> CH^n(X)_{\deg 0}$$
factors as $W\>{h}>>S^m(X_\reg)\>{\gamma_m}>> CH^n(X)_{\deg 0},$
for some morphism $h$.  
\end{points}
For example, let $C'$ be a reduced Cartier curve in $X$ or, more general,
let $(C',\iota)$ be an admissible pair. Then the homomorphism
$\eta:\Pic^0(C')\to CH^n(X)_{\deg 0}$ 
constructed in lemma~\ref{gysin} is regular. In fact, the first
condition holds true by \ref{pic} whereas the second one
follows from the dominance of $S^g(C'_\reg) \to \Pic^0(C')$, 
for $g=\dim_k(H^1(C',\sO_{C'}))$.

Recall that $k$ is called a universal domain, if its
trancendence degree over the prime field is uncountable.
\begin{thm}\label{finite} Let $X$ be a projective variety of
dimension $n$ over a universal domain $k$. Then $CH^n(X)$ is
finite dimensional if and only if 
$$
\varphi:CH^n(X)_{\deg 0}\>>> A^n(X)
$$
defines an isomorphism between $CH^n(X)_{\deg 0}$ and
{\em (the closed points of)} $A^n(X)$. 
\end{thm}
\noindent{\bf Proof}\ \\ 
Let us write $U=X_\reg$. By lemma \ref{generators} (ii) the composite 
\[S^m(U)\>>> CH^n(X)_{\deg 0}\>>> A^n(X)\]
 is always surjective for $m=2\cdot \dim(A^n(X))$. Hence, if 
$CH^n(X)_{\deg 0}\to A^n(X)$ is an isomorphism, then $CH^n(X)$ is finite 
dimensional.

So the main thrust of the theorem is the converse, that if
$CH^n(X)_{\deg 0}$ is finite dimensional, then $CH^n(X)_{\deg 0}\to
A^n(X)$ is an isomorphism. We imitate Roitman's proof of this result in
the non-singular case, and the analogous argument for the normal case in
\cite{S}; however there are extra refinements needed here, particularly in
characteristic $p>0$. 

First, we note that by \cite{LW}, proposition 4.2, the ``graphs
of rational equivalence'' 
\[\Gamma_{r,s}=S^r(U)\times_{CH^n(X)_{\deg 0}}S^s(U)\] 
decompose as a countable union of locally closed subvarieties,
for each $r,s$, and over a universal domain such a decomposition
is unique. This immediately implies that if $f_j:Z_j\to CH^n(X)_{\deg 0}$,
$j=1,2$ are regular maps, then
\[Z_1\times_{CH^n(X)_{\deg 0}}Z_2 = \{(z_1,z_2)\in Z_1\times
Z_2; \ f_1(z_1)=f_2(z_2)\}\] 
is a countable union of locally closed subvarieties of $Z_1\times Z_2$. 

Now arguing as in \cite{S}, lemma~(1.3), we first see that if $G$ is a
smooth connected commutative algebraic group, and $f:G\to CH^n(X)_{\deg 0}$ is
any regular map which is a group homomorphism, then there is a
well-defined connected component of the identity $G^0\subset \ker f$,
which is a connected algebraic subgroup of $G$, and which has countable
index in $\ker f$. Then the induced homomorphism 
$$G/G^0\>>> CH^n(X)_{\deg 0}$$
has a countable kernel. Hence, for any such homomorphism $G\to
CH^n(X)_{\deg 0}$, we can define the {\em dimension of the image of $G$}
to be the dimension of $G/G^0$. 

Next, notice that if $G_1\to CH^n(X)_{\deg
0}$ and $G_2\to CH^n(X)_{\deg 0}$  are two regular homomorphisms from
smooth connected commutative algebraic groups $G_i$ such that
${\rm image}\,G_1$ is properly contained in ${\rm image}\, G_2$, then in
fact 
\[\dim\, {\rm image}\,G_1< \dim\, {\rm image}\,G_2.\]
Indeed, we may assume the maps $G_i\to CH^n(X)_{\deg 0}$ have countable
kernel, so that we wish to assert that $\dim G_1<\dim G_2$. Now
$G_3=G_1\times_{CH^n(X)_{\deg 0}}G_2$
is a subgroup of $G_1\times G_2$ which is a countable union of
locally closed subvarieties, and hence has a connected component of the
identity which is a connected algebraic group, say $H$. Then $H\to G_i$ are
homomorphisms of algebraic groups with countable, hence finite
kernels, such that $H\to G_1$ is surjective, and the image of $H$ in $G_2$
is a strictly smaller subgroup. Thus $\dim G_1=\dim H<\dim G_2$.

Now suppose $\gamma_m$ is surjective. We claim that for any homomorphism
$$
G\>>> CH^n(X)_{\deg 0}
$$
as above, with countable kernel, we have 
$\dim G\leq \dim S^m(U).$
Indeed, 
\[G\times_{CH^n(X)_{\deg 0}}S^m(U)\] 
is a countable union of subvarieties of $G\times S^m(U)$ which projects onto
$G$, and maps to $S^m(U)$ with countable fibres. Hence some irreducible
component of this fibre product dominates $G$ under the projection, and
maps to $S^m(U)$ with finite fibres. 

We now claim that we can find a finite number of reduced complete
intersection curves $C_1,\ldots,C_s$ such that the induced homomorphism
from $\oplus \Pic^0(C_j)$ to $CH^n(X)_{\deg 0}$ is surjective. Indeed, given a
finite collection of such curves, if 
$$P=\bigoplus \Pic^0(C_j)\>>> CH^n(X)_{\deg 0}$$
is not surjective, we can find a curve $C$ of the same sort such that
$$\image( \Pic^0(C)\>>> CH^n(X)_{\deg 0})$$
is not contained in the image of $P$. Then the induced map
$$P\times \Pic^0(C)\>>> CH^n(X)_{\deg 0}$$
has strictly larger
dimensional image than that of $P$.  Since the dimension 
of the image is bounded above by $\dim S^m(U)=mn$, this process can be
repeated at most a finite number of times.

So we may assume given a surjective regular homomorphism 
$$f: A\>>> CH^n(X)_{\deg 0}$$
with countable kernel, where $A$ is a
connected smooth commutative algebraic group, and for some
Cartier curves $C_1, \ldots , C_s$ a surjective homomorphism
\begin{equation}\label{surjection}
\bigoplus_{j=1}^s \Pic^0(C_j) \> \oplus \rho_j >> A.
\end{equation}
Note that the composition $h:A\to CH^n(X)_{\deg 0}\to A^n(X)$ is
then a surjective homomorphism of algebraic groups.

We now distinguish between the case $k=\C$, and that of a general
universal domain $k$.

\begin{proof}[Proof of \ref{finite} for $k=\C$]

We first show that the surjective homomorphism $h:A\onto A^n(X)$ is an
{\em isogeny}.  Clearly $h$ induces an  injective homomorphism
$$h^*:\Omega(A^n(X))\to\Omega(A).$$  
We will use proposition~\ref{basic} to show that
$h^*:\Omega(A^n(X))\to\Omega(A)$ is an isomorphism.  Since $h^*$ is
injective, it suffices to prove that $\dim \Omega(A)\leq\dim
\Omega(A^n(X))$. 

Consider the set $\Gamma=U\times_{CH^n(X)_{\deg 0}}A$.
This is a countable union of algebraic subvarieties, and maps surjectively
to $U$ under the projection. Recalling that $U=\cup_j U_j$, we can then
find irreducible varieties $\Gamma_j\subset \Gamma$ such that
$\Gamma_j$ dominates $U_j$ under the projection $\Gamma\to U$. Then
$\pi_j:\Gamma_j\to U_j$ has countable, and hence finite, fibres. 
Let $d_j$ be the degree of $\pi_j$, and let $V_j\subset U_j$ be a dense
open subset such that $\pi_j:\pi_j^{-1}(V_j)\to V_j$ is an \'etale
covering of degree $d_j$. Let $c$ be the l.c.m. of the $d_j$, and let
$c=d_jc_j$. If $q:\Gamma\to A$ is the second projection, then consider the
morphism 
\begin{gather*}
\mu:V=\bigcup_jV_j\>>> A,\\
\mu(x)=c_j\mathop{\sum}_{y\in \pi_j^{-1}(x)}q(y)\,\,\mbox{for $x\in
V_j$.}
\end{gather*} 
One verifies at once that the diagram
\begin{equation}\label{diag-1}
\begin{CD}
\bigcup_j V_j=V \>{\mu}>> A \\
\V VV \V V f V\\
\bigcup_j U_j=U \>{c\cdot\gamma_1}>>CH^n(X)_{\deg 0}
\end{CD}
\end{equation}
commutes. 

The image of $\mu(V)$ in $CH^n(X)_{\deg 0}$ generates
$CH^n(X)_{\deg 0}$ as a group, since any 0-cycle on $X$ is rationally
equivalent to a cycle supported on $V$. Hence the subgroup of $A$
generated by $\mu(V)$ has countable index, and is also a countable
increasing union of constructible subsets, namely the images of $\mu(V)^{2m}$
under the maps
\begin{align*}
\sigma_m:A^{2m}&\>>>A,\;\;\;\; m\geq 1,\\ 
(a_1,\ldots,a_{2m})&\longmapsto
a_1+\cdots+a_m-a_{m+1}-\cdots-a_{2m}. 
\end{align*}
By dimension considerations, one of the subsets
$\sigma_m(\mu(V)^{2m})$ must be dense in $A$, and then
$\sigma_{2m}(\mu(V)^{4m})=A$. Hence the induced map on 1-forms
$$
\Omega(A)\to H^0(V^{4m},\Omega^1_{V^{4m}/\C})
$$
is injective.  Now the action of $\sigma_{2m}$ on 1-forms is given by
\[\sigma_{2m}^*(\omega)=
(\omega,\ldots,\omega,-\omega,-\omega,\cdots,-\omega).\]
This means that the map on 1-forms
$\Omega(A)\to H^0(V,\Omega^1_{V/\C})$, induced by the morphism $V\to
\mu(V)\into A$, is injective. 
 
We claim that $\image\Omega(A)\subset \Omega(A^n(X))$, so that
$\dim\Omega(A)\leq \dim\Omega(A^n(X))$. To see this, it
suffices by proposition~\ref{basic} to show that for any reduced Cartier
curve $C\subset X$ with $B=(C_{\reg})\cap V$ dense in $C$, the image of
any element of $\Omega(A)$ in $H^0(B,\Omega^1_{B/\C})$ lies in the image of
$H^0(C,\omega_C)$. Fixing base points in each component of $B$,
we obtain a morphism $\vartheta:C_{\reg}\to\Pic^0(C)$. 
If $C_i$ is any component of $C_{\reg}$, then the two induced maps 
\begin{gather*}
C_i\>>> \Pic^0(C)\>>> CH^n(X)_{\deg 0},\\
C_i\into U\>{\gamma_1}>> CH^n(X)_{\deg 0}
\end{gather*}
agree up to translation by a fixed element of $CH^n(X)_{\deg
0}$. 

Now consider the subgroup $\Gamma_C=\Pic^0(C)\times_{CH^n(X)_{\deg 0}} A$.
As before, this is a countable union of subvarieties of $\Pic^0(C)\times
A$. Hence there is a connected algebraic subgroup
$\Gamma^0_C\subset\Gamma_C$ such that $\Gamma_C/\Gamma^0_C$ is a countable
group. Further, $\Gamma_C\to\Pic^0(C)$ is surjective with countable
fibres. Hence $\Gamma^0_C\to\Pic^0(C)$ is an isogeny.
Restricting (\ref{diag-1}) one obtains a
commutative diagram
\[\begin{CD}
B \>{\mu}>> A\\
\V VV \V V f V\\
C_{\reg}\>{c\gamma_1}>> CH^n(X)_{\deg 0}
\end{CD}\]
and hence a morphism $B\to \Gamma^0_C$ such that 
\begin{points}
\item for each component $C_i$ of $C_{\reg}$, 
the composite $C_i\cap B\to \Gamma^0_C\to\Pic^0(C)$
equals the restriction of the composite $C_i\to \Pic^0(C)\longby{c\cdot
}\Pic^0(C)$, up to a translation (here $c\cdot$ denotes multiplication
by $c$)
\item $C_i\cap B\to \Gamma^0_C\to A$ agrees with $\mu$, up to a
translation. 
\end{points}
Hence, by (ii), 
$$
{\rm image}(\Omega(A)\>{\mu^*}>>\Gamma(B,\Omega^1_{C/\C}))
\subset {\rm image}\left(\Omega(\Gamma^0_C)\>>>
\Gamma(B,\Omega^1_{C/\C})\right) 
$$
while by (i), 
\[{\rm image}\,\Omega(\Gamma^0_C)={\rm
image}\,\Omega(\Pic^0(C))={\rm image}\,\Gamma(C,\omega_C).\]
 Since $C$ was arbitrary, we have verified the hypotheses of
proposition~\ref{basic}.  This completes the proof that the composite
$h:A\to CH^n(X)_{\deg 0}\to A^n(X)$ is an isogeny. 

In particular, $f:A\onto CH^n(X)_{\deg 0}$ has a finite kernel. Replacing
$A$ by $A/(\ker f)$, we may assume given a regular homomorphism $f:A\to
CH^n(X)_{\deg 0}$ which is an {\em isomorphism} of groups. Now repeating
the above arguments once more, we obtain (\ref{diag-1}) with $c=1$. By
corollary~\ref{equ-reg}, this means the group isomorphism
$f^{-1}:CH^n(X)_{\deg 0}\to A$ is a regular homomorphism, which must
factor through $\varphi:CH^n(X)_{\deg 0}\onto A^n(X)$. This forces
$\varphi$ to be an isomorphism of groups, as well. 
\end{proof}
\begin{rmk}\label{roitman}
Over the field of complex numbers the last part of the proof of \ref{finite}
is consistent with the Roitman theorem proved in \cite{BiS}.
In fact, if 
$$A\cong CH^n(X)_{\deg 0}\>>> A^n(X)$$
is surjective with finite kernel the generalization of Roitman's
theorem implies that the composite 
\[A\>{\cong}>> CH^n(X)_{\deg 0}\>>> A^n(X)\>>> J^n(X)\]
is an isomorphism on torsion subgroups, so that $CH^n(X)_{\deg 0}\to
A^n(X)$ is an injection on torsion subgroups. Hence the isogeny
$A\cong CH^n(X)_{\deg 0}\to A^n(X)$ must be an isomorphism. 
\end{rmk}

In the algebraic case we have to modify the arguments, in
particular since the lower horizontal morphism in the diagram
(\ref{diag-1}) need not to be surjective in characteristic $p>0$.

\begin{proof}[Proof of \ref{finite} for $k$ a universal
domain] \ \\
Let us write $B$ for the kernel of $h:A \to A^n(X)$, a
closed subgroup scheme of $A$, not necessarily reduced.
We may replace $A$ by $A/\kappa$, for any (zero dimensional) closed
subgroup scheme $\kappa$ of $B$ such that $\kappa(k)\subset\ker f$.

The group $B$ acts on $U\times_{A^n(X)}A$ with quotient
$U\times_{A^n(X)}A^n(X)=U$. The kernel $\sK$ of the map
$A(k) \to CH^n(X)_{\deg 0}$ consists of countably many closed
points, the induced action on $U\times_{CH^n(X)_{\deg 0}}A$
is free, and the induced map on the quotient 
$$
(U\times_{CH^n(X)_{\deg 0}}A)/\sK = U\times_{CH^n(X)_{\deg
0}}(A/\sK) \>>> U
$$
is a bijection on the closed points.

Let $V \subset U$ be an open dense subscheme, and 
let $\Gamma_j$ be a locally closed irreducible subscheme of
$V\times_{A^n(X)}A$, contained in $V\times_{CH^n(X)_{\deg 0}}A$,
and dominant over the component $V_j=V\cap U_j$ of $V$ under the first
projection. For $V$ small enough, we may assume that $\Gamma_j
\to V_j$ is finite. Let $\kappa_j\subset \sK$ be the subgroup
of elements $g$ with $g(\Gamma_j) = \Gamma_j$. Then
$\kappa_j$ is a finite group and $\Gamma_j/\kappa_j \to V_j$
is an isomorphism on the closed points. Replacing
$\Gamma_j$ by its image in $U\times_{A^n(X)}(A/\kappa_j)$ and
$A$ by $A/\kappa_j$ we may assume that $\kappa_j$ is trivial,
and thereby that $\Gamma_j \to V_j$ is purely inseparable.

Repeating this construction for the different components of $U$ we finally
reduce to the situation, where $U$ has an open dense subscheme $V$, and
where $V\times_{CH^n(X)_{\deg 0}}A$ has a closed subscheme $\Gamma$ which
is purely inseparable over $V$. 

Assume that $\Gamma \to V$ is not an isomorphism, in particular,
that the characteristic of $k$ is $p>0$. The restriction of the
group action to $B\times \Gamma$ factors as 
$$
B \times \Gamma \> \cong >> (V\times_{A^n(X)}A)\times_V\Gamma
\> pr_1 >> V\times_{A^n(X)}A
$$
and the preimage $S(\Gamma)$ of $\Gamma \subset V\times_{A^n(X)}A$
is isomorphic to $\Gamma \times_V \Gamma$. Thus $S(\Gamma)$ is a
subscheme of $B\times \Gamma$, supported in the zero section
$\{e\}\times \Gamma$. Hence $S(\Gamma)$ is contained in the
$\nu$-th infinitesimal neighbourhood $\{e\}_\nu \times \Gamma$
of the zero section, for some $\nu > 0$.

The kernel $\kappa^{(\nu')}$ of the $\nu'$-th geometric
Frobenius $F^{(\nu')}: B \to B^{(\nu')}$ is defined by the
sheaf of ideals in $\sO_B$, generated by the $p^{\nu'}$-th
powers of the generators of the sheaf of ideals $\rm \bf m$
defining $\{e\} \subset B$.
For some $\nu'>0$ it is contained in ${\rm \bf m}^\nu$
and $\{e\}_\nu$ is a subscheme of $\kappa^{(\nu')}$.

Dividing $A$ by $\kappa^{(\nu')}$, we may assume that
$S(\Gamma)=\Gamma\times_V\Gamma$ is isomorphic to $\Gamma$, and
thereby that $\Gamma$ is isomorphic to $V$.

Independent of the characteristic of $k$, we have thus reduced 
to the situation where $U$ has an open dense
subscheme $V$, for which 
$$
pr_1:V\times_{CH^n(X)_{\deg 0}}A \to V
$$
has a section, such that on projecting to $A$ we obtain a morphism 
$\mu:V\to A$ and (using the notation introduced in \ref{difference})
a commutative diagram
\begin{equation}\label{diag2}
\begin{CD}
\Pi_V \>{\mu^{(-)}}>> A \\
\V\subseteq VV \V V f V\\
\Pi_U \>{\gamma^{(-)}}>>CH^n(X)_{\deg 0}.
\end{CD}
\end{equation}
\begin{claim}\label{claim2}
There exists a surjective homomorphism $\phi:CH^n(X)_{\deg 0} \to A$
with $\mu^{(-)}=\phi \circ \gamma^{(-)}|_{\Pi_V}$. In
particular, $\phi$ is regular.
\end{claim}
\begin{proof}
Let $(C',\iota)$ be an admissible pair with $B=(\iota^{-1}(V))_\reg$
dense in $C'$. By restriction (\ref{diag2}) gives rise to a
commutative diagram
\begin{equation}\label{diag3}
\begin{CD}
\Pi_B \>{{\mu'}^{(-)}}>> A \\
\V= VV \V V f V\\
\Pi_B \>{\gamma_B^{(-)}}>>CH^n(X)_{\deg 0}
\end{CD}
\end{equation}
where $\mu'=\mu|_B$. By lemma \ref{pic} the lower horizontal
map in the diagram (\ref{diag3}) factors as
\begin{equation}\label{diag4}
\begin{CDS}
\Pi_B \>\gamma_B^{(-)}>> CH^n(X)_{\deg 0}\\
\novarr \SE E \vartheta^{(-)} E \A \eta A A \\
\noharr \Pic^0(C').
\end{CDS}
\end{equation}
Let $\Gamma_{C'}^0$ be the connected component of
$\Gamma_{C'}=\Pic^0({C'})\times_{CH^n(X)_{\deg 0}} A$ containing the origin.
$\Gamma_{C'}/\Gamma_{C'}^0$ is a countable group and $\Gamma_{C'}^0 \to
\Pic^0({C'})$ is an isogeny. Since the diagrams (\ref{diag3}) and
(\ref{diag4}) are commutative, the image of
$$
\Pi_B \> (\vartheta^{(-)}, {\mu'}^{(-)}) >> \Pic^0({C'})\times A
$$
is contained in $\Gamma_{C'}^0$. This implies that $\Gamma_{C'}^0\to
\Pic^0(C')$ must be an isomorphism. In fact, by \ref{generators2} there is
an open connected subscheme $W$ of $S^g(B)$ such that the morphism 
$\vartheta_{W} : W \to \Pic^0(C)$ is an open embedding.
On the other hand, $\vartheta_W$ factors through the isogeny
$\Gamma_{C'}^0\to \Pic^0(C')$.

Hence the morphism ${\mu'}^{(-)}$ in the diagram (\ref{diag3})
is the composite
$$
\Pi_B \> \vartheta^{(-)} >> \Pic^0(C')\cong \Gamma_{C'}^0 \> pr_2
>> A,
$$
and the condition (b) in lemma \ref{pic} holds true. 

Thereby the homomorphism $\phi$ in \ref{claim2} exists, and 
it remains to show that $\phi$ is surjective.
Equivalently, it suffices to show that the image of $\phi$ generates $A$
as a group, which will follow if we show that $\mu^{-}(\Pi_V)$ generates
$A$. But we know that $\gamma^{-}(\Pi_V)$ generates $CH^n(X)_{\deg 0}$,
and so $\mu^{-}(\Pi_V)$ generates a subgroup of countable index in $A$.
Since $k$ is a universal domain, $\mu^{-}(\Pi_V)$ generates $A$. 
\end{proof}

By claim \ref{claim2} and by the universal property for
$A^n(X)$ the regular homomorphism $\phi: CH^n(X)_{\deg 0}\to A$
factors through a homomorphism of algebraic groups $\chi:A^n(X)\to A$.
Since $\phi$ is surjective, the induced morphism $\chi$ is surjective as
well. Further, the composite $CH^n(X)_{\deg 0}\>{\phi}>> A\> f>>
CH^n(X)_{\deg 0}$ is clearly the identity, since it is so on the image of
$\Pi_V$, which is a set of generators. By the universal property of 
$\varphi:CH^n(X)_{\deg 0}\to A^n(X)$, we deduce that the composite
$A^n(X)\>{\chi}>> A\>h>> A^n(X)$ is the identity. Hence $\chi$ and $h$ are
inverse isomorphisms of algebraic groups, and $f:A\to CH^n(X)_{\deg 0}$
and $\varphi:CH^n(X)_{\deg 0}\to A^n(X)$ are both isomorphisms (of groups)
as well.  
\end{proof}
\bibliographystyle{plain}

\begin{thebibliography}{BPW}
\bibitem[Alb]{Albanese} Giacomo Albanese, {\it Corrispondenze algebriche
fra i punti di une superficie algebriche}, Ann. Scuola Norm. Pisa (II) 3
(1934) 1-26 and 149-182.
\bibitem[BPW]{BPW} L. Barbieri-Viale, C. Pedrini, C. A. Weibel,
{\it Roitman's theorem for singular complex projective surfaces},
Duke Math. J.  {\bf 84} (1996) 155-190.
\bibitem[BS]{BS} L. Barbieri-Viale, V. Srinivas,
{\it The Albanese 1-motive}, in preparation.
\bibitem[BiS]{BiS} J. Biswas, V. Srinivas, {\it Roitman's
theorem for singular projective varieties}, preprint.
\bibitem[BH]{BH} T. Bloom, M. Herrera, {\it De Rham cohomology
of an analytic space}, Invent. Math. {\bf 7} (1969) 275-296.
\bibitem[BLR]{BLR} S. Bosch, W. L\"utkebohmert, M. Raynaud,
{\it N\'eron Models}, Ergebnisse der Mathematik (3. Folge) {\bf
21}, Springer-Verlag 1990.
\bibitem[C]{C} A. Collino, {\it Torsion in the Chow group of 
codimension 2: the case of varieties with isolated singularities},
J. Pure Appl. Alg. {\bf 34} (1984) 147-153.
\bibitem[D]{D} P. Deligne, {\it Th\'eorie de Hodge, II, III},
Publ. Math. IHES {\bf 40} (1971) 5-57 and {\bf 44} (1974) 5-77.
\bibitem[EV]{EV} H. Esnault, E. Viehweg, {\it Deligne-Beilinson
cohomology}, in {\it Beilinson's Conjectures on special values of
L-functions}, Perspectives in Math. {\bf 4}, Academic Press, New
York 1988.
\bibitem[FW]{FW} G. Faltings, G. W\"ustholz, {\it Einbettungen
kommutativer algebraischer Gruppen und einige ihrer
Eigenschaften}, J. Reine Angew. Math. {\bf 354} (1984) 175-205. 
\bibitem[F]{F} H. Flenner, {\it Die S\"atze von Bertini f\"ur
lokale Ringe}, Math. Ann. {\bf 229} (1977) 97-111.
\bibitem[Ful]{Fulton} W. Fulton, {\it Intersection Theory}, Ergeb. Math.
3. Folge, Band 2, Springer-Verlag 1984.
\bibitem[Ha]{Ha} R. Hartshorne, {\it Algebraic Geometry}, Grad.
Texts in Math. {\bf 52}, Springer-Verlag 1977.
\bibitem[Ha2]{Ha2} R. Hartshorne {\it Ample subvarieties of
Algebraic Varieties}, Lect. Notes in Math. {\bf 156},
Springer-Verlag 1970. 
\bibitem[I]{Igusa} J.-I. Igusa, {\it On the Picard varieties attached to
algebraic varieties}, Amer. J. Math. 74 (1952) 1-22.
\bibitem[La]{La} S. Lang, {\it Abelian Varieties}, Interscience
Publisher, New York 1959
\bibitem[LW]{LW} M. Levine, C. A. Weibel, {\it Zero cycles and
complete intersections on singular varieties}, J. Reine Ang.
Math. {\bf 359} (1985) 106-120. 
\bibitem[L]{L} M. Levine, {\it Deligne-Beilinson cohomology for
singular varieties}, in {Algebraic K-theory, commutative algebra
and algebraic geometry (Santa Margherita Liguria, 1989)},
Contemp. Math. {\bf 126}, A.M.S., Providence (1992) 113-145. 
\bibitem[L2]{L2} M. Levine, {\it A geometric theory of the Chow
ring for singular varieties}, preprint (1983). 
\bibitem[Ma]{Matsusaka} T. Matsusaka, {\it On the algebraic construction
of the Picard variety II}, Jap. J. Math. 22 (1952), 51-62 (1953).
\bibitem[M]{M} D. Mumford, {\it Rational equivalence of 0-cycles
on algebraic surfaces}, J. Math. Kyoto Univ. {\bf 9} (1968)
195-240. 
\bibitem[R]{R} A. A. Roitman, {\it Rational equivalence of
0-cycles}, Math. USSR Sbornik {\bf 18} (1972) 571-588.
\bibitem[R2]{R2} A. A. Roitman, {\it The torsion of the group of
0-cycles modulo rational equivalence}, Ann. Math. {\bf 111}
(1980) 553-569,
\bibitem[Se]{Se} J. P. Serre, {\it Morphismes universels
et diff\'erentielles de troisi\`eme esp\`ece}, S\'em. C.
Chevalley (1958/59) 10,
\bibitem[S]{S} V. Srinivas, {\it Zero cycles on a singular
surface II}, J. Reine Ang. Math. {\bf 362} (1985) 4-27.
\bibitem[Si]{Si} Y.-T. Siu, {\it Techniques of extension of
analytic objects}, Lect. Notes in Pure and Appl. Math., {\bf 8},
Marcel-Dekker, New York (1975).
\bibitem[W]{Weil} A. Weil, {\it Collected Papers Vol.~I (1926-1951)},
Corrected Second Printing, Springer-Verlag, 1980.
\end{thebibliography}

\end{document}